\documentclass[11pt,reqno,twoside,makeidx]{amsart}
\usepackage[margin=2.7cm]{geometry}

\usepackage{caption}
\usepackage{graphicx}
\usepackage{dsfont}
\usepackage{amsthm}
\usepackage{amsmath}
\usepackage{amsfonts}
\usepackage{amscd}
\usepackage{fancyhdr}
\usepackage{enumerate}
\usepackage{enumitem}
\usepackage{url}
\usepackage{lipsum}
\usepackage{hyperref} 

\usepackage{tikz}
\usepackage{tikz-cd}
    \usetikzlibrary{arrows}
    \usetikzlibrary{matrix}
    \usepackage{pictexwd,dcpic}
    \usepackage{sidecap}
    \graphicspath{./imagine/}

\usepackage{caption}
\usepackage{subcaption}

\usepackage{cancel}
\usepackage{import}

\usepackage{amssymb}

\usepackage{mathtools}
\usepackage{layout}

\usepackage{array}
\newcolumntype{C}[1]{>{\centering\let\newline\\\arraybackslash\hspace{0pt}}m{#1}} 

\usepackage{romannum} 
\AtBeginDocument{\pagenumbering{arabic}}

\usepackage{colortbl,xcolor} 

\newcommand{\fund}[1]{ \pi_1 \left( #1 \right)}

\newcommand\restr[2]{{
  \left.\kern-\nulldelimiterspace 
  #1 
  \vphantom{|} 
  \right|_{#2} 
  }}
  
 \newcommand{\ZSU}{\mathcal{Z}\left(SU(2)\right)}

\newcommand{\Ui}{U(1)}
\newcommand{\normalsubgroup}[1]{\langle\! \langle #1 \rangle \! \rangle}
\newcommand{\R}[1]{\mathcal{R}_{U(1)} \left( #1\right)}
\newcommand{\geo}[2]{\Delta \left( #1 ,#2\right)}

\tikzset{cross/.style={cross out, draw=black, minimum size=2*(#1-\pgflinewidth), inner sep=0pt, outer sep=0pt},
cross/.default={1pt}}

\DeclareMathSymbol{\shortminus}{\mathbin}{AMSa}{"39}

\DeclareMathOperator{\ima}{Im}
\DeclareMathOperator{\Tr}{Tr}
\DeclareMathOperator\coker{coker}
\DeclareMathOperator{\rank}{rank}

\usepackage{mathrsfs}

\usepackage{multicol}
\usepackage{lscape}
\usepackage{cjhebrew}

\usepackage{bbold}

\usepackage{graphicx,nicefrac}
 \usepackage{framed}

\usepackage{algorithm}
\usepackage{algpseudocode}
\algtext*{EndFor}
\algtext*{EndIf}
\usepackage{tikz}
\makeatletter
\newtheorem*{rep@teo}{\rep@title}
\newcommand{\newreptheorem}[2]{%
\newenvironment{rep#1}[1]{%
 \def\rep@title{#2 \ref{##1}}%
 \begin{rep@teo}}%
 {\end{rep@teo}}}
 \newreptheorem{teo}{Theorem}
 \newreptheorem{cor}{Corollary}
 \newreptheorem{prop}{Proposition}
 \newreptheorem{defn}{Definition}
 \newreptheorem{lemma}{Lemma}
 \newreptheorem{condition}{Condition}

\theoremstyle{definition}
\newtheorem{defn}{Definition}[section]

\newtheorem{exmp}[defn]{Example}
\newtheorem{rmk}[defn]{Remark}
\newtheorem{notation}[defn]{Notation}

\newtheorem*{PN*}{Please Note}
\newtheorem{claim}{Claim}

\theoremstyle{plain}
\newtheorem{teo}[defn]{Theorem}
\newtheorem*{teo*}{Theorem}
\newtheorem{cor}[defn]{Corollary}
\newtheorem{lemma}[defn]{Lemma}
\newtheorem{prop}[defn]{Proposition}

\newtheorem*{prob*}{Problem}

\newtheorem*{qst*}{Question}
\newtheorem{fact}[defn]{Fact}

\newtheorem{conj}[defn]{Conjecture}

\numberwithin{equation}{section}
\makeatother

\usepackage[
backend=biber,
style=alphabetic,
sorting=ynt,
natbib=true,
sorting=anyvt
]{biblatex}

\bibliography{biblio_new.bib}

\usepackage{overpic}

\author{Giacomo Bascape}
\address{Département de Mathématiques, Université du Québec à Montréal, Canada}
\email{giacomo.bascape@icloud.com}
\title[Toroidal graph manifolds with small homology are not SU(2)-abelian]{Toroidal graph manifolds with small homology are not SU(2)-abelian}

\begin{document}
\begin{abstract}
 We show that if $Y$ is a toroidal closed graph manifold rational homology $3$-sphere with $|H_1(Y;\mathbb{Z})| \le 5$, then there exists an irreducible representation $\fund{Y} \to SU(2)$, using topological methods and avoiding the use of gauge theory. 
This answers positively to a conjecture by Baldwin and Sivek \cite[Conjecture 1.5]{InstantoLSpaceAndSplicing} in the case of graph manifolds.
\end{abstract}

\maketitle

\section{Introduction}

In this work we will consider only compact and oriented $3$-manifolds, with or without boundary.
The study of $SU(2)$-representations of $3$-manifold fundamental groups is a meeting point between topology, geometry, and gauge theory. 
A particularly illuminating example is the resolution of Property P by Kronheimer and Mrowka in \cite{KronheimerMrowkaDehnSurgeryFundamental}. In particular, they show that if $Y$ is a $3$-manifold obtained by performing a non-trivial surgery along a non-trivial knot in $S^3$, then there exists an irreducible $SU(2)$-representation $\fund{Y} \to SU(2)$. The existence of this representation implies that $\fund{Y}$ is non-trivial and therefore that $Y$ is not $S^3$, and this proves that if $S^3$ is obtained by performing a surgery along a knot, then either the surgery or the knot is trivial.
We recall that an $SU(2)$-representation is non-abelian if and only if it is irreducible.

A $3$-manifold $Y$ is called \emph{$SU(2)$-abelian} if every representation $\fund{Y}\to SU(2)$ has abelian image.
The prototypical examples of $SU(2)$-abelian manifolds are lens spaces, these are the closed $3$-manifolds whose fundamental group is cyclic. A less elementary example is contained in \cite{Cornwell2015CharacterVO}, in which it is proved that the double cover of $S^3$ branched along the knot $8_{18}$ is a closed (hyperbolic) $SU(2)$-abelian manifold. 
Understanding the $SU(2)$-abelian status of a $3$-manifold is a remarkable problem in low-dimensional topology with a deep connection to instanton Floer homology, a $3$-manifold invariant that determines the diffeomorphism class up to finitely many possibilities. As an example, in \cite[Theorem 4.6]{SivekBalwinSteinFilling} Baldwin and Sivek showed that if a $3$-manifold is $SU(2)$-abelian, then it has trivial instanton Floer homology whenever the fundamental group has the property of being cyclically finite, whose definition can be found in \cite{BoyerNicasAndrewVarietiesOfGroupsRepresentations}.

In this work we are interested in manifolds with the opposite behavior: we study manifolds that are not $SU(2)$-abelian, and therefore admitting a non-abelian $SU(2)$-representation.
Our main motivation comes from the following conjecture:
\begin{conj}[{\cite[Conjecture 1.5]{InstantoLSpaceAndSplicing}}]\label{conj: BS minore o uguale a 4}
    Let $Y$ be a toroidal rational homology $3$-sphere. If $|H_1(Y;\mathbb{Z})| \le 4$, then $Y$ is not $SU(2)$-abelian.
\end{conj}
It is worth mentioning that if $Y$ is a toroidal $3$-manifold such that $|H_1(Y;\mathbb{Z})|=1$, then there exists a non-abelian representation $\fund{Y} \to SU(2)$ as a consequence of \cite[Theorem 1.1]{ToroidalZS3HaveIrreducibleRepres}.
As an original result, we settle Conjecture \ref{conj: BS minore o uguale a 4} for toroidal graph manifolds.  
\begin{teo}\label{teo: sigma small then not SU(2)-abelian}
    Let $Y$ be a toroidal graph manifold rational homology $3$-sphere. If $|H_1(Y;\mathbb{Z})| \le 5$, then $Y$ is not $SU(2)$-abelian.
\end{teo}
Notice that in Theorem \ref{teo: sigma small then not SU(2)-abelian} the hypothesis that the manifold is toroidal cannot be weakened: for example, the connected sum $\mathbb{RP}^3 \# \mathbb{RP}^3$ is an $SU(2)$-abelian graph manifold with first homology of cardinality $4$.
A similar result to Theorem \ref{teo: sigma small then not SU(2)-abelian} can be found in \cite[Theorem 1.1]{baldwin2024smallheegaardgenussu2} in which it is shown that if $Y$ is a rational homology $3$-sphere with Heegaard genus $2$ and $|H_1(Y;\mathbb{Z})| \in \{1,3,5\}$, then $Y$ is not $SU(2)$-abelian. 

Conjecture \ref{conj: BS minore o uguale a 4} has a natural connection in the realm of Heegaard Floer homology.
In this work we shall not rely on Heegaard Floer theory, and therefore we only provide the minimal background.
Recall that a rational homology $3$-sphere $Y$ is said to be an \emph{Heegaard Floer L-space} if its (hat) Heegaard Floer group is the simplest possible, namely if $\rank \widehat{HF}(Y)=|H_1(Y;\mathbb{Z})|$. Details can be found in \cite{HFdef2} and \cite{HFdef1}. The notions of being $SU(2)$-abelian and of being a Heegaard Floer L-space are expected to be linked in the following way:
\begin{conj}[\cite{zhang2019remarks}]\label{conj: SU(2)-ab implies HF L-space}
    Let $Y$ be a rational homology $3$-sphere. If $Y$ is $SU(2)$-abelian, then it is an Heegaard Floer L-space.
\end{conj}
The relevance of Conjecture \ref{conj: SU(2)-ab implies HF L-space} stems from both examples and theoretical expectations.  
On the one hand, Conjecture \ref{conj: SU(2)-ab implies HF L-space} is known to hold for several $3$-manifolds: for instance, lens spaces and the double cover of $S^3$ branched along the knot $8_{18}$ are both $L$-spaces and, at the same time, $SU(2)$-abelian \cite[Section 8]{classofknots}.  
On the other hand, as discussed in \cite{InstantoLSpaceAndSplicing}, one expects that a rational homology $3$-sphere fails to be $SU(2)$-abelian precisely when it is not an instanton $L$-space. Moreover, \cite[Conjecture 7.24]{kronheimer2010knots} claim that a rational homology $3$-sphere is an instanton $L$-space if and only if it is a Heegaard Floer $L$-space. Summarising, the expected chain of implications is:
\[
   \text{$SU(2)$-abelian} \ \Longrightarrow \ 
   \text{instanton $L$-space} \ \Longrightarrow \ 
   \text{Heegaard Floer $L$-space}.
\]
The converse of Conjecture \ref{conj: SU(2)-ab implies HF L-space} does not hold: suppose that $K\subset S^3$ is the alternating pretzel knot of type $P(2,3,7)$. The double cover $\Sigma_2(K)$ of $S^3$ branched along $K$ is a Heegaard Fleor L-space by \cite[Proposition 3.3]{alternatingimpliesLpsace}. It is known that $\Sigma_2(K)$ is a Seifert fibered space and by \cite[Theorem 1.2]{Menangerie} it is not $SU(2)$-abelian. Therefore $\Sigma_2(K)$ is an example of a Heegaard Floer L-space that is not $SU(2)$-abelian. In general, Conjecture \ref{conj: SU(2)-ab implies HF L-space} provides a criterion to decide the $SU(2)$-abelian status of a $3$-manifold: if $Y$ is not an Heegaard Floer L-space, then we expect $Y$ to be not $SU(2)$-abelian. 

Theorem \ref{teo: sigma small then not SU(2)-abelian} has a consequence that the author finds interesting.
Let us assume that $Y$ is a toroidal rational homology $3$-sphere with $|H_1(Y;\mathbb{Z})|=5$.
If $Y$ is not a Heegaard Floer L-space, then by Conjecture \ref{conj: SU(2)-ab implies HF L-space} we expect $Y$ to be not $SU(2)$-abelian.
As an application of \cite[Theorem 7.20]{BorderedHeegaardHomology}, if $Y$ is a toroidal Heegaard Floer L-space, then $Y$ is a toroidal graph manifold. Therefore, by Theorem \ref{teo: sigma small then not SU(2)-abelian} the manifold $Y$ is not $SU(2)$-abelian.
We conclude that, if $Y$ is a toroidal $3$-manifold with $|H_1(Y;\mathbb{Z})|=5$, then we expect $Y$ to be not $SU(2)$-abelian.
Similarly, if $Y$ is a toroidal manifold with $|H_1(Y;\mathbb{Z})|=6$, then $Y$ is not a Heegaard Floer L-space by \cite[Theorem 7.20]{BorderedHeegaardHomology}. According to Conjecture \ref{conj: SU(2)-ab implies HF L-space}, we expect $Y$ to not be $SU(2)$-abelian. We summarise this by proposing the following as a refinement of \cite[Conjecture 1.5]{InstantoLSpaceAndSplicing}:
\begin{conj}
    Let $Y$ be a toroidal rational homology $3$-sphere. If $|H_1(Y;\mathbb{Z})| \le 6$, then $Y$ is not $SU(2)$-abelian.
\end{conj}

We now give an outline of the techniques that we will use in this work. Let $Y_1$ and $Y_2$ be two $3$-manifolds with torus boundary. Let $\Sigma$ be a separating torus in a closed $3$-manifold $Y$. If the torus $\Sigma \subset Y$ is such that 
\[
    \overline{Y \setminus \Sigma} = Y_1 \cup Y_2,
\]
then we write that $Y=Y_1 \cup_{\Sigma} Y_2$.
Let us now assume that $Y=Y_1\cup_{\Sigma}Y_2$. Note that in this case $\Sigma \subset Y$ coincides with $\partial Y_1 = \partial Y_2$.
We denote by $\Ui$ the subgroup of diagonal $SU(2)$-matrices. We define $\R{\Sigma}$ as $\text{Hom}(\fund{\Sigma},\Ui)$. Furthermore, we define $T(Y_i,\partial Y_i) \subset \R{\Sigma}$ as the set of representations $\fund{\Sigma}\to \Ui$ that extend to a representation $\fund{Y_i} \to SU(2)$.
As a consequence of \cite[Theorem 3.7]{Mino}, the $SU(2)$-abelian status of $Y=Y_1 \cup_{\Sigma} Y_2$ is completely determined by the intersection
\begin{align}\label{eq: una volta intersezione}
    T(Y_1,\partial Y_1)\cap T(Y_2,\partial Y_2) \subset \R{\Sigma}.
\end{align}
Details can be found in Section \ref{sec: review of T(Y,Y)}. Therefore, we prove Theorem \ref{teo: sigma small then not SU(2)-abelian} by assuming that $Y=Y_1\cup_{\Sigma}Y_2$ and by studying $T(Y_1,\partial Y_1)$, $T(Y_2,\partial Y_2)$, and the intersection in \eqref{eq: una volta intersezione}.
Roughly speaking, the above-mentioned \cite[Theorem 3.7]{Mino} implies that the $SU(2)$-abelian status of a manifold $Y$ can be determined by splitting $Y$ along a separating torus $\Sigma$ into two pieces $Y_1$ and $Y_2$, and by analyzing the $SU(2)$-representations of its two
constituents.
In the case of a graph manifold rational homology $3$-sphere, as we shall see in Corollary \ref{cor: the useful decomposition}, the separating torus $\Sigma$ can be chosen so that 
\[
    Y = Y_1 \cup_{\Sigma} Y_2,
\]
where $Y_1$ is a Seifert fibered space admitting a Seifert fibration with disk base, and $Y_2$ is a graph manifold with the same rational homology as $S^1 \times \mathbb{D}^2$.
The notation for Seifert fibered spaces can be found in Subsection \ref{subsection: SFS}.
The space $T(Y_1,\partial Y_1)$ is fully described after \cite{Mino} and the author’s Ph.D. thesis, therefore in this work we shall study the $SU(2)$-representations for graph rational solid tori. In line with this, a key result is the following:
\begin{teo}\label{teo: admissible pieces have the property P}
    Let $Y$ be a graph manifold rational homology solid torus. If $Y$ is not a solid torus, then there exists a path $\gamma \colon (-1,1) \to \R{\partial Y}$ of representations such that:
    \begin{enumerate}
        \item the path starts and ends at representations $\fund{\partial Y} \to \Ui$ that extend to an abelian representation $\fund{Y} \to SU(2)$;
        \item for every $x \in (-1,1)$ the representation $\gamma(x)$ extends to an irreducible representation $\fund{Y} \to SU(2)$;
        \item the path $\gamma$ is non-trivial in $\R{\partial Y}$ and consists of a line segment.
        \end{enumerate}
\end{teo}
In this context, we say that a path in $\R{\partial Y}$ is trivial if it is homotopic, relative to its endpoints, to a constant path. A similar property is proved in \cite[Theorem 7.1]{IntegerHomSL2CIrrResp} for the exterior of a knot $K \subset S^3$.   

Let us now highlight a couple of results that the author finds interesting and that are consequences of Theorem \ref{teo: admissible pieces have the property P}.
Let $Y_1 \neq S^1 \times \mathbb{D}^2$ and $Y_2\neq S^1 \times \mathbb{D}^2$ be two $3$-manifolds with torus boundary, we denote by
$\lambda_i \subset \partial Y_i$ the rational longitude of $Y_i$. We call $o_i \in \mathbb{N}$ the order of $\lambda_i$ in $H_1(Y_i;\mathbb{Z})$.
Let $o_1 =1$, $o_2 \in \{1,3\}$, and $Y$ be the result of gluing $Y_1$ and $Y_2$ in a way that $\geo{\lambda_1}{\lambda_2}=1$.
According to \cite[Theorem 1.1]{boyer2023slopedetectiontoroidal3manifolds} and \cite[Theorem 1.3]{boyer2023slopedetectiontoroidal3manifolds}, the manifold $Y$ is not an Heegaard Floer L-space.
The following propositions state that such a $Y$ is not $SU(2)$-abelian, as suggested by Conjecture \ref{conj: SU(2)-ab implies HF L-space}.

\begin{prop}\label{prop: sigma =1 ish}
     Let $Y_1$ and $Y_2$ be two graph manifolds with torus boundary and null-homologous rational longitudes. If neither $Y_1$ nor $Y_2$ is a solid torus and $Y=Y_1 \cup_{\Sigma} Y_2$ is a rational homology sphere such that $\Delta(\lambda_1,\lambda_2)= 1$, then $Y$ is not $SU(2)$-abelian.
\end{prop}

\begin{cor}\label{cor: o1=1 and o2=3 and delta=1}
    Let $Y_1$ and $Y_2$ be two graph manifolds with torus boundary. If neither $Y_1$ nor $Y_2$ is a solid torus, $o_1=1$, $o_2=3$, and $Y=Y_1 \cup_{\Sigma} Y_2$ is a rational homology $3$-sphere such that $\Delta(\lambda_1,\lambda_2)= 1$, then $Y$ is not $SU(2)$-abelian.
\end{cor}

\subsection*{Acknowledgements}
I am deeply grateful to my research advisors, Duncan McCoy and Steven Boyer, for their constant patience, for guiding me along the path of my Ph.D. and for the useful discussions during which many key ideas were discussed. I would also like to thank Lorenzo Campioni, Pietro Capovilla, and Simone Maletto for human support they gave me during the drafting of this work and Emanuele Ronda for the help with the titles of the sections.
Lastly, I would like to thank the anonymous referees for their helpful comments.

\section{Notation}
\label{sec: notation}
Recall that for any irreducible orientable $3$-manifold $Y$, there is minimal collection $\mathscr{T}$ of properly embedded disjoint essential tori such that each component of $Y \setminus \mathscr{T}$ is either a hyperbolic or a Seifert fibered manifold, and such a collection is unique up to isotopy (see \cite{Jaco1979SeifertFS}). The \emph{JSJ decomposition} of $Y$ is given by
\[
Y= Y_0 \sqcup \cdots \sqcup Y_n
\]
where each $Y_i$ is the closure of a component of $Y \setminus \mathscr{T}$. A manifold $Y_i$ is called a \emph{JSJ piece} of $Y$ and a torus in the collection $\mathscr{T}$ is called a \emph{JSJ torus} of $Y$.

\begin{defn}
    Let $Y$ be a compact, irreducible, and orientable $3$-manifold. We say that $Y$ is a graph manifold if every of its JSJ pieces is a Seifert fibered manifold.
\end{defn}

Let $M$ be a manifold with torus boundary. A (rational) slope on the boundary $M$ is an element $[\alpha]$ of the projective space of $H_1(\partial M;\mathbb{Q})$, where $\alpha \in H_1(\partial M;\mathbb{Q}) \setminus \{ 0 \} $. Slopes can be identified in two equivalent ways:
 \begin{itemize}
     \item a $\pm$-pair of primitive elements of $H_1(\partial M; \mathbb{Z})\equiv \fund{\partial M}$;
    \item a $\partial M$-isotopy class of essential simple closed curves on $\partial M$. 
 \end{itemize}
For further details, see \cite[Subsection 4.2]{boyer2022orderdetection}.

Let $\Sigma$ be a torus, and $\{\mu,\lambda\}$ a basis of $H_1(\Sigma;\mathbb{Z})$. We use the convention, depending on the choice of the basis $\{\mu,\lambda\}$, such that an element $\nicefrac{p}{q} \in \mathbb{Q}\cup \{\nicefrac{1}{0}\}$ corresponds to the slope $p \mu + q \lambda \in H_1(\Sigma;\mathbb{Z})$. The \emph{distance} between two slopes $\nicefrac{p}{q}$ and $\nicefrac{r}{s}$ is $\Delta(\nicefrac{p}{q},\nicefrac{r}{s})=|ps-rq|$ and it corresponds to the absolute value of the algebraic intersection number between curves representing $\nicefrac{p}{q}$ and $\nicefrac{r}{s}$.

If $Y$ is a compact $3$-manifold with toroidal boundary components $\Sigma_1, \cdots, \Sigma_n$ with a fixed basis $\{\mu_i,\lambda_i\}$ for $H_1(\Sigma_i;\mathbb{Z})$ for each $i$, then
\[
    Y\left(\Sigma_1, \cdots, \Sigma_n;\nicefrac{r_1}{s_n},\cdots \nicefrac{r_n}{s_n}\right)
\]
denotes the closed $3$-manifold obtained by performing Dehn fillings along a simple closed curve representing $\nicefrac{r_i}{s_i}$ on $\Sigma_i$ for each $i=1,\cdots, n$. Moreover, if $Y$ has torus boundary and $\alpha \subset \partial Y$ is a slope, we make the notation lighter by defining
\[
    Y\left(\alpha\right) \coloneq Y(\partial Y;\alpha).
\]

Let $Y$ be a $3$-manifold with torus boundary, we define the \emph{rational longitude} of $Y$ as the unique slope $\lambda_{Y} \subset Y$ such that $[\lambda_{Y}]$ is a torsion element of $H_1(Y; \mathbb{Z})$.
In order to make the notation lighter, if $Y_i$ is a $3$-manifold with torus boundary, then, when it is not confusing, we denote by $\lambda_i \subset \partial Y_i$ the rational longitude of $Y_i$ and by $o_i$ the order of $\lambda_i$ in $H_1(Y_i;\mathbb{Z})$.

The subgroup of $SU(2)$ made of the diagonal matrices is called $\Ui$. In order to make the notation a little lighter, when not confusing, we use the following notation:
\[
e^{i\theta} \coloneqq
\begin{bmatrix}
    e^{i\theta} & 0 \\ 0 & e^{-i\theta}
\end{bmatrix} \in \Ui \le SU(2).
\]

Let $z$ be a non-central element of $SU(2)$. We denote by $\Lambda_z \subset SU(2)$ the centralizer subgroup of the element $z$. The subgroups $\Lambda_z$ and $\Ui$ are known to be conjugate in $SU(2)$, details can be found in \cite[Lecture 13]{saveliev}. The following is well known.
\begin{fact}\label{Fact: centralizers}
Let $x,y$ be two non-central elements of $SU(2)$. The elements $x$ and $y$ commute if and only if the two centralizer subgroups $\Lambda_x$ and $\Lambda_y$ coincide.
\end{fact}

\subsection{Seifert fibered spaces}\label{subsection: SFS}
A compact space $Y$ is said to be a \emph{Seifert fibered space}, or
\emph{Seifert fibered manifold},
if there exists a collection of pairwise disjoint circles $f_\alpha \subset Y$, which are called \emph{fibers}, such that
\[
    Y = \bigcup_{\alpha} f_\alpha,
\]
and every fiber $f_\alpha$ admits a tubular neighbourhood in $Y$ that is union of fibers. We call such a tubular neighbourhood a \emph{fibered neighbourhood}.

Every fiber $f_\alpha$ is associated with an invariant called \emph{order} that counts how many times any other close by fiber wraps around $f_\alpha$.
A fiber of order one is called \emph{regular fiber} and the ones with order $2$ or greater are called \emph{singular fibers}.
A compact Seifert fibered space admits finitely many singular fibers.

\begin{defn}
    Let $Y$ be a Seifert fibered space and let $S \subset Y$ be an embedded or immersed surface. We say that $S$ is \emph{vertical} if it is a union of fibers of $Y$.
\end{defn}

Let $Y$ be a Seifert fibered space. It can be proven that 
\[
    Y/\sim \quad \text{with} \quad x\sim y \text{ if and only if }x,y \in f_\alpha,
\]
is a $2$-dimensional \emph{orbifold} $\mathcal{B}$ with a number of cone points equal to the number of the singular fibers, each of order equal to the order of the corresponding singular fiber.
We give \cite{ScottGoemtries} and \cite{orbifold} as references for the concept of orbifold.
Let $B$ be the underlying surface of the orbifold $\mathcal{B}$ and $\{p_1,\cdots, p_n\}$ the orders of its cone points. The orbifold $\mathcal{B}$ is sometimes denoted by
\[
    B(p_1,\cdots,p_n).
\]
We say that $Y$ is \emph{Seifert fibered} over the orbifold $\mathcal{B}=B(p_1,\cdots,p_n)$ and, when the Seifert fibration of $Y$ is fixed, that $B$ is the \emph{base space} of $Y$. To be even more accurate, when we say that $Y$ has base space $B$ we mean that there exists a Seifert fibration of $Y$ whose base space is $B$.

Let $n\ge 1$. A Seifert fibered manifold with $n$ singular fibers can be described from the surface $B$ and a set of fractions $\{\nicefrac{p_1}{q_n}, \cdots, \nicefrac{p_n}{q_n}\}$, where $p_i$ is the order of the $i^{th}$ singular fiber. We call the fraction $\nicefrac{p_i}{q_i}$ a \emph{Seifert coefficient} of $Y$. We can therefore write
\[
    Y= B\left(\frac{p_1}{q_1},\cdots, \frac{p_n}{q_n}\right).
\]
In this case, $Y$ is Seifert fibered over the orbifold $B(p_1,\cdots,p_n)$ and has base space $B$.
Details can be found in \cite[Section 10.3.2]{martelli}.

\section{Preliminaries and background on  \texorpdfstring{$T(Y,\partial Y)$}{T(Y,Y)}}
\label{sec: review of T(Y,Y)}

This section is devoted to recalling the fundamental definitions and some properties of $T(Y,\partial Y)$ for a $3$-manifold with torus boundary $Y$. We use \cite{Mino} as the main reference for this. 

Let $\Sigma$ be a torus, we define the space $\R{\Sigma}$ as $\text{Hom}(\fund{\Sigma},\Ui)$.
Let $\{x,y\} $ be an ordered basis of $\fund{\Sigma}$. We give the space $\R{\Sigma}$ coordinates $(\theta,\psi)$ according to this basis in the following way:
the point $(\theta,\psi) \in \R{\Sigma}$ corresponds to the unique representation $\fund{\Sigma} \to \Ui$ such that
\begin{equation}\label{eq: general coordinates}
    x \mapsto e^{i\theta} \quad \text{and} \quad y\mapsto e^{i \psi}.
\end{equation}
The correspondence in \eqref{eq: general coordinates} implies that the space $\R{\Sigma}$ is diffeomorphic to $S^1 \times S^1$.

Let $Y$ be a $3$-manifold with boundary, define $R(Y)$ as $\text{Hom}(\fund{Y},SU(2))$. The \emph{representation space of $Y$ relative to $\partial Y$}, denoted by $T(Y,\partial Y)$, is defined as the set of representations $\fund{\partial Y} \to \Ui$ that extend to an $SU(2)$-representation of $\fund{Y}$.
By definition, there is an inclusion $T(Y,\partial Y) \subset \R{\partial Y}$. We recall that a non-abelian $SU(2)$-representation is also called \emph{irreducible}.

\begin{defn} \label{defn: A1 H1 e P1}
    Let $Y$ be a $3$-manifold with torus boundary, let $\iota \colon \partial Y \to Y$ be the topological natural inclusion. We define the sets $A(Y), H(Y)$, and $P(Y)$ as:
    \begin{alignat*}{2}
        A(Y)& \coloneqq \left\{ \eta \in \R {\partial Y} \,\middle|\, \exists \rho \in R(Y) \text{ such that } \restr{\rho}{\iota_\ast \fund{\partial Y}}\equiv \eta\text{ and $\rho$ is abelian} \right\}, \\
        H(Y)& \coloneqq \left\{ \eta \in \R{\partial Y} \,\middle|\, \exists \rho \in R(Y) \text{ such that } \restr{\rho}{\iota_\ast \fund{\partial Y}}\equiv \eta\text{ and $\rho$ is irreducible} \right\}, \\
        P(Y)& \coloneqq \Big\{ \eta \in \R{\partial Y} \,\Big|
            \begin{aligned}[t]
            &\exists \rho \in R(Y) \text{ such that } \restr{\rho}{\iota_\ast \fund{\partial Y}}\equiv \eta,\\
            &\text{$\eta$ is central, and $\rho$ is abelian and non-central} \Big\}.
    \end{aligned}
    \end{alignat*}
\end{defn}
Equivalently, a representation $\fund{\partial Y} \to \Ui$ is in $A(Y)$ (resp. $H(Y)$) if and only if it extends to an abelian (resp. an irreducible) representation $\fund{Y} \to SU(2)$. Similarly, a representation $\fund{\partial Y} \to \ZSU$ is in $P(Y)$ if and only if it extends to an abelian representation $\fund{Y} \to SU(2)$ whose image is not in $\ZSU$. 
Notice that $P(Y) \subset A(Y)$ and $A(Y)\cup H(Y)=T(Y,\partial Y) \subseteq \R{\partial Y}$.

\begin{cor}[{\cite[Corollary 3.9]{Mino}}]\label{cor: A(Y) with rational longitude}
    Let $Y$ be a $3$-manifold with torus boundary. Then
    \[
    A(Y)= \left\{\eta \in \R{\partial Y} \middle| \eta(\lambda_Y)^{o_Y}=1\right\} \subset \R{\partial Y}.
\]
\end{cor}

If $Y=Y_1 \cup_{\Sigma} Y_2$, then in $Y$ the tori $\partial Y_1$  and $\partial Y_2$ are both equal to $\Sigma$. Therefore, $T(Y_1,\partial Y_1)$ and $T(Y_2,\partial Y_2)$ are both subsets of the variety $\R{\Sigma}$. This implies that the intersection $T(Y_1,\partial Y_1) \cap T(Y_2,\partial Y_2)$ is well defined in $\R{\Sigma}$. Hence we have the following result:

\begin{teo}[{\cite[Theorem 3.7]{Mino}}]\label{teo: M SU(2)-abeliano se e solo se i pezzi sono empty}
    Let $Y_1$ and $Y_2$ be two $3$-manifolds with torus boundary. The manifold $Y=Y_1 \cup_{\Sigma} Y_2$ is $SU(2)$-abelian if and only if $H(Y_1)\cap H(Y_2)$, $H(Y_1) \cap A(Y_2)$, $A(Y_1) \cap H(Y_2)$, and $P(Y_1) \cap P(Y_2)$ are empty.
\end{teo}

\begin{cor}[{\cite[Corollary 3.11]{Mino}}]\label{cor: Y SU(2) allora i due pezzi sono SU(2)-abelian}
    If the manifold $Y = Y_1 \cup_{\Sigma} Y_2$ is $SU(2)$-abelian, then the manifold $Y_1(\lambda_2)$
    is $SU(2)$-abelian.
\end{cor}

The next Lemma was previously proven in \cite{LinHolonomyPerturbation} with a different notation, we report the proof for the sake of completeness.

\begin{lemma}\label{lemma: translation}
    Let $Y$ be a rational homology solid torus with torus boundary such that the rational longitude $\lambda_Y$ has odd order $o_Y$. Let $\xi$ be a slope such that $\{\xi,\lambda_Y \}$ is a basis for $\fund{\partial Y}$. Let $\R{\partial Y}$ be given with coordinates $(\theta,\psi)$ according to this ordered basis as in \eqref{eq: general coordinates}.
    The sets $A(Y)$ and $H(Y)$ are invariant, as subsets of $\R{\partial Y}$, by translation of vector $(\pi,0)$ with respect to the $(\theta,\psi)$-coordinates.
    \begin{proof}
        Since $Y$ is a rational homology solid torus, we have that
        \[
            H_1(Y;\mathbb{Z}) = \mathbb{Z} \oplus \text{Torsion}\quad \text{and} \quad \frac{H_1(Y;\mathbb{Z})}{\text{Torsion}}= \mathbb{Z}.
        \]
        Let $p \colon \fund{Y} \to \mathbb{Z}$ be the composition of the abelianisation map $\fund{Y} \twoheadrightarrow H_1(Y;\mathbb{Z})$ and the quotient map $H_1(Y;\mathbb{Z}) \twoheadrightarrow H_1(Y;\mathbb{Z})/\text{Torsion}$. The homomorphism $p$ is a composition of two surjective homomorphisms and therefore it is surjective.
        
        Let $z \in \mathbb{Z}$ be a generator of $\ima p$. Therefore, every element of $\ima p$ is of the form $n \cdot z$, with $n \in \mathbb{Z}$.
        A standard homology computation shows that $p(\xi)=  o_Y \cdot z \in \ima p$ and that $p(\lambda_Y)=0$.
        We define the homomorphism
        \[
            \omega \colon \frac{H_1(Y;\mathbb{Z})}{\text{Torsion}}=\ima{p} \to \{\pm 1\} \quad \text{by} \quad n \cdot z \mapsto (-1)^n.
        \]
        Let $\chi \colon \fund{Y} \to \{\pm 1\} = \ZSU$ be the composition $\omega \circ p$. We notice that 
        \[
        \chi(\xi)=(\omega \circ p) (\xi) = \omega (o_Y \cdot z)=(-1)^{o_Y}=-1 \quad \text{and that} \quad \chi(\lambda_Y)=1.
        \]
        
        Let $\eta$ be a representation in $A(Y)$ (resp. $H(Y)$), then there exists an abelian (resp. irreducible) representation $\rho \in R(Y)$ such that $\restr{\rho}{\fund{\partial Y}} \equiv \eta$. Let us assume that $\eta \in \R{\partial Y}$ has coordinates $(\theta_0,\psi_0) \in [0,2\pi]^2$.
        We define $\rho' \colon \fund{Y} \to SU(2)$ by the multiplication
        \[
            x \mapsto \chi(x) \rho(x).
        \]
        Clearly $\rho' \in R(Y)$ and it is abelian (resp. irreducible) if and only if $\rho$ is abelian (resp. irreducible). Thus, $\restr{\rho'}{\fund{\partial Y}}$ is a point of $A(Y)$ (resp. $H(Y)$). We notice that
        \[
        \rho'(\lambda_Y)=\rho(\lambda_Y)=e^{i\psi_0} \quad \text{and} \quad\rho'(\xi)=\chi(\xi)\rho(\xi)=-\rho(\xi) = e^{i \pi} e^{i\theta_0}=e^{i (\pi + \theta_0)}.
        \]
        Thus, if $(\theta_0,\psi_0)$ is a point of $A(Y)$ (resp. $H(Y)$), then $(\theta_0 + \pi,\psi_0)$ is a point of $A(Y)$ (resp. $H(Y)$).
    \end{proof}
\end{lemma}

\begin{defn}\label{defn: jewelled representation}
    Let $G$ be a group and $H$ a linear group \footnote{A linear group is a group of invertible matrices over a given field and with the operation of matrix multiplication.}. For a given matrix $M$, we denote by $M^\top$ the transpose of a $M$. Let $\rho \colon G \to H$ be a representation. We define the \emph{jewelled representation} $\rho^\dagger : G \to H$ as
    \[\rho^\dagger(x)=\left(\rho(x)^{-1}\right)^\top = \left(\rho(x)^{\top}\right)^{-1}.
    \]
\end{defn}

In the following we assume that $Y$ is a $3$-manifold with torus boundary and that an ordered basis $\{x_1,x_2\}$ for $\fund{\partial Y}$ is given. Therefore, we parameterize $\R{\partial Y}$ with coordinates $(\theta,\psi)$ as in \eqref{eq: general coordinates} according to this ordered basis.

\begin{cor}\label{cor: jewelled repr}
    Let $Y$ be a $3$-manifold with torus boundary. If $(\theta_0,\psi_0) \in T(Y,\partial Y) \subset \R{\partial Y}$, then $(2\pi - \theta_0,2\pi-\psi_0) \in T(Y,\partial Y) \subset \R{\partial Y}$. In particular, if $(\theta_0,\psi_0)$ is in either $H(Y)$, $A(Y)$, or $P(Y)$, then $(2\pi - \theta_0,2\pi-\psi_0)$ is in either $H(Y)$, $A(Y)$, or $P(Y)$ respectively.
    \begin{proof}
        Let $\eta \colon \fund{\partial Y} \to \Ui$ be the representation corresponding to $(\theta_0,\psi_0)$. Let $\rho\colon \fund{Y} \to SU(2)$ be the extension of $\eta$.
        The jewelled representation $\rho^\dagger$ is such that
        \[
        \restr{\rho^\dagger}{\fund{\partial Y}} = \eta^\dagger.
        \]
        In particular, $\eta^\dagger \in T(Y,\partial Y)$. Moreover, since $\eta(x_1)$ and $\eta(x_2)$ are diagonal matrices, we have that $\eta(x_1)^\top=\eta(x_1)$ and $\eta(x_2)^\top=\eta(x_2)$. This implies that
        \[
            \eta^\dagger(x_1)=\left(e^{i \theta_0}\right)^{-1}=e^{i (2\pi - \theta_0)} \quad \text{and} \quad  \eta^\dagger(x_2)=\left(e^{i \psi_0}\right)^{-1}=e^{i (2\pi - \psi_0)}.
        \]
        Therefore, $(2\pi-\theta_0,2\pi -\psi_0) \in T(Y,\partial Y) \subset \R{\partial Y}$. It is easy to see that $\rho^\dagger$ is abelian (resp. irreducible, resp. abelian and non-central) if and only if $\rho$ is abelian, (resp. irreducible, resp. abelian and non-central).
    \end{proof}
\end{cor}

\section{The manifolds \texorpdfstring{$C_2$}{C3} and \texorpdfstring{$C_3$}{C2}}
\label{sec: two important pieces}

In this section we introduce the manifolds $C_2$ and $C_3$, and we prove some of their basic properties. These two manifolds will play a central role in this work, as they form the constituents of the decomposition of Corollary \ref{cor: the useful decomposition}, where we prove that a graph manifold rational homology $3$-sphere can be obtained as a union of Seifert fibered spaces with torus boundary together with copies of $C_2$ and $C_3$. 

\begin{figure}
    \centering
    \includegraphics[width=0.3\linewidth]{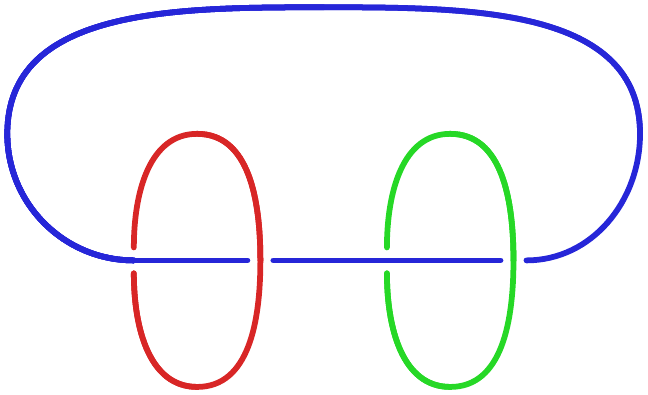}
    \caption{The connected sum of two Hopf links.}
    \label{fig:Link}
\end{figure}

Let $L \subset S^3$ be the $3$-components link obtained as the connected sum of two Hopf links. The link $L$ is represented in Figure \ref{fig:Link}.
We define $C_3$ as the closed manifold $S^3 \setminus \nu(L)$, where $\nu(L)$ is an open tubular neighborhood of the link $L$. The manifold $C_3$ is sometimes called \emph{3D pair of pants}.

Another way to define $C_3$ is as follows. Let $\mathbb{D}^2$ be a closed disk of radius $1$. Let $D_1$ and $D_2$ be two open disks of $\mathbb{D}^2$, each of radius $\nicefrac{1}{4}$ and centered in $(0,\nicefrac{1}{2})$ and $(0,-\nicefrac{1}{2})$ respectively. The manifold $C_3$ is
\[
    C_3 = \left(\mathbb{D}^2 \times S^1\right) \setminus (D_1 \cup D_2) \times S^1 = \left(\mathbb{D}^2 \setminus (D_1 \cup D_2) \right) \times S^1.
\]

The manifold $C_3$ is a Seifert fibered manifold with fundamental group equal to
\begin{equation}\label{eq: fundamental group C3}
    \fund{C_3} = \left\langle x_1,x_2,x_3,h \middle| [x_j,h], x_1x_2=x_3\right\rangle,
\end{equation}
where $h$ is a regular fiber of the Seifert fibration and $[x,y]$ denotes the commutator element.

The boundary of $C_3$ is made of three torus components, namely $\Sigma_1$, $\Sigma_2$, and $\Sigma_3$. For $j \in \{1,2,3\}$, the ordered basis $\{x_j,h\}$ generates the fundamental group $\fund{\Sigma_j}$. We endow $\R{\Sigma_j}$ with coordinates $(\theta_j,\psi_j)$ as in \eqref{eq: general coordinates}:
the point $(\theta_j,\psi_j) \in \R{\Sigma_j}$ corresponds to the unique representation $\fund{\Sigma_j} \to \Ui$ such that
\begin{equation}\label{eq: coordinates for C3}
    x_j \mapsto e^{i\theta_j} \quad \text{and} \quad h\mapsto e^{i \psi_j}.
\end{equation}

\begin{notation}
    When using $C_3$, we call $h$ the regular fiber of $C_3$. When $h$ is considered in $\Sigma_j$, we call it $h_j$ for $j \in \{1,2,3\}$.
\end{notation}

We define $C_2$ as the twisted I-bundle over a once-punctured M\"obius band.
The manifold $C_2$ is a Seifert fibered manifold with fundamental group equal to
\begin{equation}\label{eq: fundamental group C2}
    \fund{C_2} = \left\langle x_1,x_2,z,h \middle| [x_j,h], x_1x_2z^2, zhz^{-1}h\right\rangle,
\end{equation}
where $h$ is the regular fiber of the fibration.

The manifold $C_2$ has two toroidal boundary components, namely $\Sigma_1'$, $\Sigma_2'$. For $j \in \{1,2\}$, the ordered basis $\{x_j,h\}$ generates the fundamental group $\fund{\Sigma_j'}$. As before, we give coordinates $(\theta_j,\psi_j)$ to the variety $\R{\Sigma_j'}$ as in \eqref{eq: general coordinates}:
the point $(\theta_j,\psi_j) \in \R{\Sigma_j'}$ corresponds to the unique representation $\fund{\Sigma_j'} \to \Ui$ such that
\begin{equation}\label{eq: coordinates for C2}
    x_j \mapsto e^{i\theta_j} \quad \text{and} \quad h\mapsto e^{i \psi_j}.
\end{equation}

\begin{notation}
    When using $C_2$ we call $h$ the regular fiber of $C_2$. When $h$ is considered in $\Sigma_j'$, we call it $h_j$ for $j \in \{1,2\}$.
\end{notation}

Let $Y_1$ and $Y_2$ be two manifolds with torus boundary. When we write $Y=Y_1 \cup_{\Sigma_1} C_3 \cup_{\Sigma_2} Y_2$, we mean that $Y$ is a manifold with torus boundary admitting two embedded tori $\Sigma_1$ and $\Sigma_2$ such that
\[
    \overline{Y \setminus (\Sigma_1 \cup \Sigma_2)}= Y_1 \cup C_3 \cup Y_2.
\]
In particular $\partial Y = \Sigma_3 \subset \partial C_3$. Figure \ref{Figure: the manifold Y1 C3 Y2} shows the manifold $Y=Y_1 \cup_{\Sigma_1} C_3 \cup_{\Sigma_2} Y_2$.

\begin{figure}[t]
            \centering
            \begin{overpic}[width=0.4\textwidth]{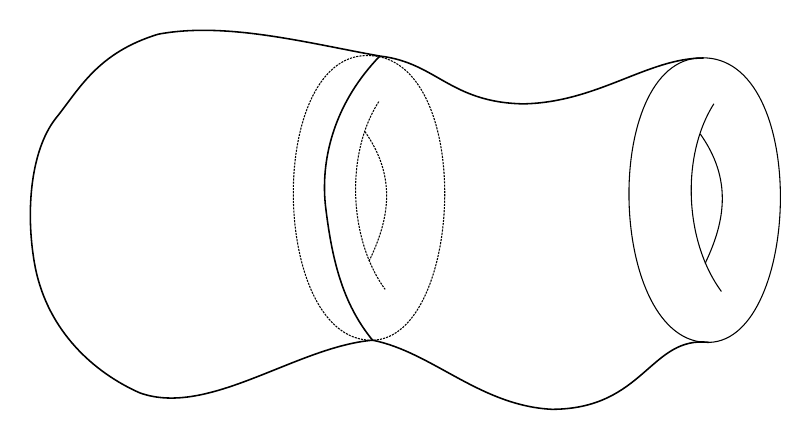}
            \put (21,27){$Y_1$}
            \put (63,25){$C_2$}
            \put (42,5){$\Sigma_1$}
            \end{overpic}
            \caption{The manifold $Y=Y_1 \cup_{\Sigma_1'} C_2$.}
            \label{Figure: the manifold Y1 C2}
\end{figure}

Similarly, when we write $Y=Y_1 \cup_{\Sigma_1'} C_2$, we mean that $Y$ is a manifold with torus boundary admitting an embedded torus $\Sigma_1'$ such that
\[
    \overline{Y \setminus \Sigma_1'}= Y_1 \cup C_2.
\]
Figure \ref{Figure: the manifold Y1 C2} shows the manifold $Y=Y_1 \cup_{\Sigma_1'} C_2$.
In particular $\partial Y = \Sigma_2' \subset \partial C_2$.

\begin{rmk}[{\cite[Lemma 7.1]{LSpaceTautFoliationsAndGraphManifolds}}]\label{rmk: order of toroidal}
    Let $Y=Y_1 \cup_{\Sigma} Y_2$ be a rational homology $3$-sphere.
    We obtain that
        \begin{align*}
        \left|H_1(Y;\mathbb{Z})\right|=o_1o_2t_1t_2 \geo{\lambda_1}{\lambda_2},
        \end{align*}
    where $o_i$ is the order of $\lambda_i$ in $H_1(Y_i,\mathbb{Z})$ and $t_i$ is the order of the torsion subgroup of $H_1(Y_i,\mathbb{Z})$. In particular $o_i$ divides $t_i$.
\end{rmk}

\begin{lemma}\label{lemma orders and torsion of C3}
    Let $Y_1$, $Y_2$, and $Y=Y_1 \cup_{\Sigma_1} C_3 \cup_{\Sigma_2} Y_2$ be rational homology solid tori with torus boundary. If $H_1(Y;\mathbb{Z})$ has the torsion subgroup of order one, then both $H_1(Y_1;\mathbb{Z})$ and $H_1(Y_2;\mathbb{Z})$ have torsion subgroup of order one. In particular, $\lambda_1$ and $\lambda_2$ are both null-homologous.
        \begin{proof}
            We denote by $\lambda_Y \subset \Sigma_3 = \partial Y$ the rational longitude of $Y$, by $o_Y$ its order in $H_1(Y;\mathbb{Z})$ and by $t_Y$ the order of the torsion subgroup of $H_1(Y;\mathbb{Z})$. 
            By hypothesis $t_Y=1$. We recall that $o_Y$ divides $t_Y$. Therefore $o_Y=1$.
            Similarly, for $i \in \{1,2\}$, we denote by $\lambda_i$, $o_i$, and $t_i$ the rational longitude of $Y_i$, its order in $H_1(Y_i;\mathbb{Z})$ and the order of the torsion subgroup of $H_1(Y_i;\mathbb{Z})$ respectively.
            
            The product $\geo{\lambda_1}{h_1}\geo{\lambda_2}{h_2}$ is either zero or non zero. We split the proof in these two cases.
            
            Let us assume that $\geo{\lambda_1}{h_1}\geo{\lambda_2}{h_2}=0$. Without loss of generality, we assume that $\geo{\lambda_1}{h_1}=0$.
            A standard homology computation shows that the regular fiber $h_3 \subset \Sigma_3$ is the rational longitude of $Y$, and therefore $\lambda_Y=h_3$ in $\Sigma_3=\partial Y$.
            With an abuse of notation, we call $x_3 \subset \Sigma_3$ the slope of $\Sigma_3$ whose homotopy class is
            \[
                x_3 \in \fund{C_3},
            \]
            where $\fund{C_3}$ is presented as in \eqref{eq: fundamental group C3}.
            Therefore $\geo{h_3}{x_3}=\geo{\lambda_Y}{x_3}=1$. Since $\geo{\lambda_Y}{x_3} \neq 0$, the manifold $Y(x_3)$ is a rational homology $3$-sphere.
            Remark \ref{rmk: order of toroidal} states that
            \[
            \left|H_1(Y(x_3);\mathbb{Z})\right| =o_Yt_Y\geo{\lambda_Y}{x_3}=1.
            \]
            
            It is straightforward to see that $C_3(\Sigma_3;x_3)$ is homeomorphic to $S^1 \times S^1 \times [0,1]$. Therefore,
             there exists a diffeomorphism $\varphi \colon \partial Y_1 \to \partial Y_2$ such that
            \[
            Y(x_3) = Y_1 \cup_{\varphi} Y_2.
            \]
            Again, Remark \ref{rmk: order of toroidal} implies that
            \[
            1 = \left|H_1(Y(x_3);\mathbb{Z})\right| = \left|H_1(Y_1 \cup_{\varphi} Y_2;\mathbb{Z})\right| = o_1t_1o_2t_2 \geo{\lambda_1}{\lambda_2}.
            \]
            This latter implies that $t_1=t_2=1$ and, since $o_i$ divides $t_i$, that $o_1=o_2=1$.
            
            Let us assume that $\geo{\lambda_1}{h_1}\geo{\lambda_2}{h_2} \neq 0$. For $i \in \{1,2\}$, let  $n_i,m_i$ be integers such that
            \[
                \lambda_1=n_1 x_1 + m_1 h_1 \subset \Sigma_1 \quad \text{and} \quad \lambda_2=n_2 x_2 + m_2 h_2 \subset \Sigma_2.
            \]
            Since $\geo{\lambda_1}{h_1}\geo{\lambda_2}{h_2}=n_1n_2 \neq 0$, the integers $n_1$ and $n_2$ are not null.
            Let us set $g=\gcd(n_1,n_2)$. Let $\sigma \subset \Sigma_3$ be the slope
            \[
                \sigma\coloneqq \frac{n_1 n_2}{gk} x_3 + \frac{n_1 m_2+n_2 m_1}{gk} h_3 \subset \Sigma_3 \quad \text{with} \quad k = \gcd \left( \frac{n_1n_2}{g}, \frac{n_1 m_2+n_2 m_1}{g}\right).
            \]
            It is easy to see that in $H_1(Y;\mathbb{Z})$, the element  $(o_1o_2k) \cdot \sigma$ is trivial. This implies that $\sigma \subset \Sigma_3$ is the rational longitude of $Y$, and therefore $\sigma=\lambda_Y$. We point out that $\geo{\lambda_1}{h_1}=n_1$, $\geo{\lambda_2}{h_2}=n_2$, and $\geo{\lambda_Y}{h_3}=\nicefrac{n_1n_2}{gk} \neq 0$.

            Since $\geo{\lambda_Y}{h_3} \neq 0$, the longitude $\lambda_Y \subset \partial \Sigma_3 = \partial Y$ is not the regular fiber $h_3 \subset \partial \Sigma_3 = \partial Y$. This implies that $Y(h_3)$ is a rational homology $3$-sphere. In particular
            \[
                Y(h_3)=Y_1(h_1) \# Y_2(h_2).
            \]
            This implies that $ H_1(Y(h_3);\mathbb{Z}) = H_1(Y_1(h_1);\mathbb{Z}) \oplus H_1(Y_2(h_2);\mathbb{Z})$. 
            Therefore, by Remark \ref{rmk: order of toroidal}, we obtain that
            \begin{align*}
            \left| H_1(Y(h_3);\mathbb{Z}) \right| & = \left| H_1(Y_1(h_1);\mathbb{Z}) \right| \cdot \left| H_1(Y_2(h_2);\mathbb{Z}) \right| \\
                o_Yt_Y\geo{\lambda_Y}{h_3} &= o_1t_1 \geo{\lambda_1}{h_1} o_2t_2 \geo{\lambda_2}{h_2} \\
                 \frac{n_1n_2}{gk} &= o_1t_1 n_1 o_2t_2 n_2 \\
                1 & = o_1t_1 o_2t_2 gk.
            \end{align*}
     Therefore, $t_1=t_2=o_1=o_2=1$.
        \end{proof}
\end{lemma}

\begin{lemma}\label{lemma: rational longitude of C2}
    Let $Y_1$ be a $3$-manifold with torus boundary and $Y=Y_1 \cup_{\Sigma_1'} C_2$. The regular fiber $h_2 \subset \Sigma_2$ is the rational longitude of $Y$. In particular, its order $o_Y$ in $H_1(Y;\mathbb{Z})$ is computed as follows:
    \[
    o_Y = \begin{cases}
        1 & \text{if $\geo{\lambda_1}{h_1} = 0$ and $o_1=1$,} \\
        2 & \text{otherwise}.
    \end{cases}
    \]
    \begin{proof}
        The presentation of $\fund{C_2}$ as in \eqref{eq: fundamental group C2} implies that $2[h]=2[h_2]$ is trivial in $H_1(C_2;\mathbb{Z})$. Therefore, $2[h_2]$ is trivial in $H_1(Y;\mathbb{Z})$. This implies that $h_2$ is the rational longitude of $Y$ and it has order either $1$ or $2$.

        If $\geo{\lambda_1}{h_1} = 0$ and $o_1=1$, then $[h_2]=[h_1]=[\lambda_1]$ in $H_1(Y;\mathbb{Z})$. Since $\lambda_1 \subset \Sigma_1$ is nullhomologous in $H_1(Y_1;\mathbb{Z})$, the element $[h_2]=[h_1]=[\lambda_1]$ is nullhomologous in $H_1(Y;\mathbb{Z})$. This implies that the order of the rational longitude of $Y$ is one

        If either $\geo{\lambda_1}{h_1} \neq 0$ or $o_1 \neq 1$, then $[h_2]=[h_1]=[h]$ is non-trivial in both $H_1(Y_1;\mathbb{Z})$ and $H_1(C_2;\mathbb{Z})$. Therefore the order of the rational longitude is not one. As the order is either $1$ or $2$, we get the conclusion.
    \end{proof}
\end{lemma}

\begin{rmk}\label{rmk: bases for C3}
    Let us consider $C_3$.
    If $\xi_1 \subset \Sigma_1$ and $\xi_2 \subset \Sigma_2$ are two slopes such that $\geo{\xi_1}{h}=\geo{\xi_2}{h}=1$, then $\{\xi_1,h\}$ and $\{\xi_2,h\}$ are bases for $\fund{\Sigma_1}$ and $\fund{\Sigma_2}$. In particular we can assume that $\fund{C_3}$ is presented as in \eqref{eq: fundamental group C3} and that
    \[
        [\xi_1]=x_1 \in \fund{\Sigma_1} \quad \text{and} \quad [\xi_2]=x_2 \in \fund{\Sigma_2}.
    \]
\end{rmk}

A similar consideration can be done for $C_2$ as well:
\begin{rmk}\label{rmk: bases for C2}
    Let us consider $C_2$.
    If $\xi_1 \subset \Sigma_1'$ is a slope such that $\geo{\xi_1}{h}=1$, then $\{\xi_1,h\}$ is a basis for $\fund{\Sigma_1'}$. Furthermore, we can assume that $\fund{C_2}$ is presented as in \eqref{eq: fundamental group C2} and that $[\xi_1]=x_1 \in \fund{\Sigma_1'}$.
\end{rmk}

\section{The topological setup}
\label{section: the topology set up}

We recall that we have made explicit the notation for Seifert fibered spaces in Subsection \ref{subsection: SFS}. We will prove in this section the existence of a decomposition into ‘simple’ pieces of graph manifolds rational homology $3$-spheres. That is, we will show that a graph manifold rational homology $3$-sphere admits a system of tori that decomposes it into a collection of copies of $C_3$, copies of $C_2$, and $\mathbb{D}^2(\nicefrac{p_\ast}{q_\ast})$;
where $\mathbb{D}^2(\nicefrac{p_\ast}{q_\ast})$ indicates a Seifert fibered space with disk base space. We emphasise that a solid torus $S^1 \times \mathbb{D}^2$ has disk base space, therefore solid tori are contained in the class $\mathbb{D}^2(\nicefrac{p_\ast}{q_\ast})$.

Let $Y$ be a graph manifold rational homology $3$-sphere and let $Y_i$ be a JSJ piece. Through the calculation of the first Betti number, we conclude that the Seifert fibered space $Y_i$ has a planar base space, i.e. either a punctured sphere $S^2$ or a punctured projective plane $\mathbb{RP}^2$. For details see \cite[Section 2.2]{boyer2017foliations}.

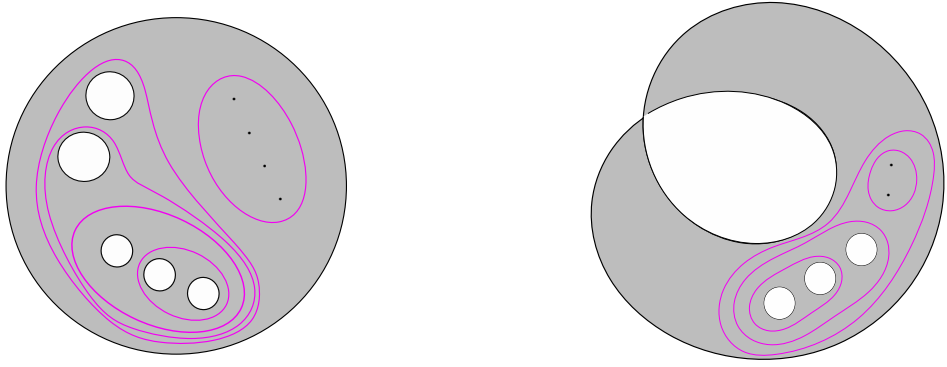
\begin{figure}[t]
    \centering
     \begin{subfigure}[b]{0.45\textwidth}
     \centering
     \begin{overpic}[width=0.7\textwidth]{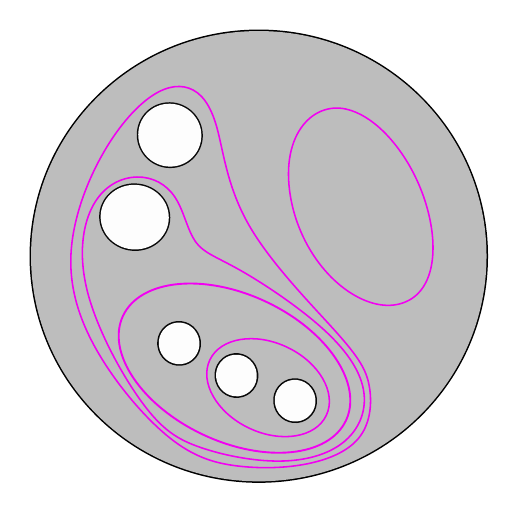}
            \put(77,46){\circle*{1}}
            \put(73,54.6){\circle*{1}}
            \put(69,63.2){\circle*{1}}
            \put(65,72){\circle*{1}}%
            \end{overpic}
            \caption{The case where $Y$ fibers over a multi-punctured $S^2$.}
    \end{subfigure} 
    \quad
    \begin{subfigure}[b]{0.45\textwidth}
    \centering
            \begin{overpic}[width=0.7\textwidth]{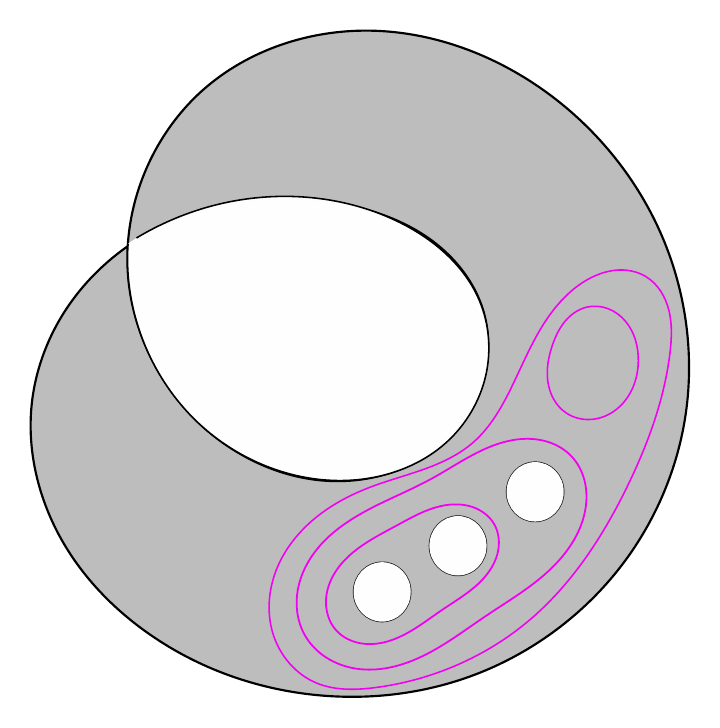}
            \put(80.5,46.5){\circle*{1}}
            \put(81.25,54){\circle*{1}}
            \end{overpic}
            \caption{The case where $Y$ fibers over a multi-punctured $\mathbb{RP}^2$.}
    \end{subfigure}
    \caption{The vertical tori of Lemma \ref{lemma: system for S2} and Lemma \ref{lemma: system for RP2}. The tori are shown in pink and the black dots represent the singular fibers.}
    \label{Figure: vertical tori for C3}
    \label{Figure: vertical tori for C2}
\end{figure}

\begin{lemma}\label{lemma: system for S2}
    Let $Y$ be a Seifert fibered space with boundary and base space a punctured $S^2$. There exists a collection of disjoint embedded tori $\tau \subset Y$ such that $\overline{Y \setminus \tau}$ is the union of a Seifert fibered space and a (possibly empty) collection of copies of $C_3$.
        \begin{proof}
        If $Y$ has one boundary component, then $Y=\mathbb{D}^2(\nicefrac{p_{\ast}}{q_{\ast}})$ and $\tau = \emptyset$.
        If $Y$ has two boundary components, then there exists a vertical torus $\Sigma \subset Y$ such that
    \[
        \overline{Y \setminus \tau} = C_3 \cup \mathbb{D}^2(\nicefrac{p_{\ast}}{q_{\ast}}),
    \]
        where $\mathbb{D}^2(\nicefrac{p_{\ast}}{q_{\ast}})$ is a Seifert fibered space with disk base space. In particular, $\mathbb{D}^2(\nicefrac{p_{\ast}}{q_{\ast}})=S^1\times \mathbb{D}^2$ if $Y$ has at most one singular fiber.
        Let us assume that $Y$ has $n \ge 3$ boundary components and $m \ge 0$ singular fibers.
        Let $\text{sgm} \colon \mathbb{N} \to \{1,0\}$ be the signum function.
        It is easy to see that there exists a system $\tau\subset Y$ of $n-2+ \text{sgm}(m)$ vertical tori such that
        \[
            \overline{Y \setminus \tau} = 
            \begin{dcases}
                \bigcup_{n-2} C_3 & \text{if $m=0$} \\
                \left(\bigcup_{n-1} C_3\right) \cup \mathbb{D}^2(\nicefrac{p_{\ast}}{q_{\ast}}) & \text{if $m \ge 1$}.
            \end{dcases}
        \]
        For example, Figure \ref{Figure: vertical tori for C3} shows the system of $5$ tori in the case $n=6$ and $m \ge 1$.
        \end{proof}
\end{lemma}

\begin{lemma}\label{lemma: system for RP2}
    Let $Y$ be a Seifert fibered space with boundary and base space a punctured real projective plane $\mathbb{RP}^2$. There exists a collection of disjoint embedded tori $\tau \subset Y$ such that $\overline{Y \setminus \tau}$ is the union of a Seifert fibered space with disk base, a copy of $C_2$, and a (possibly empty) collection of copies of $C_3$.
        \begin{proof}
            If $Y$ has one boundary component and no singular fibers, then
            according to \cite[Theorem 10.4.19]{martelli} $Y$ admits a Seifert fibration over a disk with two singular fibers, each of order $2$. In this case, $Y= \mathbb{D}^2(\nicefrac{p_\ast}{q_\ast})$ and $\tau = \emptyset$.
            If $Y$ has one boundary component and $m \ge 1$ singular fibers, then there exists a vertical torus $\Sigma \subset Y$ such that
            \[
                \overline{Y \setminus \Sigma} = C_2 \cup M,
            \]
            where $M= \mathbb{D}^2 (\nicefrac{p_\ast}{q_\ast})$. If $Y$ has two boundary components and no singular fibers, then
            $Y = C_2$ and $\tau = \emptyset$.
            
            Similarly to Lemma \ref{lemma: system for S2}, if $Y$ has $n \ge 2$ boundary components and $m \ge 0$ singular fibers, then there exists a system $\tau \subset Y$ of $n-1+\text{sgm}(m)$ tori such that
        \[
            \overline{Y \setminus \tau} = C_2 \cup
            \begin{dcases}
                \bigcup_{n-1} C_3 & \text{if $m=0$} \\
                \left(\bigcup_{n} C_3\right) \cup\mathbb{D}^2(\nicefrac{p_{\ast}}{q_{\ast}}). & \text{if $m \ge 1$}.
            \end{dcases}
        \]
        For example, Figure \ref{Figure: vertical tori for C2} shows a system of disjoint vertical tori $\tau \subset Y$ in the case that $Y$ has four boundary components and a positive number of singular fibers.
        \end{proof}
\end{lemma}

\begin{cor}\label{cor: the useful decomposition}
        Let $Y$ be a graph manifold rational homology $3$-sphere, there exists a system of embedded disjoint tori $\tau \subset Y$ and three finite (potentially empty) sets $I,J,K$ such that
        \[
            \overline{Y \setminus \tau} = \left(\bigcup_{k \in K} C_3\right) \cup \left(\bigcup_{j \in J} C_2\right) \cup \left(\bigcup_{i \in I} \mathbb{D}^2(\nicefrac{p_{\ast,i}}{q_{\ast,i}})\right),
        \]
        where $\mathbb{D}^2(\nicefrac{p_{\ast,i}}{q_{\ast,i}})$ is a Seifert fibered space with disk base.
        \begin{proof}
        We recall that a solid torus has disk base space and so it is of the form $\mathbb{D}^2(\nicefrac{p_{\ast,i}}{q_{\ast,i}})$.
        Let us suppose that $Y$ is a Seifert fibered space. Since $Y$ a is rational homology $3$-sphere, its base space is either $S^2$ or $\mathbb{RP}^2$. The conclusion holds by removing a fibered solid torus and according to Lemma \ref{lemma: system for S2} and Lemma \ref{lemma: system for RP2}.

        Let us suppose that $Y$ has a nontrivial JSJ decomposition. Let $Y_i$ be a JSJ piece. As we said before, the Seifert fibered space $Y_i$ has a base space that is either a punctured $S^2$ or a punctured $\mathbb{RP}^2$. The conclusion holds by applying Lemma \ref{lemma: system for S2} and Lemma \ref{lemma: system for RP2} to all JSJ pieces $Y_i$.
        \end{proof}
\end{cor}

\begin{cor}\label{cor: C3 obstruction longitudes}
    Let $Y_1$ and $Y_2$ two rational homology solid tori.
    Let $Y=Y_1\cup_{\Sigma_1}C_3 \cup_{\Sigma_2}Y_2$ be such that $\Delta(\lambda_{1},h_1)=\Delta(\lambda_{2},h_2)=0$. The manifold $Y$ is not a rational homology solid torus.
    \begin{proof}
    By Mayer Vietoris, we obtain that the following sequence is exact:
    \[
        \cdots \to H_1(\Sigma_1;\mathbb{Q})\oplus H_1(\Sigma_2;\mathbb{Q}) \overset{\iota}{\to} H_1(Y_1;\mathbb{Q}) \oplus H_1(Y_2;\mathbb{Q}) \oplus H_1(C_3;\mathbb{Q}) \overset{j}{\to} H_1(Y;\mathbb{Q})\to 0.
    \]
    The conditions $\Delta(\lambda_{1},h_1)=0 $ and $\Delta(\lambda_{2},h_2)=0$ imply that $\rank \iota =3$. This implies that $\dim \ker j = \dim \ima \iota=3$. Since $ H_1(Y_1;\mathbb{Q}) \oplus H_1(Y_2;\mathbb{Q}) \oplus H_1(C_3;\mathbb{Q})$ is a vector space of dimension five, we obtain that $\rank j =2$ and $\dim H_1(Y;\mathbb{Q}) =2$. This implies that $Y$ is not a rational homology solid torus.
    \end{proof}
\end{cor}

\begin{cor}\label{cor: C2 obstruction longitudes}
    Let $Y_1$ be a rational homology solid torus.
    If $Y=Y_1\cup_{\Sigma_1'}C_2 $ is such that $\Delta(\lambda_{1},h_1)=0$, then $Y$ is not a rational homology solid torus. 
    \begin{proof}
    By Mayer Vietoris, we obtain that the following sequence is exact:
    \[
        \cdots \to H_1(\Sigma_1';\mathbb{Q}) \overset{\iota}{\to} H_1(Y_1;\mathbb{Q})  \oplus H_1(C_2;\mathbb{Q}) \overset{j}{\to} H_1(Y;\mathbb{Q})\to 0.
    \]
    The condition $\Delta(\lambda_{1},h_1)=0 $ implies that $\rank \iota =1$. Thus, $\dim \ker j = \dim \ima \iota=1$. Since $ H_1(Y_1;\mathbb{Q}) \oplus H_1(C_2;\mathbb{Q})$ is a vector space of dimension three, $\rank j =2$ and $\dim H_1(Y;\mathbb{Q}) =2$. This implies that $Y$ is not a rational homology solid torus.
    \end{proof}
\end{cor}

\section{The recipe for \texorpdfstring{$C_3$}{C3}}

Let $Y=Y_1 \cup_{\Sigma_1}C_3 \cup_{\Sigma_2} Y_2$. In this section we describe $T(Y,\partial Y) \subset \R{\Sigma_3}$ taking as inputs the spaces $T(Y_1,\partial Y_1) \subset \R{\partial Y_1}$ and $T(Y_2,\partial Y_2) \subset \R{\partial Y_2}$. We recall that, if $M$ is a manifold, then
\[
    R(M) \coloneqq \text{Hom} \left( \fund{M};SU(2)\right).
\]
In this section the group $\fund{C_3}$ is considered to be presented as in \eqref{eq: fundamental group C3}.
We also recall that $\partial C_3 = \Sigma_1 \cup \Sigma_2 \cup \Sigma_3$ and that, for $j \in \{1,2,3\}$ the space $\R{\Sigma_j}$ comes with coordinates $(\theta_j,\psi_j)$ as in \eqref{eq: coordinates for C3}. 

In this section we work with coordinates on the torus $S^1 \times S^1 = [0,2\pi]^2/_\sim$. Consequently, the identities are considered modulo $2\pi$. For instance, when in Lemma \ref{lemma: C3 abelian extension} we write $\theta_1+\theta_2=\theta_3$, this has to to be read as an identity modulo $2\pi$.

\begin{lemma}\label{lemma: C3 abelian extension}
      For $j \in \{1,2,3\}$, let $\eta_j \colon \fund{\Sigma_j} \to \Ui$ be the representation with coordinates $(\theta_j,\psi_j)$ in $\R{\Sigma_j}$. There exists an abelian representation $\rho \in R(C_3)$ such that
    \[
        \restr{\rho}{\fund{\Sigma_1}} \equiv \eta_1, \quad \restr{\rho}{\fund{\Sigma_2}} \equiv \eta_2, \quad \text{and} \quad \restr{\rho}{\fund{\Sigma_3}} \equiv \eta_3
    \]
    if and only if $\psi_1=\psi_2=\psi_3$ and $\theta_1+\theta_2=\theta_3$.
    \begin{proof}
        If such an abelian representation $\rho \colon \fund{C_3} \to SU(2)$ exists, then $\rho(h)=\eta_j(h)$ for $j \in \{1,2,3\}$. Hence, $\psi_1=\psi_2=\psi_3$. Moreover, the presentation of $\fund{C_3}$ as in \eqref{eq: fundamental group C3} implies that
        \[
        \eta_1(x_1)\eta_2(x_2)=\rho(x_1)\rho(x_2)=\rho(x_1x_2)=\rho(x_3)=\eta_3(x_3).
        \]
        Hence, $e^{i\theta_1}e^{i\theta_2}=e^{i\theta_3}$ and then $\theta_1+\theta_2=\theta_3$.

        Conversely, let us assume that $\psi_1=\psi_2=\psi_3$ and $\theta_1+\theta_2=\theta_3$. We define $\rho \in R(C_3)$ as
        \[
        \rho(x_1)=e^{i\theta_1}, \quad \rho(x_2)=e^{i\theta_2}, \quad \rho(x_3)=e^{i\theta_3}, \quad \text{and,} \quad \rho(h)=e^{i\psi_1}.
        \]
        The representation $\rho$ is abelian and $\restr{\rho}{\fund{\Sigma_j}} \equiv \eta_j$ for every $j \in \{1,2,3\}$.
    \end{proof}
\end{lemma}

\begin{defn}
    
For two reals $a,b \in [-2,2]$ we define the interval $I(a,b)$ as
\begin{align}\label{eq: interval I}
I(a,b) \coloneqq \left(ab- \frac{1}{2}\left(\sqrt{(4-a^2)(4-b^2)} \right),ab+ \frac{1}{2}\left(\sqrt{(4-a^2)(4-b^2)} \right) \right) \subseteq (-2,2).
\end{align}
\end{defn}
In \cite[Proposition 7.2]{Mino} the author proved that, for given $a,b,c \in [-2,2]$, there exist two matrices $A,B \in SU(2)$ such that $\Tr A =a$, $\Tr B =b$, $\Tr AB =c$ and $AB \neq BA$ if and only if $c \in I(a,b)$. We recall that $\Tr e^{i \theta}= 2\cos \theta$.

\begin{lemma}\label{lemma: C3 extension irreducible}
    If there exists an irreducible representation $\rho \in R(C_3)$, then $\rho(h) \in \ZSU$ and $\Tr \rho(x_3) \in I \left( \Tr\rho(x_1),\Tr \rho(x_2)\right)$.
    
    Conversely, let $\varepsilon \in \{0,\pi\}$ and $\theta_1,\theta_2,\theta_3 \in [0,2\pi]$. If
    \[
    \psi_1=\psi_2=\psi_3=\varepsilon \quad \text{and}  \quad 
    2\cos \theta_3 \in I(2\cos\theta_1,2\cos \theta_2),
    \]
    then there exists an irreducible representation $\rho \in R(C_3)$ such that $\restr{\rho}{\fund{\Sigma_j}}$ is conjugate to the representation $\eta_j\in \R{\Sigma_j}$ with coordinates $(\theta_j,\psi_j)$ for $j \in \{1,2,3\}$. 
    \begin{proof}
        If such an irreducible representation $\rho \in R(C_3)$ exists, then, since $h$ is central in $\fund{C_3}$, we obtain that $\rho(h) \in \ZSU=\{\pm 1\}$ by Fact \ref{Fact: centralizers}.
        Since $\rho$ is irreducible, $\rho(x_1)\rho(x_2) \neq \rho(x_2)\rho(x_1)$, because $\rho$ would be abelian instead.
        Thus, by \cite[Proposition 7.2]{Mino}, $\Tr \rho(x_3) \in I(\Tr \rho(x_1),\Tr \rho(x_2))$.

        We prove now the converse. For $j \in \{1,2,3\}$, let $\eta_j\in \R{\Sigma_j}$ be a representation with coordinates $(\theta_j,\psi_j)$. We recall that $\Tr \eta(x_j) = 2\cos \theta_j$. By \cite[Proposition 7.2]{Mino}, since $2\cos \theta_3 \in I(2\cos\theta_1,2\cos \theta_2)$, there exist three matrices $X_1,X_2,X_3 \in SU(2)$ such that $X_3 = X_1X_2$, $\Tr X_j=2\cos \theta_j$, and $X_1X_2 \neq X_2X_1$. We define the irreducible representation $\rho \colon \fund{C_3} \to SU(2)$ by
        \[
            \rho(x_j)=X_j, \quad \text{and} \quad \rho(h)=e^{i\varepsilon}.
        \]
        Since $\Tr \eta_j(x_j)=\Tr X_j$ and $\Tr \eta_j(h)=\pm 1$, we get that $\restr{\rho}{\fund{\Sigma_j}} \equiv \eta_j$ up to conjugation.
    \end{proof}
\end{lemma}

\begin{defn}\label{defn: A+B}
    Let $\Sigma_1$, $\Sigma_2$, and $\Sigma_3$ be the boundary components of $C_3$. Let $\R{\Sigma_j}$ comes with coordinates $(\theta_j,\psi_j)$ as in \eqref{eq: coordinates for C3}. Let $A \subset \R{\Sigma_1}$ and $B \subset \R{\Sigma_2}$ be two subsets. We define $A+B$ as the subset of $\R{\Sigma_3}$ defined by
    \[
    A+B \coloneqq \left\{ (\theta_1+\theta_2,\psi) \in \R{\Sigma_3} \middle| (\theta_1,\psi) \in \R{\Sigma_1} \text{ and } (\theta_2,\psi) \in \R{\Sigma_2}\right\}.
    \]
\end{defn}
\begin{figure}[t]
            \centering
            \begin{overpic}[width=0.5\textwidth]{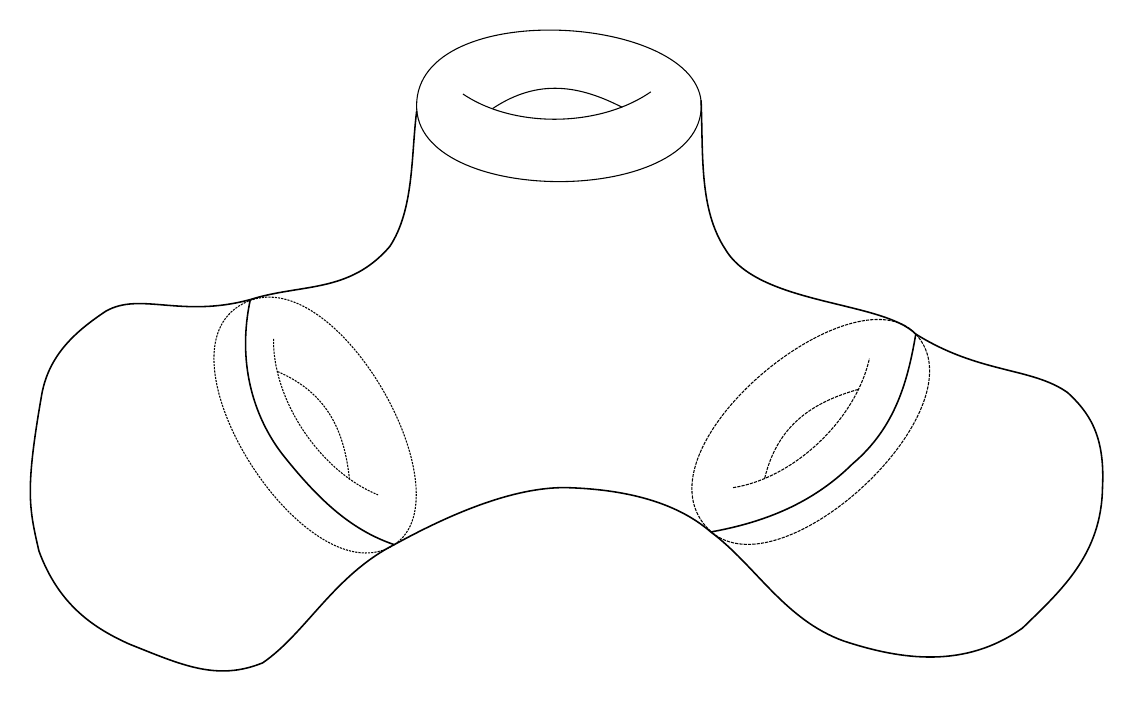}
            \put (12,17){$Y_1$}
            \put (46,31) {$C_3$}
            \put (36,10) {$\Sigma_1$}
            \put (57,12) {$\Sigma_2$}
            \put (83,14) {$Y_2$}
            \end{overpic}
            \caption{The manifold $Y= Y_1 \cup_{\Sigma_1} C_3 \cup_{\Sigma_2} Y_2$.}
            \label{Figure: the manifold Y1 C3 Y2}
\end{figure}
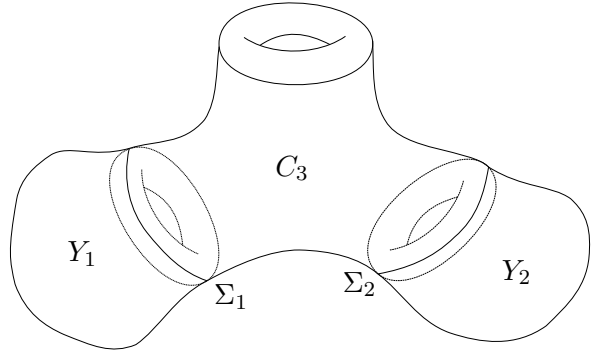
From the rest of the section, $Y_1$ and $Y_2$ will be two manifolds with torus boundary and $Y$ the manifold $Y_1 \cup_{\Sigma_1} C_3 \cup_{\Sigma_2} Y_2$ represented in Figure \ref{Figure: the manifold Y1 C3 Y2}.

\begin{defn}\label{defn H3,0 e H3,pi}
    Let $Y=Y_1 \cup_{\Sigma_1} C_3 \cup_{\Sigma_2} Y_2$ be as before. For $j \in \{1,2\}$, we consider $T(Y_j,\partial Y_j)$ as a subset of $\R{\Sigma_j}$.
    Let $\varepsilon \in \{0,\pi\}$. We define $S_j^\varepsilon $ as $ T(Y_j,\partial Y_j) \cap \{\psi_j = \varepsilon\} \subset \R{\Sigma_j}$. We define the interval $J^\varepsilon \subset (-2,2)$ as
    \[
        J^\varepsilon= \bigcup_{\substack{(\theta_1,\varepsilon) \in S_1^\varepsilon \\ (\theta_2,\varepsilon) \in S_2^\varepsilon}} I(2\cos \theta_1,2\cos\theta_2).
    \]
    Finally, we define $H_{3,\varepsilon} \subseteq \R{\Sigma_3}$ as $H_{3,\varepsilon} \coloneqq \left\{ (\theta_3,\varepsilon) \in \R{\Sigma_3} \middle| 2\cos \theta_3 \in J^\varepsilon \right\}$.
\end{defn}

\begin{rmk}
    \label{rmk: When theta1 theta2 is in piZ}
    According to Definition \ref{defn H3,0 e H3,pi}, the interval $J^\varepsilon$ is an open subinterval of $(-2,2)$. Therefore, it does not contain $\pm 2$. Hence, the sets $H_{3,0}$ and $H_{3,\pi}$ do not contain any of the central representations $\fund{\Sigma_3}\to \ZSU$. Moreover if $2\cos \theta_3 \in I(2\cos \theta_1,2 \cos \theta_2)$, then $\theta_1,\theta_2 \notin \pi \mathbb{Z}$.
\end{rmk}

\begin{teo}\label{teo: alg for C3}
    Let $Y=Y_1 \cup_{\Sigma_1} C_3 \cup_{\Sigma_2} Y_2$. Let $A(Y_j),H(Y_j)$, and $P(Y_j)$ be considered as subsets of $\R{\Sigma_j}$. The space $T(Y,\partial Y)\subset \R{\Sigma_3}$ is equal to $A(Y)\cup H(Y)$, where \[A(Y)=A(Y_1) + A(Y_2)\] and
    \begin{gather*}
        H(Y)= \left(A(Y_1) + H(Y_2)\right) \cup \left(H(Y_1) + A(Y_2)\right) \cup \left(H(Y_1) + H(Y_2)\right)\\ \cup (P(Y_1) + P(Y_2)) \cup H_{3,0} \cup H_{3,\pi}.
     \end{gather*}
    \begin{proof}
        We first prove the statement $A(Y)=A(Y_1)+A(Y_2)$. Let $\eta_3 \in A(Y) \subseteq T(Y,\partial Y)$. Hence, there exists an abelian representation $\rho \in R(Y)$ such that $\restr{\rho}{\fund{\Sigma_3}}\equiv \eta_3$. We define $\eta_1$ (resp. $\eta_2$) as the restriction $ \restr{\rho}{\fund{\Sigma_1}}$ (resp. $ \restr{\rho}{\fund{\Sigma_2}}$). Since $\rho$ is abelian, $\ima \rho \subseteq \Ui$ and therefore $\ima \eta_j \subseteq \Ui$. Thus, $\eta_j \in \R{\Sigma_j}$. Let $(\theta_j,\psi_j)$ be the coordinates of $\eta_j$ in $\R{\Sigma_j}$, Lemma \ref{lemma: C3 abelian extension} implies that $\psi_1=\psi_2=\psi_3$ and $\theta_3 = \theta_1 + \theta_2$. Thus, $(\theta_3,\psi_3) \in A(Y_1)+A(Y_2)$. This implies that $A(Y) \subseteq A(Y_1)+A(Y_2)$.

        Conversely, let $\eta_3 \in A(Y_1)+A(Y_2)$ be of coordinates $(\theta_3,\psi_3)$. By the definition of the set $A(Y_1) + A(Y_2)$, we obtain that $\theta_3=\theta_1+\theta_2$ with $(\theta_1,\psi_3) \in A(Y_1) \subset \R{\Sigma_1}$ and $(\theta_2,\psi_3) \in A(Y_2) \subset \R{\Sigma_2}$. Let $\eta_1$ and $\eta_2$ be the representations corresponding to these two points.
        Lemma \ref{lemma: C3 abelian extension} implies that there exists an abelian representation $\rho_c \in R(C_3)$ such that $\restr{\rho_c}{\fund{\Sigma_j}} \equiv \eta_j$.
        Let $\rho_1 \in R(Y_1)$ (resp. $\rho_2 \in R(Y_2)$) be an abelian representation such that $\restr{\rho_1}{\fund{\Sigma_1}} \equiv \eta_1$ (resp. $\restr{\rho_2}{\fund{\Sigma_2}} \equiv \eta_2$). Thus, there exists an abelian representation $\rho \in R(Y)$ such that
        \[
            \restr{\rho}{\fund{Y_1}} \equiv \rho_1, \quad \restr{\rho}{\fund{Y_2}} \equiv \rho_2 \quad \text{and} \quad \restr{\rho}{\fund{C_3}}\equiv \rho_c,
        \]
        In particular, $\restr{\rho}{\fund{\Sigma_3}}\equiv\restr{\rho_c}{\fund{\Sigma_3}}\equiv\eta_3$. Thus, $\eta_3 \in A(Y)$ and $A(Y_1)+A(Y_2) \subseteq A(Y)$.

        We focus on $H(Y)$.
        If $\eta_3 \in P(Y_1) + P(Y_2)$, then, by Lemma \ref{lemma: C3 abelian extension}, there exists an abelian representation $\rho_c \in R(C_3)$ such that $\restr{\rho_c}{\fund{\Sigma_1}}$ and $\restr{\rho_c}{\fund{\Sigma_2}}$ are central and they extend to two abelian non-central representations $\rho_1 \colon\fund{Y_1} \to SU(2)$ and $\rho_2 \colon \fund{Y_2} \to SU(2)$, respectively. Hence, there exists a representation $\rho \in R(Y)$ such that
        \[
            \restr{\rho}{\fund{Y_1}} \equiv \rho_1, \quad \restr{\rho}{\fund{Y_2}} \equiv \rho_2, \quad \text{and} \quad \restr{\rho}{\fund{C_3}}\equiv \rho_c,
        \]
        If $\rho$ is irreducible, then $\eta_3 \in H(Y)$ and $P(Y_1) + P(Y_2) \subseteq H(Y)$. Let us assume that $\rho$ is abelian, this implies that $\ima \rho$, $\ima \rho_1$, $\ima \rho_2$, and $\ima \rho_c$ are all in $\Ui$.
        We define
        \[
            M \coloneq \frac{\sqrt{2}}{2} \begin{bmatrix}
                i & i \\ i & -i
            \end{bmatrix} \in SU(2).
        \]
        Let us compute:
        \[
            Me^{i x}M^{-1} = 
            \begin{bmatrix}
                \frac{e^{i x}+e^{-i x}}{2}  & \frac{e^{ix} -e^{-ix}}{2}\\ \frac{e^{i x} - e^{-i x}}{2}  & \frac{e^{i x} +e^{-ix}}{2}
            \end{bmatrix} = \begin{bmatrix}
                \cos x  & i \sin x \\ i \sin x  & \cos x
            \end{bmatrix}.
        \]
        Notice that if $e^{i x} \in SU(2)$ is non-central, meaning that $x \notin \{0,\pi\}$, then the matrix $Me^{i x}M^{-1}$ is non-diagonal.
        Therefore, the product $Me^{i x}M^{-1}$ does not commute with $e^{i x}$. 
        We define the representation $\rho' \in R(Y)$ as
        \[
            \rho'(x)= \begin{dcases}
                M \rho(x)M^{-1} \quad \text{if } x \in \fund{Y_1} \\
                \rho(x) \quad \text{otherwise}.
            \end{dcases}
        \]
        We notice that, since $\rho(\fund{\Sigma_1})\subset \ZSU$, the representation $\rho'$ is well-defined.
        Since neither $\rho_1$ nor $\rho_2$ is central, $\rho'$ is irreducible. Hence $\eta_3 \in H(Y)$ and $P(Y_1) + P(Y_2) \subseteq H(Y)$.
        
        If $\eta_3 \in \left(A(Y_1) + H(Y_2)\right) \cup \left(H(Y_1) + A(Y_2)\right) \cup \left(H(Y_1) + H(Y_2)\right)$, then by Lemma \ref{lemma: C3 abelian extension}, there exists an abelian representation $\rho_c \in R(Y)$ such that $\restr{\rho_c}{\fund{\Sigma_3}}=\eta_3$ and either $\restr{\rho}{\fund{\Sigma_1}} \in H_1$ or $\restr{\rho}{\fund{\Sigma_2}} \in H_2$. This implies that there exists a $\rho \in R(Y)$ such that $\restr{\rho}{\fund{\Sigma_3}}=\eta_3$ and either $\restr{\rho}{\fund{Y_1}}$ or $\restr{\rho}{\fund{Y_2}}$ is irreducible. Thus, $\rho$ is irreducible and $\eta_3 \in H(Y)$. This implies that $\left(A(Y_1) + H(Y_2)\right) \cup \left(H(Y_1) + A(Y_2)\right) \cup \left(H(Y_1) + H(Y_2)\right) \subseteq H(Y)$.

        Let $\varepsilon \in \{0,\pi\}$ and $\eta_3 \in \eta_3 \in H_{3,\varepsilon}$. Let $(\theta_3,\varepsilon)$ be the coordinates of $\eta_3$ in $\R{\Sigma_3}$. By definition of the set $H_{3,\varepsilon}$, there exist $\theta_1$ and $\theta_2$ such that 
        \begin{align}
            \label{eq: la uso una volta coseni}
            2\cos \theta_3 \in I(2\cos \theta_1,2\cos \theta_2), \quad (\theta_1,\varepsilon) \in T(Y_1,\partial Y_1), \quad \text{and} \quad (\theta_2,\varepsilon) \in T(Y_2,\partial Y_2).
        \end{align}
        Let $\eta_1 \in \R{\Sigma_1}$ and $\eta_2 \in \R{\Sigma_2}$ be the representations corresponding to the coordinates $(\theta_1,\varepsilon)$ and $(\theta_2,\varepsilon)$ respectively.
        The identity \eqref{eq: la uso una volta coseni} states that $\eta_1$ and $\eta_2$ extend to $\rho_1 \colon \fund{Y_1} \to SU(2)$ and $\rho_2 \colon \fund{Y_2} \to SU(2)$. In particular, for $j \in \{1,2\}$ we have that $\restr{\rho_j}{\fund{\Sigma_j}} \equiv \eta_j$.
        
        Lemma \ref{lemma: C3 extension irreducible} implies that there exists an irreducible representation $\rho_c \in R(C_3)$ and three matrices $X_1,X_2,X_3 \in SU(2)$ such that
        \[
            \restr{\rho_c}{\fund{\Sigma_j}}=X_j \eta_j X_j^{-1}.
        \]
        Up to conjugation, we can assume that $\restr{\rho_c}{\fund{\Sigma_3}}\equiv \eta_3$. We define the representation $\rho \colon \fund{Y} \to SU(2)$ as
        \[
        \rho(x)= 
        \begin{dcases}
            X_1 \rho_1(x) X_1^{-1} \quad \text{if $x \in \fund{Y_1}$} \\
            X_2 \rho_2(x) X_2^{-1} \quad \text{if $x \in \fund{Y_2}$} \\
            \rho_c(x) \quad \text{if $x \in \fund{C_3}$}.
        \end{dcases}
        \]
        We notice that, for $j \in \{1,2\}$, we have that $\restr{\rho_c}{\fund{\Sigma_j}} \equiv X_j \eta_jX_j^{-1}$ and
        \[
            \restr{\left(X_j\rho_j X_j^{-1}\right)}{\fund{\Sigma_j}} \equiv X_j \restr{\rho_j}{\fund{\Sigma_j}}X_j^{-1} \equiv X_j \eta_j X_j^{-1} \equiv \restr{\rho_c}{\fund{\Sigma_j}}.
        \]
        This implies that the representation $\rho$ is well-defined.
        The representation $\rho$ is irreducible and $\restr{\rho}{\fund{\Sigma_3}}\equiv \eta_3$. Thus, $\eta_3 \in H(Y)$ and hence $H_{3,0} \cup H_{3,\pi} \subset H(Y)$. This concludes the proof that
        \begin{multline*}
            \left(A(Y_1) + H(Y_2)\right) \cup \left(H(Y_1) + A(Y_2)\right) \cup \left(H(Y_1) + H(Y_2)\right) \cup (P(Y_1) + P(Y_2)) \cup \\ \cup H_{3,0} \cup H_{3,\pi} \subseteq H(Y).
        \end{multline*}

        Conversely, let $\eta_3 \in H(Y)$. Thus, there exists an irreducible extension $\rho \in R(Y)$ of $\eta_3$.
        Let us focus on $\rho_c \coloneq \restr{\rho}{\fund{C_3}}$.
        If $\rho_c$ is irreducible, then Lemma \ref{lemma: C3 extension irreducible} implies that $\eta_3 \in H_{3,0}\cup H_{3,\pi}$.
        We conclude that
        \[
            H_{3,0} \cup H_{3,\pi} \subseteq H(Y).
        \]
        
        Let us assume that $\rho_c$ is abelian. If $\restr{\rho}{\fund{Y_1}}$ is irreducible, then Lemma \ref{lemma: C3 abelian extension} implies that $\eta_3 \in (H(Y_1)) + A(Y_2)) \cup (H(Y_1) + H(Y_2))$. If $\restr{\rho}{\fund{Y_1}}$ is abelian and $\restr{\rho}{\fund{Y_2}}$ is irreducible, then Lemma \ref{lemma: C3 abelian extension} implies that $\eta_3 \in (A(Y_1) + H(Y_2))$. Finally, if both $\restr{\rho}{\fund{Y_1}}$ and $\restr{\rho}{\fund{Y_2}}$ are both abelian, then, as we saw before, since $\rho$ is irreducible, $\eta_3 \in P(Y_1) + P(Y_2)$. This concludes that
        \begin{multline*}
       H(Y) \subseteq\\ \left(A(Y_1) + H(Y_2)\right) \cup \left(H(Y_1) + A(Y_2)\right) \cup \left(H(Y_1) + H(Y_2)\right) \cup (P(Y_1) + P(Y_2)) \cup H_{3,0} \cup H_{3,\pi}.
        \end{multline*}
    \end{proof}
\end{teo}

\begin{figure}[t]
    \centering
    \begin{subfigure}[b]{0.4\textwidth}
        \begin{tikzpicture}[>=latex,scale=1.45]
            \draw[->] (-.5,0) -- (2.5,0) node[right] {$\theta_1$};
            \foreach \x /\n in {1/$\pi$,2/$2\pi$} \draw[shift={(\x,0)}] (0pt,2pt) -- (0pt,-2pt) node[below] {\tiny \n};
            \draw[->] (0,-.5) -- (0,2.5) node[below right] {$\psi_1$};
            \foreach \y /\n in {1/$\pi$,2/$2\pi$}
            \draw[shift={(0,\y)}] (2pt,0pt) -- (-2pt,0pt) node[left] {\tiny \n};
            \node[below left] at (0,0) {\tiny $0$};
             \draw (0,0) rectangle (2,2);
             \draw (1,0) -- (1,2); 
             \draw (0,1) -- (2,1);
             \draw [blue, thick] (0,0)--(2,2);
             \draw [red, thick] (0.5,0.5)--(0,2);
             \draw [red, thick] (1.5,1.5)--(2,0);
             \filldraw[cyan] (1,1) circle (1pt);
        \end{tikzpicture}
        \caption{$T(Y_1,\partial Y_1) \subset \R{\Sigma_1}$.}
    \end{subfigure}
    \quad
    \begin{subfigure}[b]{0.4\textwidth}
        \begin{tikzpicture}[>=latex,scale=1.45]
        \draw[->] (-.5,0) -- (2.5,0) node[right] {$\theta_2$};
        \foreach \x /\n in {1/$\pi$,2/$2\pi$} \draw[shift={(\x,0)}] (0pt,2pt) -- (0pt,-2pt) node[below] {\tiny \n};
        \draw[->] (0,-.5) -- (0,2.5) node[below right] {$\psi_2$};
        \foreach \y /\n in {1/$\pi$,2/$2\pi$}
        \draw[shift={(0,\y)}] (2pt,0pt) -- (-2pt,0pt) node[left] {\tiny \n};
        \node[below left] at (0,0) {\tiny $0$};
         \draw (0,0) rectangle (2,2);
         \draw (1,0) -- (1,2); 
         \draw (0,1) -- (2,1);
            \draw [blue, thick] (1,0)--(1,2);
            \draw [red, thick] (1,1.5)--(0,1);
            \draw [red, thick] (1,0.5)--(2,1);
            \filldraw[cyan, thick] (1,2) circle (1pt);
            \filldraw[cyan, thick] (1,1) circle (1pt);
            \filldraw[cyan, thick] (1,0) circle (1pt);
    \end{tikzpicture}
        \caption{$T(Y_2,\partial Y_2) \subset \R{\Sigma_2}$.}
    \end{subfigure} \\[10pt]
    \begin{subfigure}[t]{0.4\textwidth}
        \begin{tikzpicture}[>=latex,scale=1.45]
            \draw[->] (-.5,0) -- (2.5,0) node[right] {$\theta_3$};
            \foreach \x /\n in {1/$\pi$,2/$2\pi$} \draw[shift={(\x,0)}] (0pt,2pt) -- (0pt,-2pt) node[below] {\tiny \n};
            \draw[->] (0,-.5) -- (0,2.5) node[below right] {$\psi_3$};
            \foreach \y /\n in {1/$\pi$,2/$2\pi$}
            \draw[shift={(0,\y)}] (2pt,0pt) -- (-2pt,0pt) node[left] {\tiny \n};
            \node[below left] at (0,0) {\tiny $0$};
             \draw (0,0) rectangle (2,2);
             \draw (1,0) -- (1,2); 
             \draw (0,1) -- (2,1);
             \draw[blue, thick] (1,0) --(2,1);
             \draw[blue, thick] (0,1)-- (1,2);

             \draw[orange, thick] (0.5,1.5)--(1,0);
             \draw[orange, thick] (1.5,0.5)--(1,2);

             \draw[magenta, thick] (0,1-1/3) --(2,1+1/3);
             \draw[magenta, thick] (0,1+1/3) --(0.5,1.5);
             \draw[magenta, thick] (2,1-1/3) --(1.5,0.5);

             \draw[red, thick] (0,1.2)--(0.5,1.5);
             \draw[red, thick] (2,1.2)--(5/6,0.5);
             
             \draw[red, thick] (7/6,1.5)--(0,0.8);
             \draw[red, thick] (2,0.8)--(1.5,0.5);
             \filldraw[cyan] (2,1) circle (1pt);
             \filldraw[cyan] (0,1) circle (1pt);
             
    \filldraw[cyan] (3,0.2) circle (2pt) node[anchor=west]{$\;=P_1+P_2$};
    \filldraw[blue] (3,0.6) circle (2pt) node[anchor=west]{$\;=A_1+A_2$};
    \filldraw[orange] (3,1) circle (2pt) node[anchor=west]{$\;=H_1+A_2$};
    \filldraw[magenta] (3,1.4) circle (2pt) node[anchor=west]{$\;=A_1+H_2$};
    \filldraw[red] (3,1.8) circle (2pt) node[anchor=west]{$\;=H_1+H_2$};
        \end{tikzpicture}
        \caption{$T(Y,\partial Y) \subset \R{\Sigma_3}$.}
    \end{subfigure}
    \caption{An application of Theorem \ref{teo: alg for C3}. In the first two pictures, $A(Y_1)$ and $A(Y_2)$ are in blue, $H(Y_1)$ and $H(Y_2)$ are in red, and the sets $P(Y_1)$ and $P(Y_2)$ are in light blue.}
    \label{Figure: example of alg C3}
\end{figure}
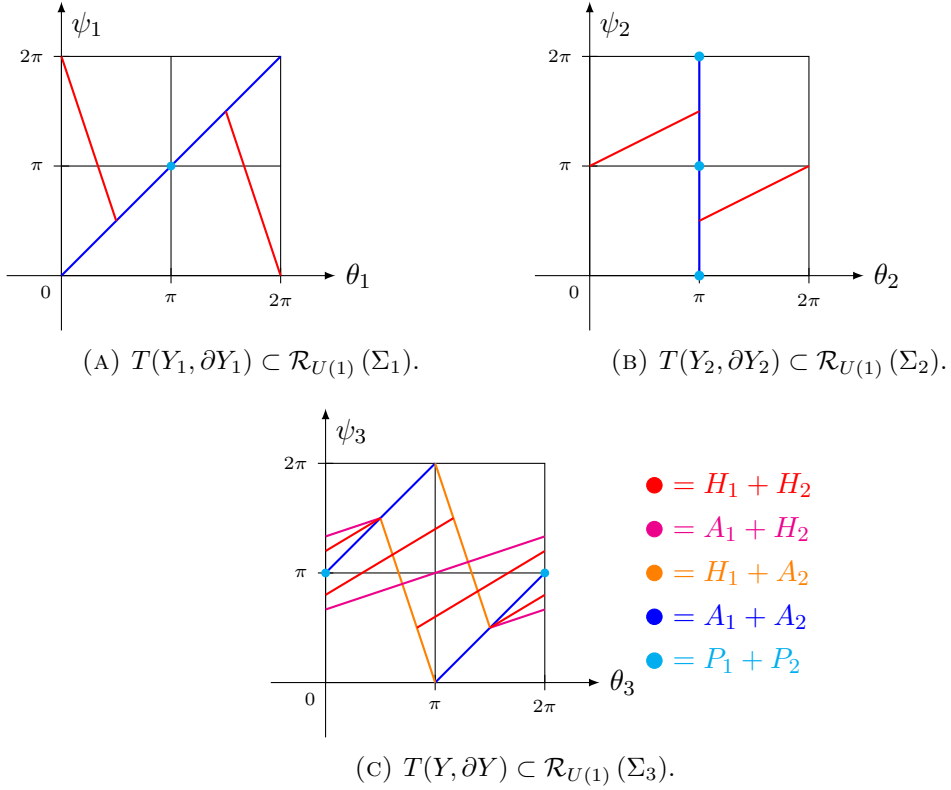

We give Figure \ref{Figure: example of alg C3} as an example how Theorem \ref{teo: alg for C3} works. In the considered case, we notice that $S_1^0=\{(0,0)\}$ and $S_2^\pi=\{(0,\pi),(\pi,\pi)\}$, thus $J^0$ and $J^\pi$ are empty, and therefore $H_{3,0}$ and $H_{3,\pi}$ are empty as well.

\section{The recipe for \texorpdfstring{$C_2$}{C2}}

Let $Y=Y_1 \cup_{\Sigma_1'}C_2 $. In this section we describe $T(Y,\partial Y) \subset \R{\Sigma_2'}$ taking as input the space $T(Y_1,\partial Y_1) \subset \R{\partial Y_1'}$. In this section the group $\fund{C_2}$ is considered to be presented as in \eqref{eq: fundamental group C2}.

We recall that $\partial C_2 = \Sigma_1' \cup \Sigma_2' $ and, in particular, the space $\R{\Sigma_j'}$ comes with coordinates $(\theta_j,\psi_j)$ as in \eqref{eq: coordinates for C2}. We recall that, as explained in Section \ref{sec: two important pieces}, the manifold $C_2$ is a Seifert fibered space with regular fiber $h \in \fund{C_2}$. For $j \in \{1,2\}$ the torus $\Sigma_j'$ contains a regular fiber, which we called $h_j \in \fund{\Sigma_j'}$. In particular, $h_1=h=h_2$ in $\fund{C_2}$.

As in the previous section, we will work with coordinates on the torus $S^1 \times S^1 = [0,2\pi]^2/_\sim$. Consequently, the identities are considered modulo $2\pi$.

\begin{lemma}{\label{lemma: extension abelian C2}}
    The representation $\eta_1 \colon \fund{\Sigma_1'}\to SU(2)$ extends to an abelian representation of $\fund{C_2}$ if and only if $\eta_1(h_1)^2=1$.
    \begin{proof}
        Let us assume that there exists an abelian representation $\rho \in R(C_2)$ such that $\restr{\rho}{\fund{\Sigma_1'}}\equiv\eta_1$. Since $\rho$ is abelian, it factors through the first homology $H_1(C_2;\mathbb{Z})$. Notice that the presentation \eqref{eq: fundamental group C2} implies that the homology class of the regular fiber $h$ has order $2$ in $H_1(C_2;\mathbb{Z})$. Hence, $\rho(h)^2=\eta_1(h_1)^2=1$.

        Conversely, up to conjugation, we assume that $\ima \eta_1 \subseteq \Ui$. Let $(\theta_1,\psi_1)$ be the coordinates of $\eta_1$ in $\R{\Sigma_1'}$. Since $\eta_1(h)^2=1$, we obtain that $\psi_1 \in \pi \mathbb{Z}$. We define the representation $\rho \in R(C_2)$ as
        \[
        \rho(x_1)=e^{i \theta_1}, \quad \rho(x_2)=e^{i (\pi -\theta_1)}, \quad \rho(h)=e^{i\psi_1}, \quad \text{and} \quad \rho(z)=e^{i \frac{\pi}{2}}.
        \]
        The representation $\rho$ is abelian and it is an extension of $\eta_1$.
    \end{proof}
\end{lemma}

\begin{rmk}\label{rmk: the matrices H and K for C2}
    We define $Z \in SU(2)$ as
    \[
         Z \coloneq \frac{\sqrt{2}}{2} 
        \begin{bmatrix}
            0 & 1+i \\
            -1+i & 0
        \end{bmatrix}.
    \]
    Let $e^{i\alpha} \in SU(2)$ be a non-central matrix.
    Then, $Z^2=-1$, $Ze^{i\alpha}Z^{-1}e^{i\alpha}=1$ and $Z e^{i\alpha}\neq e^{i\alpha}Z$.
\end{rmk}

\begin{lemma}\label{lemma: irr extension C2, translazione}
    Let $\eta_1 \in \R{\Sigma_1'}$ be the representation of coordinates $(\theta_1,\psi_1)$ with $\psi_1 \notin \pi \mathbb{Z}$. There exists an irreducible representation $\rho \colon \fund{C_2} \to SU(2)$ such that $\restr{\rho}{\fund{\Sigma_1'}} \equiv \eta_1$, $\ima \restr{\rho}{\fund{\Sigma_2'}} \subseteq \Ui$. Furthermore, $\restr{\rho}{\fund{\Sigma_2'}}$ has coordinates $(\pi -\theta_1,\psi_1) \in \R{\Sigma_2}$.
    \begin{proof}
        Let $Z \in SU(2)$ be the matrix in Remark \ref{rmk: the matrices H and K for C2}. We define the representation $\rho \in R(C_2)$ as
        \[
        \rho(x_1)=e^{i\theta_1} \quad \rho(x_2)=e^{i(\pi-\theta_1)} \quad \rho(h)=e^{i\psi_1}, \quad \text{and} \quad \rho(z)=Z.
        \]
        Since $\psi_1 \notin \pi \mathbb{Z}$, the element $\rho(h)$ is non-central and Remark \ref{rmk: the matrices H and K for C2} implies that $\rho(z)\rho(h)\rho(z)^{-1} \rho(h)=1$ and $\rho(z)\rho(h) \neq \rho(h)\rho(z)$. Therefore, $\rho$ is a well defined representation of $\fund{C_2}$ and is irreducible.
        The representation $\rho$ is an extension of $\eta_1$. Moreover, $\restr{\rho}{\fund{\Sigma_2'}}$ has image in $\Ui$ and it has coordinates $(\pi - \theta_1,\psi_1)$ in $\R{\Sigma_2}$.
    \end{proof}
\end{lemma}

\begin{lemma}\label{lemma: irr extension C2, central h}
     Let $\eta_1 \in \R{\Sigma_1'}$ be of coordinates $(\theta_1,\psi_1)$. Let us assume that $\psi_1 \in \pi \mathbb{Z}$ and $\theta_1 \notin \pi \mathbb{Z}$. For every $\theta_2 \in (0,\pi)\cup (\pi,2\pi)$ there exists an irreducible representation $\rho \colon \fund{C_2} \to SU(2)$ such that $\restr{\rho}{\fund{\Sigma_1'}}\equiv \eta_1$ and $\Tr \rho(x_2) = 2\cos \theta_2$.
     \begin{proof}
        We define the representation $\rho \colon \fund{C_2} \to SU(2)$ as
        \[
        \rho(x_1)=e^{i\theta_1} \quad \rho(x_2)=\begin{bmatrix}
            \cos \theta_2 & i\sin\theta_2 \\ i\sin\theta_2 & \cos \theta_2
        \end{bmatrix} \quad \rho(h)=e^{i\psi_1}, \quad \text{and} \quad \rho(z)=Z',
        \]
        where $Z'\in SU(2)$ is such that $(Z')^2 =\rho(x_2)^{-1}\rho(x_1)^{-1}$. Such a matrix $Z'$ exists as the product $\rho(x_2)^{-1}\rho(x_1)^{-1}$ is contained in its centralizer, which is a divisible subgroup of $SU(2)$. Details can be found in \cite[Proposition 3.19]{Lam_1999}.
        Since $\psi_1 \in \pi \mathbb{Z}$, the image $\rho(h)$ is central in $SU(2)$. This implies that $[\rho(x_j),\rho(h)]=1$ for $j \in \{1,2\}$.
        Since $\theta_2 \notin \pi \mathbb{Z}$, the matrix $\rho(x_2)$ is not diagonal. This implies that $\rho(x_1)\rho(x_2)\neq \rho(x_2)\rho(x_1)$, and hence $\rho$ is irreducible. We conclude that $\restr{\rho}{\fund{\Sigma_1'}}\equiv \eta_1$ and $\Tr \rho(x_2) = 2\cos \theta_2$.
     \end{proof}
\end{lemma}

\begin{lemma}\label{lemma: almost the receipt for C2}
     Let $Y=Y_1 \cup_{\Sigma_1'} C_2$. Let the set $T(Y_1,\partial Y_1) =A(Y_1) \cup H(Y_1)$ be considered as a subset of $\R{\Sigma_1'}$. Let $(\theta_1,\psi_1) \in T(Y_1,\partial Y_1)$ and $\varepsilon \in \{0,\pi\}$, then
     \begin{enumerate}
         \item If $\psi_1 \notin \pi \mathbb{Z}$, then $(-\theta_1+\pi,\psi_1) \in H(Y) \subset \R{\Sigma_2'}$
          \item If $\psi_1 = \varepsilon$ and $(\theta_1,\psi_1) \in A(Y_1)$, then $\{\psi_2 = \varepsilon\} \subset A(Y) \subset \R{\Sigma_2'}$,
          \item If $\psi_1 = \varepsilon$ and $(\theta_1,\psi_1) \in H(Y_1)$, then $\{\psi_2 = \varepsilon\} \subset H(Y) \subset \R{\Sigma_2'}$,
          \item If $\psi_1 = \varepsilon$ and $\theta_1 \notin \pi \mathbb{Z}$, then $\{\psi_2 = \varepsilon\}\setminus \{0,\pi\}^2 \subset H(Y) \subset \R{\Sigma_2'}$.
     \end{enumerate}
     \begin{proof}
     Let $\eta_1 \colon \fund{\Sigma_1'}\to SU(2)$ be the representation corresponding to $(\theta_1,\psi_1)\in \R{\Sigma_1'}$.
     As $(\theta_1,\psi_1) \in T(Y_1,\partial Y_1)$, there exists a representation $\rho_1 \colon \fund{Y_1} \to SU(2)$ such that $\restr{\rho_1}{ \fund{\Sigma_1'}}$ has these coordinates.
     
     Let us assume that $\psi_1 \notin \pi \mathbb{Z}$. Lemma \ref{lemma: irr extension C2, translazione} implies that there exists an irreducible representation $\rho_c \colon\fund{C_2} \to SU(2)$ such that
     $\restr{\rho_c}{\fund{\Sigma_1'}} \in \R{\Sigma_1'}$ has coordinates $(\theta_1,\psi_1)$ and $\restr{\rho_c}{\fund{\Sigma_2'}}$ has coordinates $(-\theta_1+\pi,\psi_1)$.
     Therefore, there exists an irreducible representation $\rho \colon \fund{Y} \to SU(2)$ such that
     \[
        \restr{\rho}{\fund{Y_1}}\equiv \rho_1 \quad \text{and} \quad \restr{\rho}{\fund{C_2}}\equiv \rho_c.
     \]
     This proves the first statement. Similarly,
     the second and the third are consequences of Lemma \ref{lemma: extension abelian C2}.
    Let us focus on the fourth one. Let $\theta_2 \in (0,\pi) \cup (\pi,2\pi)$ and $\psi_1=\varepsilon$. 
    As a consequence of Lemma \ref{lemma: irr extension C2, central h}, there exists an irreducible representation $\rho_c \colon \fund{C_2} \to SU(2)$ such that $\restr{\rho_c}{\fund{\Sigma_1}'}$ has coordinates $(\theta_1,\varepsilon)$ and $\Tr \rho_c(x_2) = 2\cos \theta_2$. Therefore, there exists an irreducible representation $\rho \colon \fund{Y} \to SU(2)$ such that
    \[
        \restr{\rho}{\fund{Y_1}} \equiv \rho_1, \quad \restr{\rho}{\fund{C_2}}\equiv \rho_c, \quad \text{and}\quad \Tr \rho(x_2) = 2\cos \theta_2.
    \]
    Up to conjugation, we can assume that $\restr{\rho}{\fund{\Sigma_2'}}$ has image in $\Ui$. In particular, $\rho(x_2)=e^{i \theta_2}$ up to conjugation. This implies that $(\theta_2,\varepsilon) \in H(Y) \subset \R{\Sigma_2'}$. We conclude that, for every $\theta_2 \in (0,\pi) \cup (\pi,2\pi)$, the point $(\theta_2,\varepsilon)$ is in $ H(Y) \subset \R{\Sigma_2'}$. This concludes the fourth statement.
     \end{proof}
\end{lemma}

\begin{figure}[t]
    \centering
    \begin{subfigure}[b]{0.4\textwidth}
        \begin{tikzpicture}[>=latex,scale=1.55]
            \draw[->] (-.5,0) -- (2.5,0) node[right] {$\theta_1$};
            \foreach \x /\n in {1/$\pi$,2/$2\pi$} \draw[shift={(\x,0)}] (0pt,2pt) -- (0pt,-2pt) node[below] {\tiny \n};
            \draw[->] (0,-.5) -- (0,2.5) node[below right] {$\psi_1$};
            \foreach \y /\n in {1/$\pi$,2/$2\pi$}
            \draw[shift={(0,\y)}] (2pt,0pt) -- (-2pt,0pt) node[left] {\tiny \n};
            \node[below left] at (0,0) {\tiny $0$};
             \draw (0,0) rectangle (2,2);
             \draw (1,0) -- (1,2); 
             \draw (0,1) -- (2,1);
             \draw [blue, thick] (0,0)--(2,2);
             \draw [red, thick] (0.5,0.5)--(0,2);
             \draw [red, thick] (1.5,1.5)--(2,0);
             
        \end{tikzpicture}
        \caption{$T(Y_1,\partial Y_1) \subset \R{\Sigma_1}$.}
    \end{subfigure}
    \qquad
    \begin{subfigure}[b]{0.4\textwidth}
        \begin{tikzpicture}[>=latex,scale=1.55]
        \draw[->] (-.5,0) -- (2.5,0) node[right] {$\theta_2$};
        \foreach \x /\n in {1/$\pi$,2/$2\pi$} \draw[shift={(\x,0)}] (0pt,2pt) -- (0pt,-2pt) node[below] {\tiny \n};
        \draw[->] (0,-.5) -- (0,2.5) node[below right] {$\psi_2$};
        \foreach \y /\n in {1/$\pi$,2/$2\pi$}
        \draw[shift={(0,\y)}] (2pt,0pt) -- (-2pt,0pt) node[left] {\tiny \n};
        \node[below left] at (0,0) {\tiny $0$};
         \draw (0,0) rectangle (2,2);
         \draw (1,0) -- (1,2); 
         \draw (0,1) -- (2,1);
            \draw[blue, thick] (0,0.01)--(2,0.01);
            \draw[blue, thick] (0,1.01)--(2,1.01);
            \draw[blue, thick] (0,2.01)--(2,2.01);
            \draw[red, thick] (0,1-0.01)--(2,1-0.01);
            \draw[red, thick] (0,2-0.01)--(2,2-0.01);
            \draw[red, thick] (0,0-0.01)--(2,0-0.01);
            
            \draw[red, thick] (1,0)--(0,1);
            \draw[red, thick] (2,1)--(1,2);
            \draw[red, thick] (0.5,0.5)--(1,2);
            
            
    \end{tikzpicture}
        \caption{$T(Y,\partial Y) \subset \R{\Sigma_2}$.}
    \end{subfigure}
    \caption{Example of Theorem \ref{teo: alg for C2}. In the pictures, $A(Y_1)$ and $A(Y)$ are in blue, $H(Y_1)$ and $H(Y)$ are in red.}
    \label{Figure example C2}
\end{figure}
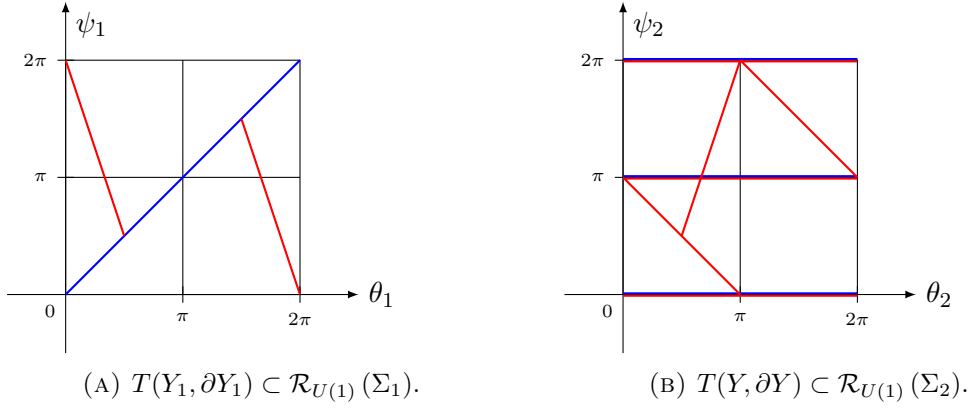

Here is an observation on the next theorem: let $Y=Y_1 \cup_{\Sigma_1'} C_2$, the trivial representation $t\colon \fund{Y_1} \to SU(2)$ is abelian and $\restr{t}{\fund{\Sigma_1}}$ has coordinates $(0,0)\in \R{\Sigma_1'}$, hence $\restr{t}{\Sigma_1} \in A(Y_1)$. This implies that $A(Y_1) \cap \{\psi_1=0\} \subset \R{\Sigma_1'}$ is never empty. We give Figure \ref{Figure example C2} as an example of how Theorem \ref{teo: alg for C2} works.

\begin{teo}\label{teo: alg for C2}
    Let $Y=Y_1 \cup_{\Sigma_1} C_2$. Let the set $T(Y_1,\partial Y_1) =A(Y_1) \cup H(Y_1)$ be considered as a subset of $\R{\Sigma_1'}$. As a subset of $\R{\Sigma_2'}$, the set $T(Y,\partial Y)$ is equal to $A(Y) \cup H(Y)$. Furthermore,
        \begin{align}\label{eq: la uso una volta. A(Y) in C2 alg}
        A(Y) = \begin{dcases}
            \{\psi_2 = 0\} & \text{if } A(Y_1) \cap \{\psi_1=\pi\} = \emptyset \\
            \{\psi_2 = 0\} \cup \{\psi_2=\pi\} & \text{otherwise}.
        \end{dcases}
         \end{align}
        and
        $H_{Y,t}  \cup H_{Y,0} \cup H_{Y,\pi} \subseteq H(Y)$, where
        \begin{align*}
             H_{Y,t} \coloneq \left\{(- \theta_1+\pi,\psi_1) \in \R{\Sigma_2'} \middle| (\theta_1,\psi_1) \in T(Y_1,\partial Y_1), \psi_1 \notin \pi \mathbb{Z}\right\},
            \end{align*}
and, for $\varepsilon \in \{0,\pi\}$,
            \begin{align*}
             H_{Y,\varepsilon} \coloneq \begin{dcases}
                    \{\psi_2 = \varepsilon\} \subset \R{\Sigma_2'}& \text{if } H(Y_1) \cap \{\psi_1=\varepsilon\} \neq \emptyset \\
                    \{\psi_2 = \varepsilon \}\setminus \{0,\pi\}^2 \subset \R{\Sigma_2'} & \text{if } T(Y_1,\partial Y_1) \cap \left( \{\psi_1=\varepsilon\} \setminus \{0,\pi\}^2\right) \neq \emptyset \\
           \emptyset & \text{otherwise}.
            \end{dcases}
        \end{align*}
        \begin{proof}
            Let us start by showing the conclusion for $A(Y)$. Let $S$ be the right-hand side of \eqref{eq: la uso una volta. A(Y) in C2 alg}. The second statement of Lemma \ref{lemma: almost the receipt for C2} implies that $S \subseteq A(Y)$. Conversely, let $\eta \in A(Y) \subset \R{\Sigma_2'}$. Thus, there exists an abelian representation $\rho \colon \fund{Y} \to SU(2)$ such that $\restr{\rho}{\fund{\Sigma_2}}=\eta$. In particular, both
            $\restr{\rho}{\fund{C_2}}$ and $\restr{\rho}{\fund{Y_1}}$ are abelian. Lemma \ref{lemma: extension abelian C2} implies that
            \[
            \restr{\rho}{\fund{C_2}}(h) = \restr{\rho}{\fund{\Sigma_1'}}(h_1)= \pm 1.
            \]
            Thus, since the representation $\restr{\rho}{\fund{\Sigma_1'}}$ extends to an abelian representation of $\fund{Y_1}$, the set $A(Y_1)$ intersect either the line $\{\psi_1 =0\} \subset\R{\Sigma_1'}$ or the line $\{\psi_1 =\pi\} \subset\R{\Sigma_1'}$. Therefore, $A(Y) \subseteq S$. We conclude that $A(Y)=S$.

            Let us focus on $H(Y)$.
            Let  $\eta =(-\theta_1+\pi,\psi_1)$ be a point of $H_{Y,t} \subset \R{\Sigma_2'}$. In particular, by definition of $H_{Y,t}$, the point $(\theta_1,\psi_1) \in \R{\Sigma_1'}$ is in $T(Y_1,\partial Y_1)$.  According to Lemma \ref{lemma: irr extension C2, translazione} there exists an irreducible representation $\rho_c \colon \fund{C_2} \to SU(2)$ such $\restr{\rho_c}{\fund{\Sigma_2'}}=\eta$ and $\restr{\rho_c}{\fund{\Sigma_1'}}$ is the point with coordinates $(\theta_1,\psi_1) \in \R{\Sigma_1'}$. Since $(\theta_1,\psi_1) \in T(Y_1,\partial Y_1)$, then the restriction $\restr{\rho_c}{\fund{\Sigma_1'}}$ extends to $\fund{Y_1}$ by definition. This implies that $\eta$ extends to an irreducible representation $\fund{Y}\to SU(2)$. Therefore, $H_{Y,t} \subseteq H(Y)$.

            Similarly, if $\eta \in H_{Y,\varepsilon}$, then by Lemma \ref{lemma: irr extension C2, central h} there exists a representation $\rho \in R(Y)$ such that either $\restr{\rho}{\fund{C_2}}$ or $\restr{\rho}{\fund{Y_1}}$ is irreducible. Hence, $\rho$ is irreducible and $\restr{\rho}{\fund{\Sigma_2'}}\equiv \eta$. Therefore, $\eta \in H(Y)$ and $H_{Y,\varepsilon} \subseteq H(Y)$. We conclude that
            \[
                H_{Y,t} \cup H_{Y,0}\cup H_{Y,\pi} \subseteq H(Y) \subseteq \R{\Sigma_2'}.
            \]
        \end{proof}
\end{teo}

\section{Structural Behaviour}

In this section, we prove Theorem \ref{teo: admissible pieces have the property P}, which gives a further description of $T(Y,\partial Y)$ when $Y$ is a graph manifold rational homology $3$-sphere with torus boundary.
Moreover, we prove Proposition \ref{prop: sigma =1 ish} and Corollary \ref{cor: o1=1 and o2=3 and delta=1} as a consequence of Theorem \ref{teo: admissible pieces have the property P}.

\begin{defn}\label{defn: the line}
Let $Y$ be a $3$-manifold with torus boundary.
Let $\R{\partial Y}$ come with coordinates $(\theta,\psi)$.
Let $L \subset \R{\partial Y}$ be defined as
\[
    L_{\nicefrac{r}{s}}+c \coloneq \left\{(\theta,\psi)\in \R{\partial Y} \middle| r \theta + s \psi_1 \equiv_{2\pi} c \right\},
\]
where $\nicefrac{r}{s}\in \mathbb{Q} \cup \{\nicefrac{1}{0}\}$ and $c \in [0,2\pi]$.
We say that $L \subset \R{\partial Y}$ is a \emph{line} of the torus $\R{\partial Y}$, and that $L$ has slope $\nicefrac{r}{s}$.
\end{defn}

Let $\Sigma$ be a torus. We say that a path $\gamma \colon (-1,1)\to \R{\Sigma}$ is \emph{trivial} if it is homotopic relative to the set
\[
    \left\{\lim_{x \mapsto \pm 1} \gamma(x)\right\} \subset \R{\Sigma}
\]
to a constant path.
For instance, a path that starts from a point, travels for a finite amount of time, and returns to the starting point on its track is a trivial path.

\begin{defn}\label{defn: property P strong}
    Let $Y$ be a $3$-manifold with torus boundary. We say that $Y$ is \emph{LIPA} if there exists a nontrivial path $\gamma \colon (-1,1)\to H(Y)$ such that $\gamma \not \subset A(Y)$,
    \[
        \lim_{x \mapsto \pm 1} \gamma(x) \in A(Y)
    \]
    and there exists a line of $\R{\partial Y}$ containing $\gamma$. Such a path is called a \emph{LIPA-path} of $Y$.
\end{defn}

Roughly speaking, we say that a manifold is \emph{LIPA} if it admits a \textbf{L}inear \textbf{P}ath connecting two \textbf{A}belian points and consisting of \textbf{I}rreducible representations. The name \emph{LIPA} comes from the acronym of these four words.

Notice that if a manifold $Y$ with torus boundary is LIPA, then $H(Y) \subset \R{\partial Y}$ is nonempty and in particular $Y \neq S^1 \times \mathbb{D}^2$. In this section we are going to prove that if $Y$ is a graph manifold rational homology solid torus and $Y \neq S^1 \times \mathbb{D}^2$, then $Y$ is LIPA.

Let $n \ge 3$, it is easy to see that for every $1\le i \le n-1$ there exists a pair of vertical tori $\Sigma_1,\Sigma_2 \subset \mathbb{D}^2(\nicefrac{p_1}{q_1},\cdots, \nicefrac{p_n}{q_n})$ such that
\begin{equation}
\label{eq: una volta, dec con m fibers}
    \overline{\mathbb{D}^2(\nicefrac{p_1}{q_1},\cdots, \nicefrac{p_n}{q_n}) \setminus \left(\Sigma_1\cup \Sigma_2\right)} = \mathbb{D}^2(\nicefrac{p_1}{q_1},\cdots, \nicefrac{p_i}{q_i}) \cup C_3 \cup \mathbb{D}^2(\nicefrac{p_{i+1}}{q_{i+1}},\cdots \nicefrac{p_n}{q_n}).
\end{equation}
This implies that if $Y$ is a graph manifold rational homology solid torus, then we can choose a system of disjoint tori $\tau \subset Y$ as in Corollary \ref{cor: the useful decomposition} such that
\[
\overline{Y\setminus \tau} = \left(\bigcup_{k\in K} C_3\right) \cup \left(\bigcup_{j \in J} C_2\right) \cup \left(\bigcup_{i \in I}  D_i\right),
\]
where $D_i$ indicates a Seifert fibered space with disk base space and at most two singular fibers.

\begin{lemma}\label{lemma: the manifold D(p1,p2) has the property P}
    The Seifert fibered $Y = \mathbb{D}^2(\nicefrac{p_1}{q_1},\nicefrac{p_2}{q_2})$ with $2 \le p_1 \le p_2$ is LIPA.
    \begin{proof}
        Without loss of generality, assume that $q_1$ and $q_2$ are odd. We consider $\fund{Y}$ presented as follows:
        \[
            \fund{Y} = \langle x_1,x_2,h \,|\, x_1^{p_1}h^{q_1},x_2^{p_2}h^{q_2}, \text{$h$ central} \rangle.
        \]
        Let $\mu$ be the slope $\partial Y$ such that $[\mu]=x_1x_2 \in \fund{Y}$. For details see \cite[Definition 2.1.1]{Mino} or \cite[Section 2.2]{LspaceIntervalGraphManifold}.
        
        By Fact \ref{Fact: centralizers}, if $\rho \in R(Y)$ is irreducible, then $\rho(h) \in \ZSU$. We define $H(Y)_\pi \subset H(Y)$ as
        \[
            H(Y)_\pi = \left\{\eta \in H(Y) \middle| \eta(h)=e^{i\pi}=-1\right\}.
        \]
        According to \cite[Lemma 7.8]{Mino}, the following map is surjective and two-to-one:
        \[
    f \colon H(Y)_\pi \to  J_\pi(p_1,p_2) \quad \text{with} \quad
       \eta \mapsto \Tr\eta(\mu),
   \]
   where
   \[
    J_\pi(p_1,p_2) = \bigcup_{\substack{k_1,k_2 \in \mathbb{Z} \\ k_i \text{ odd} }} I \left(2 \cos\left( \frac{\pi k_1}{p_1}\right), 2 \cos\left( \frac{\pi k_2}{p_2}\right) \right),
   \]
   and the interval $I \subseteq (-2,2)$ is as defined in \eqref{eq: interval I}.
   According to \cite[Lemma 7.6]{Mino}, the interval $J_\pi(p_1,p_2)$ is nonempty.
   For $k_1,k_2$ odd, we write $I_{k_1,k_2}$ for the interval $I\left(2\cos (\nicefrac{\pi k_1}{p_1}),2\cos (\nicefrac{\pi k_2}{p_2})\right)$. Since $J_\pi(p_1,p_2)$ is nonempty, at least one interval $I_{k_1,k_2}$ is nonempty. Let us assume that the interval $I_{k_1,k_2}$ is nonempty. The latter is connected and the preimage $f^{-1}I_{k_1,k_2}$ has two connected components in $H(Y)_{\pi}$. Thus, $f^{-1}I_{k_1,k_2}$ defines two paths
   \[
   \gamma_{k_1,k_2}^{1},\gamma_{k_1,k_2}^{2} \colon (-1,1) \to H(Y)_{\pi}.
   \]
   We are going to prove that these two paths are LIPA-paths for $Y$ by verifying the conditions of Definition \ref{defn: property P strong}.
    
    Let $\alpha\in \{1,2\}$, the path $\gamma_{k_1,k_2}^{\alpha}$ is nontrivial and contained in the line
    \[
        L_{\pi} \coloneq \left\{ \eta \in \R{\partial Y} \middle| \eta(h)=e^{i\pi}\right\}.
    \]
    Let $\lambda \subset \partial Y$ be the rational longitude of $Y$.
    According to \cite[Lemma 5.1]{Mino}, the distance $\geo{\lambda}{h} \neq0$ and therefore $\lambda_Y \neq h$.
    As the rational longitude of $Y$ is not the regular fiber $h$, the path $\gamma_{k_1,k_2}^{\alpha}$ is not contained in $A(Y)$ by Corollary \ref{cor: A(Y) with rational longitude}.

   The function $f$ above shows that
   \[
    \lim_{x \mapsto \pm 1} \Tr\left(\left.\gamma_{k_1,k_2}^{\alpha}\middle(x\right)(\mu)\right) \in \partial I_{k_1,k_2} \quad \text{with} \quad \alpha \in \{1,2\}.  
    \]
    We recall that
    the endpoints of the interval $I_{k_1,k_2}$ are
    \begin{equation*}
    \left\{ 2 \cos \left( \pi \left( \frac{k_1}{p_1} + \frac{k_2}{p_2}\right) \right), 2 \cos \left( \pi \left( \frac{k_1}{p_1} - \frac{k_2}{p_2}\right)\right) \right\}.
\end{equation*}
This implies that if $\eta \in \R{\partial Y}$ is a representation such that
\[
    \eta \in \left\{ \lim_{x \mapsto \pm 1} \gamma_{k_1,k_2}^{\alpha}(x) \right\}_{\alpha=1,2},
\]
then $\eta(h)=-1$ and $\eta(\mu)=e^{i\theta}$, where
\[
    \theta \in \left\{\pi \left( \frac{k_1}{p_1} + \frac{k_2}{p_2}\right),\pi \left( \frac{k_1}{p_1} - \frac{k_2}{p_2}\right),2\pi -\pi \left( \frac{k_1}{p_1} + \frac{k_2}{p_2}\right),2\pi-\pi \left( \frac{k_1}{p_1} - \frac{k_2}{p_2}\right)\right\}.
\]
Let $g \coloneq \gcd(p_1,p_2)$. As a consequence of \cite[Lemma 5.1]{Mino}, we obtain that
\[
    \lambda_Y= \mu \frac{p_1p_2}{o_Yg} + h \frac{q_1p_2 +q_2p_1}{o_Yg} \subset \partial Y.
\]
It is straightforward to prove that if $\eta$ is as above, then $\eta(\lambda_Y)^{o_Y}=1$. This implies that $\eta \in A(Y)$ by Corollary \ref{cor: A(Y) with rational longitude}. Hence, for $\alpha \in \{1,2\}$
\[
    \lim_{x \mapsto \pm 1} \gamma_{k_1,k_2}^{\alpha} \in A(Y).  
\]
    \end{proof}
\end{lemma}

\begin{lemma}\label{lemma: the sum is empty then}
    Let $Y=Y_1 \cup_{\Sigma_1}C_3 \cup_{\Sigma_2} Y_2$. If $H(Y_2) \neq \emptyset$ and the sum $A(Y_1)+H(Y_2) \subset \R{\Sigma_3}$ is empty, then $\Delta(h_1,\lambda_{1})=0$. 
    \begin{proof}
        Let $n,m \in \mathbb{Z}$ be integers such that $\lambda_{1}=nx_1+mh_1$ in $\fund{\Sigma_1}$. Thus, $|n| = \Delta(\lambda_{1},h_1)$. By Corollary \ref{cor: A(Y) with rational longitude},
        \[
            A(Y_1)= \left\{ (\theta_1,\psi_1) \in \R{\Sigma_1} \middle| o_1 \left( n\theta_1+m\psi_1\right)=0\right\}.
        \]
        Let $(\theta_0,\psi_0) \in H(Y_2)$, such a point exists since $H(Y_2) \neq \emptyset$ by hypothesis. If $A(Y_1) \cap \{\psi_1= \psi_0\} \neq \emptyset$, then $A(Y_1) + H(Y_2) \neq \emptyset$ by Definition \ref{defn: A+B}. Therefore, since $A(Y_1) + H(Y_2)$ is empty by hypothesis, then
        $
            A(Y_1) \cap \{\psi_1= \psi_0\} = \emptyset.
        $
        This implies that $|n|=\Delta(\lambda_{1},h_1)=0$.
    \end{proof}
\end{lemma}

\begin{lemma}\label{C3 has the property P}
    Let $Y=Y_1 \cup_{\Sigma_1} C_3 \cup_{\Sigma_2} Y_2$ be a graph manifold rational homology solid torus. Let us assume that $Y$ is not a solid torus. For $j \in \{1,2\}$, if $Y_j$ is either LIPA or a solid torus, then the manifold $Y$ is LIPA.
    \begin{proof}
        If $Y_j$ is a solid torus, then, since $Y$ is a graph manifold, $\Delta(\lambda_{j},h_j) \neq 0$. Hence, if $Y_1$ and $Y_2$ are two solid tori, the manifold $Y$ is a Seifert fibered space with disk base space. Since $Y$ is not a solid torus by assumption, the conclusion holds by Lemma \ref{lemma: the manifold D(p1,p2) has the property P}. Therefore, we assume that either $Y_1$ or $Y_2$ is not a solid torus. Without loss of generality, we assume that $Y_2$ is not a solid torus and, by hypothesis, this implies that $Y_2$ is LIPA.

        Since $Y_2$ is LIPA, the set $H(Y_2)$ is nonempty. If $\geo{\lambda_1}{h_1} =0$, then, since $Y$ is a graph manifold, $Y_1$ is not a solid torus and, by hypothesis, it is LIPA. Thus, $H(Y_1) \neq \emptyset$ and, according to Corollary \ref{cor: C3 obstruction longitudes}, we obtain that $\geo{\lambda_2}{h_2}\neq 0$. Therefore, up to swapping $Y_1$ and $Y_2$, we can assume that $\geo{\lambda_1}{h_1}\neq 0$. According to Lemma \ref{lemma: the sum is empty then} the sum $A(Y_1)+H(Y_2) \subset\R{\Sigma_3}$ is not empty.
        Let $\gamma : (-1,1)\to H(Y_2)$ be an LIPA-path of $Y_2$ in $\R{\Sigma_2}$. Thus,
        \[
            \lim_{x \to \pm 1} \gamma(x) \in A(Y_2).
        \]
        As shown in the proof of Lemma \ref{lemma: the sum is empty then}, since $\geo{\lambda_1}{h_1}\neq 0$, for every $x \in (-1,1)$
        \[
            A(Y_1) + \gamma(x) \neq \emptyset.
        \]
        Therefore, according to Theorem \ref{teo: alg for C3},
        \[
        A(Y_1) + \gamma(x) \in A(Y_1)+H(Y_2) \subseteq H(Y) \subseteq \R{\Sigma_3}.
        \]
        Thus, the sum $A(Y_1) + \ima{\gamma}$ defines a path $\delta \colon (-1,1)\to H(Y)$. 
        Furthermore, we obtain that
        \[
            \left\{\lim_{x \to \pm 1} \delta(x)\right\}\in \lim_{x \to \pm 1} \left(A(Y_1) + \gamma(x)\right) \subseteq A(Y_1)+A(Y_2) = A(Y).
        \]
        Since $\gamma \not \subset A(Y_1) $, we obtain that $\delta \not \subset A(Y)$. This implies that $Y$ is LIPA.
        \end{proof}
\end{lemma}

\begin{lemma}\label{C2 has the property P}
    Let $Y=Y_1 \cup_{\Sigma_1'} C_2$ be a graph rational homology solid torus. If $Y_1$ either is LIPA or it is a solid torus, then the manifold $Y$ is  LIPA.
    \begin{proof}
        We assume $\R{\Sigma_1'}$ and $\R{\Sigma_2'}$ to be parameterized as in \eqref{eq: coordinates for C2}.
        Let $n,m \in \mathbb{Z}$ be integers such that
        \[
            \lambda_1 = n x_1 + m h_1 \in \fund{\Sigma_1}.
        \]
        Thus, $|n|=\Delta(\lambda_1,h_1)$. Corollary \ref{cor: A(Y) with rational longitude} implies that
        \[
            A(Y_1)= \left\{(\theta_1,\psi_1)\in \R{\Sigma_1'} \middle| o_1 (n \theta_1+m\psi_1)=0\right\}
        \]
        According to Corollary \ref{cor: C2 obstruction longitudes}, $|n|=\Delta(\lambda_1,h_1)\neq 0$. By Theorem \ref{teo: alg for C2}, the set $A(Y)\subset \R{\Sigma_2}$ is equal to $\{\psi_2=0\}\cup\{\psi_2=\pi\}$. We define $H'$ as
        \begin{align*}
            H'&= \left\{(\pi -\theta_1,\psi_1)\in \R{\Sigma_2'} \middle| (\theta_1,\psi_1) \in A(Y_1)\right\} \\
            &= \left\{(\pi -\theta_1,\psi_1)\in \R{\Sigma_2'} \middle| o_1 (n\theta_1+m\psi_1)=0\right\} \\
            &=\left\{(\theta_2,\psi_2) \in \R{\Sigma_2'} \middle| o_1\left( n\theta_2 -m\psi_2 - \pi n \right)=0\right\} .
         \end{align*}
        According to Theorem \ref{teo: alg for C2}, $H' \subseteq H(Y)$.
        In fact $H'$ is the subset of $ H_{Y,t}$ of Theorem \ref{teo: alg for C2}. The set $H'$ is made of $o_1$-many parallel lines each of slope $\nicefrac{-n}{m}$.
        Since $n \neq 0$, Theorem \ref{teo: alg for C2} implies that
        \[
            H' \not\subset A(Y)= \left\{ (\theta_2,\psi_2) \in \R{\Sigma_2'} \middle| \psi_2=0,\pi \right\}
        \]
        Moreover, each line in $H'$ intersects both lines of $A(Y)$. Hence, $H'\subset H_{Y,t}$ contains a LIPA path. Therefore, $Y$ is LIPA.
    \end{proof}
\end{lemma}

\begin{repteo}{teo: admissible pieces have the property P}
    Let $Y$ be a graph manifold rational homology solid torus. If $Y$ is not a solid torus, then there exists a path $\gamma \colon (-1,1) \to \R{\partial Y}$ of representations such that:
    \begin{enumerate}
        \item the path starts and ends at representations $\fund{\partial Y} \to \Ui$ that extend to an abelian representation $\fund{Y} \to SU(2)$;
        \item for every $x \in (-1,1)$ the representation $\gamma(x)$ extends to an irreducible representation $\fund{Y} \to SU(2)$;
        \item the path $\gamma$ is non-trivial in $\R{\partial Y}$ and consists of a line segment.
    \end{enumerate}
    \begin{proof}
    We can prove the conclusion by showing that $Y$ is LIPA.
    As we explained before, Corollary \ref{cor: the useful decomposition} implies the existence of a system of disjoint tori $\tau$ such that $\overline{Y \setminus \tau}$ is a union of $C_3$'s, $C_2$'s, solid tori, and Seifert fibered spaces with two singular fibers and disk base space. Lemma \ref{lemma: the manifold D(p1,p2) has the property P}, Lemma \ref{C3 has the property P}, and Lemma \ref{C2 has the property P} imply the conclusion. 
    \end{proof}
\end{repteo}

As a consequence of Theorem \ref{teo: admissible pieces have the property P}, we show Proposition \ref{prop: sigma =1 ish} and Corollary \ref{cor: o1=1 and o2=3 and delta=1}.

\begin{figure}[t]
    \centering
    \begin{subfigure}[b]{0.4\textwidth}
    \centering
        \begin{tikzpicture}[>=latex,scale=1.45]
            \draw[->] (-.5,0) -- (2.5,0) node[right] {$\theta$};
            \foreach \x /\n in {1/$\pi$,2/$2\pi$} \draw[shift={(\x,0)}] (0pt,2pt) -- (0pt,-2pt) node[below] {\tiny \n};
            \draw[->] (0,-.5) -- (0,2.5) node[below right] {$\psi$};
            \foreach \y /\n in {1/$\pi$,2/$2\pi$}
            \draw[shift={(0,\y)}] (2pt,0pt) -- (-2pt,0pt) node[left] {\tiny \n};
            \node[below left] at (0,0) {\tiny $0$};
             \draw (0,0) rectangle (2,2);
             \draw (1,0)--(1,2);
             \draw (0,1)--(2,1);

             \draw [blue, thick] (0,0)--(0,2);
             \draw [blue, thick] (2,0)--(2,2);
             \draw [cyan, thick] (0,0)--(2,0);
              \draw [cyan, thick] (0,2)--(2,2);
             \draw [orange, thick] (0,0)--(1,2);
                \draw[orange, dashed, thick] (1,0)--(2,2);
             \draw[red, thick] (0,0)--(2,1);
             \draw[red, dashed, thick] (0,1)--(2,2);
        \end{tikzpicture}
        \caption{Case $\Delta(\lambda_1,\lambda_2)o_1o_2 = 1$.}
        \label{Figure: sigma small. s=1}
    \end{subfigure}
    \quad
    \begin{subfigure}[b]{0.4\textwidth}
    \centering
        \begin{tikzpicture}[>=latex,scale=1.45]
            \draw[->] (-.5,0) -- (2.5,0) node[right] {$\theta$};
            \foreach \x /\n in {1/$\pi$,2/$2\pi$} \draw[shift={(\x,0)}] (0pt,2pt) -- (0pt,-2pt) node[below] {\tiny \n};
            \draw[->] (0,-.5) -- (0,2.5) node[below right] {$\psi$};
            \foreach \y /\n in {1/$\pi$,2/$2\pi$}
            \draw[shift={(0,\y)}] (2pt,0pt) -- (-2pt,0pt) node[left] {\tiny \n};
            \node[below left] at (0,0) {\tiny $0$};
             \draw (0,0) rectangle (2,2);
             \draw (1,0)--(1,2);
             \draw (0,1)--(2,1);

             \draw [blue, thick] (0,0)--(0,2);
             \draw [blue, thick] (2,0)--(2,2);
             
             \draw [cyan, thick] (0,0)--(2,0);
             \draw [cyan, thick] (0,2/3)--(2,2/3);
             \draw [cyan, thick] (0,4/3)--(2,4/3);
              \draw [cyan, thick] (0,2)--(2,2);
              
             \draw [orange, thick] (0,0)--(2,1/3);
             \draw [orange, dashed, thick] (0,1)--(2,1+1/3);
             \draw [red, thick] (0,4/3)--(1,2);
             \draw [red, dashed, thick] (1,4/3)--(2,2);
             \draw [red, dotted ,thick] (1,2/3)--(2,0);
        \end{tikzpicture}
        \caption{Case $\Delta(\lambda_1,\lambda_2)=o_2=1$ and $o_1=3$.}
        \label{Figure o1=1 e o2=3}
    \end{subfigure}
\\
    \begin{subfigure}[b]{0.4\textwidth}
    \centering
    \begin{tikzpicture}[>=latex,scale=1.45]
            \draw[->] (-.5,0) -- (2.5,0) node[right] {$\theta$};
            \foreach \x /\n in {1/$\pi$,2/$2\pi$} \draw[shift={(\x,0)}] (0pt,2pt) -- (0pt,-2pt) node[below] {\tiny \n};
            \draw[->] (0,-.5) -- (0,2.5) node[below right] {$\psi$};
            \foreach \y /\n in {1/$\pi$,2/$2\pi$}
            \draw[shift={(0,\y)}] (2pt,0pt) -- (-2pt,0pt) node[left] {\tiny \n};
            \node[below left] at (0,0) {\tiny $0$};
             \draw (0,0) rectangle (2,2);
             \draw (1,0)--(1,2);
             \draw (0,1)--(2,1);

             \draw [blue, thick] (0,0)--(0,2);
              \draw [blue, thick] (1,0)--(1,2);
             \draw [blue, thick] (2,0)--(2,2);
             \draw [cyan, thick] (0,0)--(2,0);
              \draw [cyan, thick] (0,2)--(2,2);
             \draw [orange, thick] (0,0)--(1,1);
             \draw [red, thick] (0.2,0)--(0.8,2);
             \draw [red, thick] (1.2,0)--(1.8,2);
        \end{tikzpicture}
        \caption{Case $\Delta(\lambda_1,\lambda_2)=o_2=1$ and $o_1=2$.}
        \label{Figure o1=1 e o2=2}
    \end{subfigure}
    \caption{The sets $A(Y_1)$ and $A(Y_2)$ are in blue and cyan respectively. The paths of irredubile representations for $Y_1$ and $Y_2$ are in orange and red. The corresponding translations are dashed.}
    \end{figure}
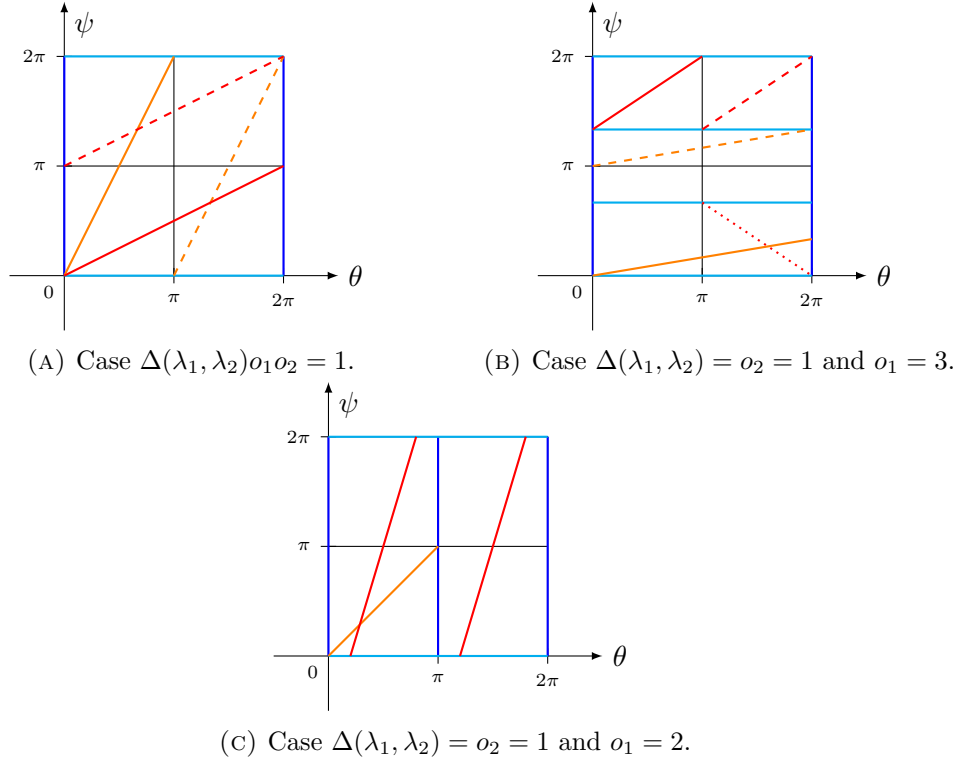

\begin{repprop}{prop: sigma =1 ish}
     Let $Y_1$ and $Y_2$ be two graph manifold with torus boundary and null-homologous longitudes. If neither $Y_1$ nor $Y_2$ is a solid torus and $Y=Y_1 \cup_{\Sigma} Y_2$ is a rational homology $3$-sphere such that $\Delta(\lambda_1,\lambda_2)= 1$, then $Y$ is not $SU(2)$-abelian.
     \begin{proof}
    According to Theorem \ref{teo: admissible pieces have the property P}, the manifolds $Y_1$ and $Y_2$ are both LIPA.
    Since $\geo{\lambda_1}{\lambda_2}=1$, the ordered set $\{\lambda_1,\lambda_2\}$ is a basis for $\fund{\Sigma}$. We give the space $\R{\Sigma}$ coordinates $(\theta,\psi)$ according to this ordered basis as in \eqref{eq: general coordinates}. According to Corollary \ref{cor: A(Y) with rational longitude}, the set $\{\theta=0\} \subset\R{\Sigma}$ (resp. $\{\psi=0\}\subset\R{\Sigma}$) is equal to $A(Y_1)$ (resp. $A(Y_2)$).
         For $i \in \{1,2\}$, let $\gamma_i \colon (-1,1)\to \R{\Sigma}$ be a LIPA-path of $Y_i$. Let $\gamma_1^t$ (resp. $\gamma_2^t$) the translation of $\gamma_1$ (resp. $\gamma_2$) by $(0,\pi)$ (resp. $(\pi,0)$). Lemma \ref{lemma: translation} implies that $\gamma_i^t$ is a LIPA path for $Y_i$.

         We call $\Gamma_1$ and $\Gamma_1^t$ the images of $\gamma_1$ and $\gamma_1^t$. We define $\Gamma_2$ and $\Gamma_2^t$ similarly. 
         If either $\Gamma_1 \cup \Gamma_1^t$ intersects $A(Y_2)$ or $\Gamma_2 \cup \Gamma_2^t$ intersects $A(Y_1)$, then either
         $H(Y_1)\cap A(Y_2) \neq \emptyset$ or $A(Y_1)\cap H(Y_2) \neq \emptyset$ respectively. In both cases, the conclusion holds by Theorem \ref{teo: M SU(2)-abeliano se e solo se i pezzi sono empty}.

         Let us assume that $\Gamma_1 \cup \Gamma_1^t$ and $\Gamma_2 \cup \Gamma_2^t$ do not intersect $A(Y_2)=\{\theta=0\}$ and $A(Y_1)=\{\psi=0\}$ respectively.
         Since $\Gamma_1$ is the image of a LIPA-path of $Y_1$, it connects the points
         \[
            (0,\psi_1) \in A(Y_1) \quad \text{and} \quad (0,\psi_2) \in A(Y_1).
         \]
         Therefore $\Gamma_1^t$ connects the points
         \[
            (0,\psi_1+\pi) \in A(Y_1) \quad \text{and} \quad (0,\psi_2+\pi) \in A(Y_1).
         \]
         Since we assumed that $\Gamma_1 \cup \Gamma_1^t$ is disjoint from $A(Y_2)=\{\psi=0\}$, we can assume that either $\psi_1,\psi_2 \in [0,\pi]$ or $\psi_1,\psi_2\in [\pi,2\pi]$. Similarly, $\Gamma_2$ is the image of a LIPA-path of $Y_2$, then it connects the points
         \[
            (\theta_1,0) \in A(Y_2) \quad (\theta_2,0) \in A(Y_2),
         \]
         and we can assume that either $\theta_1,\theta_2 \in [0,\pi]$ or $\theta_1,\theta_2 \in [\pi,2\pi]$.
         Figure \ref{Figure: sigma small. s=1} shows the union
         \[
            \left(A(Y_1) \cup \Gamma_1\cup \Gamma_1^t\right) \cup \left(A(Y_2) \cup \Gamma_2\cup \Gamma_2^t\right) \subset \R{\Sigma},
         \]
         in $(\theta,\psi)$-coordinates.
         In particular, Figure \ref{Figure: sigma small. s=1} shows that $H(Y_1) \cap H(Y_2) \neq \emptyset$. We conclude that $Y=Y_1 \cup_{\Sigma} Y_2$ is not $SU(2)$-abelian, by Theorem \ref{teo: M SU(2)-abeliano se e solo se i pezzi sono empty}.
     \end{proof}
\end{repprop}

Let $Y$ be a $3$-manifold as in Proposition \ref{prop: sigma =1 ish}, if $Y_1$ and $Y_2$ are \emph{integer} homology solid tori, then $Y$ is an integer homology $3$-sphere. Proposition \ref{prop: sigma =1 ish} implies that $Y$ is not $SU(2)$-abelian. This is also proven in \cite{ToroidalZS3HaveIrreducibleRepres} in which the authors proved that every toroidal integer homology $3$-sphere is not $SU(2)$-abelian. In particular, they showed that if a rational homology solid torus $Y$ whose longitude $\lambda_Y$ is integrally null-homologous admits a similar LIPA-path whenever $Y(\lambda_Y)$ has non-trivial instanton Floer homology (see \cite[Theorem 3.5]{ToroidalZS3HaveIrreducibleRepres}).

\begin{repcor}{cor: o1=1 and o2=3 and delta=1}
    Let $Y_1$ and $Y_2$ be two graph manifolds. If neither $Y_1$ nor $Y_2$ is a solid torus, $o_1=1$, $o_2=3$, and $Y=Y_1 \cup_{\Sigma} Y_2$ is a rational homology $3$-sphere such that $\Delta(\lambda_1,\lambda_2)= 1$, then $Y$ is not $SU(2)$-abelian.
    \begin{proof}
        The proof is similar to Proposition \ref{prop: sigma =1 ish}.
        Since $\geo{\lambda_1}{\lambda_2}=1$, the set $\{\lambda_1,\lambda_2\}$ is an ordered basis for $\fund{\Sigma}$. We give the space $\R{\Sigma}$ coordinates $(\theta,\psi)$ according to this ordered basis as in \eqref{eq: general coordinates}.
        According to Corollary \ref{cor: A(Y) with rational longitude} implies that $A(Y_1)=\{\psi=0\} \subset \R{\Sigma}$ and
        \[
            A(Y_2)=\bigcup_{k =0,1,2}\left\{\psi= \frac{2\pi k}{3}\right\}  \subset \R{\Sigma}.
        \]
        For $i \in \{1,2\}$ let $\gamma_i,\gamma_i^t,\Gamma_i,\Gamma_i^t$ be as in the proof of Proposition \ref{prop: sigma =1 ish}. If either $\Gamma_1 \cup \Gamma_1^t$ intersects $A(Y_2)$ or $\Gamma_2 \cup \Gamma_2^t$ intersects $A(Y_1)$, then either
         $H(Y_1)\cap A(Y_2) \neq \emptyset$ or $A(Y_1)\cap H(Y_2) \neq \emptyset$ respectively. In both cases, the conclusion holds by Theorem \ref{teo: M SU(2)-abeliano se e solo se i pezzi sono empty}.

         Let us assume that $\Gamma_1 \cup \Gamma_1^t$ and $\Gamma_2 \cup \Gamma_2^t$ do not intersect $A(Y_2)$ and $A(Y_1)$ respectively.
         Since $\Gamma_2$ is the image of a LIPA-path of $Y_2$, it connects the points
         \[
            \left(\theta_1,\frac{ 2 \pi k}{3}\right) \in A(Y_2) \quad \left(\theta_2,\frac{2 \pi k}{3}\right) \in A(Y_2),
         \]
         with $k \in \{0,1,2\}$.
         Since we assumed that $\Gamma_2 \cup \Gamma_2^t$ is disjoint from $A(Y_1)=\{\theta=0\}$, we can assume that either $\theta_1,\theta_2 \in [0,\pi]$ or $\theta_1,\theta_2 \in [\pi,2\pi]$.
         Similarly, $\Gamma_1$ is the image of a LIPA-path of $Y_2$, then is connects the points
         \[
            (0,\psi_1) \in A(Y_1) \quad \text{and} \quad (0,\psi_2) \in A(Y_1).
         \]
         Since we assumed that $\Gamma_1$ is disjoint from $A(Y_2)$, then there exists an $k' \in \{0,1,2\}$ such that
         \[
            \psi_1,\psi_2 \in \left[\frac{2\pi k'}{3}, \frac{2 \pi (k'+1)}{3}\right].
         \]
         Similarly, $\Gamma_1^t=\Gamma_1 +(0,\pi)$ connects the points
         \[
            (0,\psi_1+\pi) \in A(Y_1) \quad \text{and} \quad (0,\psi_2+\pi) \in A(Y_1).
         \]
         Therefore, since we assume that $\Gamma_1^t$ is disconnect from $A(Y_2)$, we have that there exists an $k'' \in \{0,1,2\}$ such that
         \[
            \psi_1+\pi,\psi_2 + \pi \in \left[\frac{2\pi k''}{3}, \frac{2 \pi (k''+1)}{3}\right].
         \]
         This implies that there exists a $k ''' \in \{0,\cdots,5\}$ such that
         \[
         \psi_1,\psi_2 \in \left[ \frac{\pi k'''}{3}, \frac{\pi (k'''+1)}{3} \right].
         \]
         If $\Gamma_1\cup \Gamma_1^t$ interests $\Gamma_2\cup \Gamma_2^t$, then $Y=Y_1 \cup_{\Sigma} Y_2$ is not $SU(2)$-abelian by Theorem \ref{teo: M SU(2)-abeliano se e solo se i pezzi sono empty}. Let us assume that this intersection is empty. Without loss of generality, we can assume that $\Gamma_1$ connects the points
         \[
            (0,\psi_1) \in A(Y_1) \quad \text{and} \quad (0,\psi_2) \in A(Y_1), \quad \text{with} \quad \psi_1,\psi_2 \in \left[0, \frac{\pi}{3}\right]
         \]
        Since $(\Gamma_1\cup \Gamma_1^t) \cap(\Gamma_2\cup \Gamma_2^t)$ is empty by assumption, we conclude that $\Gamma_2$ connects the points
        \[
            \left(\theta_1,\frac{4 \pi}{3}\right) \in A(Y_2) \quad \text{and} \quad \left(\theta_2,2\pi\right) \in A(Y_1), \quad \text{with} \quad \theta_1,\theta_2 \in \left[0, \pi \right].
         \]
         Figure \ref{Figure o1=1 e o2=3} shows the union
         \[
            \left(A(Y_1) \cup \Gamma_1\cup \Gamma_1^t\right) \cup \left(A(Y_2) \cup \Gamma_2\cup \Gamma_2^t\right) \subset \R{\Sigma},
         \]
         in $(\theta,\psi)$-coordinates. We define
         \[
            \Gamma_2^\dagger = \left\{ (2\pi -\theta, 2\pi - \psi) \middle| (\theta,\psi) \in \Gamma_2 \right\}.
         \]
         Corollary \ref{cor: jewelled repr} implies that $\Gamma_2^\dagger \subset H(Y_2)$, this is the red dotted line drawn in Figure \ref{Figure o1=1 e o2=3}. In particular, $\Gamma_2^\dagger$ connects the points
         \[
            \left(2\pi -\theta_1,\frac{2 \pi}{3}\right) \in A(Y_2) \quad \text{and} \quad \left(2 \pi-\theta_2,0\right) \in A(Y_1), \quad \text{with} \quad \theta_1,\theta_2 \in \left[0, \pi \right].
         \]
         This implies, as is shown in Figure \ref{Figure o1=1 e o2=3}, that the set $\Gamma_2^\dagger$ intersects $\Gamma_1\cup \Gamma_1^t$. This implies that $H(Y_1)\cap H(Y_2) \neq \emptyset$. The conclusion is given by Theorem \ref{teo: M SU(2)-abeliano se e solo se i pezzi sono empty}.
    \end{proof}
\end{repcor}

Let $Y=Y_1 \cup_{\Sigma} Y_2$ be a toroidal manifold either as in Proposition \ref{prop: sigma =1 ish} or as in Corollary \ref{cor: o1=1 and o2=3 and delta=1}.
We denote by $\mathcal{L}(Y_1)$ and $\mathcal{L}(Y_2)$ the \emph{L-space intervals} of $Y_1$ and $Y_2$ respectively, details can be found in \cite{LspaceIntervalGraphManifold}.
As is well-documented, the rational longitude of a $3$-manifold with torus boundary is not an L-space slope. In particular $\lambda_2 \notin \mathcal{L}(Y_2)$. Moreover,
according to \cite[Theorem 1.3]{boyer2023slopedetectiontoroidal3manifolds}, the slope $\lambda_2 \subset \Sigma=\partial Y_1$ is NLS-detected for $Y_1$. This, in particular, means that $\lambda_2 \notin \mathcal{L}^\circ(Y_1)$. Therefore,
\[
    \lambda_2 \notin \mathcal{L}^\circ(Y_1) \cup \mathcal{L}^\circ(Y_2) \subset \{\text{slopes of } \Sigma\}.
\]
By \cite[Theorem 1.1]{boyer2023slopedetectiontoroidal3manifolds}, the manifold $Y=Y_1 \cup_{\Sigma} Y_2$ is not an L-space.
This implies that Proposition \ref{prop: sigma =1 ish} and Corollary \ref{cor: o1=1 and o2=3 and delta=1} are further evidences of Conjecture \ref{conj: SU(2)-ab implies HF L-space}.

\section{Auxiliary results}
\label{sec. utensilis}
In this section we introduce the tools needed for the proof of Theorem \ref{teo: sigma small then not SU(2)-abelian}. 
Specifically, we define a weighted $\mu$ for graph manifolds with torus boundary and show that, when the first homology group has small order, the $C_2$ pieces in the decomposition of Corollary \ref{cor: the useful decomposition} can be excluded.

Let $Y$ be a graph manifold rational homology solid torus. We recall that $\R{\partial Y}$ is a torus and that a rational line of $\R{\partial Y}$ is homeomorphic to a circle $S^1$. We know that $T(Y,\partial Y) = A(Y) \cup H(Y)$.
A line $L \subset \R{\partial Y}$ is \emph{rational} if $L$ is not dense. It it straightforward to see that if $\xi \subset \partial Y$ is a slope and $\alpha \in \ZSU= \{\pm 1\} \subset SU(2)$, then
\[
    \left\{ \eta \in \R{\partial Y} \middle| \eta(\xi) = \alpha\right\}
\]
is a rational line of $\R{\partial Y}$.

Let $\{x_1,x_2\}$ be a basis for $\fund{\partial Y}$, we parameterize $\R{\partial Y}$ with coordinates according to this basis as in \eqref{eq: general coordinates}.
We consider the torus $\R{\partial Y}$ as the product of two unitary circles according to these coordinates:
\[
    \R{\partial Y} = \{e^{i \theta}\} \times \{e^{i \psi}\} = S^1 \times S^1.
\]
We equip the torus $\R{\partial Y}$ with the natural Lebesgue measure $m$ such that $m(\R{\partial Y}) = (2\pi)^2$.
We recall that, for $\nicefrac{p}{q} \in \mathbb{Q} \cup \{\nicefrac{1}{0}\}$ and $c \in [0,2\pi]$, the line
\[
    L=L_{\nicefrac{p}{q}}+c= \left\{ (\theta,\psi) \in \R{\partial Y} \middle| p \theta + q\psi = c \right\}
\]
has slope $\nicefrac{p}{q}$.
Therefore, if the line $L \subset \R{\partial Y}$ has slope $\nicefrac{p}{q} \in \mathbb{Q} \cup \{\nicefrac{1}{0}\}$, then $L$ has length, with respect to the measure $m$, equal to $2\pi\sqrt{p^2 + q^2}$.
We equip the line $L$ with the Lebesgue measure $h$
induced by the measure $m$ on $\R{\partial Y}$. We recall that
$h$ is such that $h(L)=2\pi \sqrt{p^2 + q^2}$.

\begin{defn}\label{degn weight mu}
Let $Y$ be a graph manifold rational homology solid torus with the property that there exists be a finite set of rational lines $\{L_i\}_{i=1}^n$ of $\R{\partial Y}$ whose union contains $H(Y)$.
We define a weight $\mu$ on the $L_i$ as follows:
\[
    \mu(L_i)= \begin{dcases}
        \frac{h(H(Y) \cap L_i)}{h(L_i)} & \text{if $L_i$ contains a LIPA-path of $Y$} \\
        0 & \text{otherwise}.
    \end{dcases}
\]
Moreover, we define
\[
    \mu(Y) \coloneq \max_{i=1 , \cdots, n} \left\{ \mu (L_i)\right\}.  
\]
\end{defn}

Let $Y$ be a Seifert fibered space with disk base space.
As a consequence of Fact \ref{Fact: centralizers}, the set $H(Y)$ is contained in
\[
    \left\{ \eta \in \R{\partial Y} \middle | \eta(h) \in \ZSU\right\},
\]
which is the union of two rational lines of $\R{\partial Y}$. Here $h \subset \partial Y$ is a regular fiber of $Y$.
In particular, the weight $\mu(Y)$ is well defined for Seifert fibered spaces with disk base space. 
It can be shown that if $Y$ is a graph manifold rational solid torus, then there exists a finite collection of rational lines in $\R{\partial Y}$ whose union contains $H(Y)$. Consequently, the weight $\mu(Y)$ is well defined also in this case.
However, in this work we will only consider the weight $\mu$ for Seifert fibered spaces with disk base space, which, as mentioned above, is already well defined; hence we shall not provide a proof of this fact.

Roughly speaking, the weight $\mu(L)$ measures the portion of the line $L$ that is covered by a LIPA-path and $\mu(Y)$ is the weight of the longest LIPA-path within $H(Y)$.
Clearly, if $Y=S^1 \times \mathbb{D}^2$, then $H(Y) = \emptyset$ and therefore $\mu(Y) = 0$. Conversely, if $\mu(Y) \neq 0$, then $H(Y) \neq \emptyset$ and therefore $Y \neq S^1 \times \mathbb{D}^2$. 

\begin{lemma}
    The weight $\mu(Y)$ does not depend on the choice of the basis $\{x_1,x_2\}$.
    \begin{proof}
        Let $L_0 \subset \R{\partial Y}$ be a line such that $\mu(L_0) \neq 0$. We shall prove that $\mu(L_0)$ does not depend on the basis $\{x_1,x_2\}$. Without loss of generality, we assume that $L_0$ is of slope $\nicefrac{0}{1}$. This means that there exists a $c \in [0,2\pi]$ such that
        \[
            L_0 \coloneq \left\{(\theta_1,\psi_1) \in \R{\partial Y} \middle| \theta_1=0, \psi_1=c\right\}.
        \]
        Let $\{y_1,y_2\}$ be a second basis for $\fund{\partial Y}$.
        Let $(\theta_2,\psi_2)$ be the coordinates with respect to the basis $\{y_1,y_2\}$. We denote by $\mu'$ the weight computed using the second basis. Let $n_1,m_1,n_2,m_2$ be integers such that 
        \[
            x_1= n_1y_1+m_1y_2 \quad \text{and} \quad x_2= n_2y_1+m_2y_2.
        \]
        This implies that the following is a homeomorphism:
        \[
            f\colon \R{\partial Y} \to \R{\partial Y} \quad \text{with} \quad \left(\begin{matrix}
                \theta_1 \\
                \psi_1
            \end{matrix}\right) \mapsto \left(\begin{matrix}
                n_1 \theta_2 + m_1 \psi_2 \\
                n_2 \theta_2 + m_2 \psi_2
            \end{matrix}\right).
        \]
        The line $f(L_0)$ is of slope $\nicefrac{m_1}{m_2} \in \mathbb{Q} \cup \{\nicefrac{1}{0}\}$. We denote by $h_0$ the Lebesgue measure on $L_0$ and $h_f$ the Lebesgue measure on $f(L_0)$.
        Notice that $h_0 (L_0)=2\pi$ and $h_f(f(L_0))=2\pi \sqrt{m_1^2+m_2^2}$.
        Let $S \subset L_0$ be a measurable set, it is easy to see that
        \[
            h_f\left(f(S)\right) = \sqrt{m_1^2 + m_2^2} \cdot h_0 (S). 
        \]
        Thus,
        \[
            \mu'(f(L_0)) = \frac{h_f\left( H(Y)\cap f(L_0)\right)}{h_f \left( f(L_0) \right)} = \frac{\sqrt{m_1^2+m_2^2} \cdot h_0\left( H(Y)\cap f(L_0)\right)}{2\pi \sqrt{m_1^2+m_2^2}} = \mu(L_0).
        \]
    \end{proof}
\end{lemma}

\begin{exmp}
    The author proved in \cite{Mino} the following identities:
    \begin{itemize}
        \item $\mu \left( \mathbb{D}^2(\nicefrac{2}{1},\nicefrac{2}{1})\right)= 1$,
        \item $\mu \left( \mathbb{D}^2(\nicefrac{2}{1},\nicefrac{3}{1})\right)= \nicefrac{2}{3}$,
        \item $\mu \left( \mathbb{D}^2(\nicefrac{2}{1},\nicefrac{4}{q})\right)= \nicefrac{1}{2}$, where $q \in \mathbb{Z} \setminus 2\mathbb{Z}$.
    \end{itemize}
\end{exmp}

    \begin{lemma}\label{lemma: weight of D2(p1,p2)}
    Let $Y=\mathbb{D}^2 (\nicefrac{p_1}{q_1},\nicefrac{p_2}{q_2})$ with $2 \le p_1 \le p_2$. Then $\mu(Y) \ge \nicefrac{1}{2}$. Moreover, $\mu(Y)= \nicefrac{1}{2}$ if and only if $(p_1,p_2)=(2,4)$. Furthermore, if $(p_1,p_2) \neq (2,4)$, then $\mu(Y) \ge \nicefrac{2}{3}$.
    \begin{proof}
        This is a consequence of \cite[Lemma 7.6, Lemma 7.9, Lemma 7.15]{Mino}.
    \end{proof}
    \end{lemma}

    \begin{lemma}\label{lemma: eta in H(Y)-H(Y) in SFS with 2 sf}
    Let $Y=\mathbb{D}^2(\nicefrac{p_1}{q_1},\nicefrac{p_2}{q_2})$, where $p_i \ge 2$.
    If $\eta \in \R{\partial Y}$ is such that $\eta \in \overline{H(Y)} \setminus H(Y)$, then $\eta \in A(Y)$. Furthermore, if $\eta \in \overline{H(Y)} \setminus H(Y)$ and $\eta \in \{0,\pi\}^2 \subset \R{\partial Y}$, then $\eta \in A(Y)$ and $g \coloneq \gcd(p_1,p_2) \ge 2$.
    \begin{proof}
    $\eta \in \R{\partial Y}$ is such that $\eta \in \overline{H(Y)} \setminus H(Y)$, then $\eta \in A(Y)$ according to \cite[Proposition 7.2, Lemma 7.8]{Mino}.
    
        Without loss of generality, we can assume that $q_1$ and $q_2$ are odd. We consider $Y=Y_1 \cup_{\Sigma_1} C_3 \cup_{\Sigma_2} Y_2$, where $Y_1$ and $Y_2$ are solid tori.
        Therefore, for $i \in \{1,2\}$, we obtain that $A(Y_i)=T(Y_i,\partial Y_i) \subset \R{\Sigma_i}$.
        We consider $\fund{C_3}$ presented as in \eqref{eq: fundamental group C3} and, for $i \in \{1,2,3\}$, the space $\R{\Sigma_i}$ with coordinates $(\theta_i,\psi_i)$ as in \eqref{eq: coordinates for C3}.
        According to Theorem \ref{teo: alg for C3}, $H(Y)=H_{3,0}\cup H_{3,\pi}$. We recall that, for $i \in \{1,2\}$, Corollary \ref{cor: A(Y) with rational longitude} implies that
        \begin{align}
            \label{eq: per dim equazioni di A(Yi)}
            A(Y_i)= \left\{(\theta_i,\psi_i) \in \R{\Sigma_i} \middle| p_i \theta_i + q_i \psi_1 = 0 \right\}
        \end{align}
        Let $\eta$ be the representation with coordinates $(\varepsilon_1\pi,\varepsilon_2\pi) \in \R{\Sigma_3}$, where $\varepsilon_i \in \{0,1\}$, and let us assume that $\eta \in \overline{H(Y)} \setminus H(Y)$.

        As an application of Theorem \ref{teo: alg for C3}, if $\eta \in \overline{H(Y)}\setminus H(Y))$, then there exist $\theta_1,\theta_2 \in [0,2\pi]$ such that $\theta_i \notin \pi\mathbb{Z}$,
        \begin{align}
            \label{eq: la uso una volta, endpoints}
            \pm 2=\Tr \eta(x_1x_2) \in \partial I(2\cos \theta_1,2\cos \theta_2) = \left\{2\cos (\theta_1+\theta_2),2\cos (\theta_1-\theta_2)\right\},
        \end{align}
        and $(\theta_i,\varepsilon_2 \pi) \in A(Y_i) \cap \{\psi_i = \varepsilon_2\pi\}$. Here the interval $I$ is defined in \eqref{eq: interval I}.

        In particular, the identity in \eqref{eq: la uso una volta, endpoints} implies that either $\theta_1+\theta_2 = \pi \mathbb{Z}$ or $\theta_1-\theta_2 \in \pi \mathbb{Z}$. We notice that, if $(\theta_2,\varepsilon_2 \pi) \in A(Y_2)$, then according to Corollary \ref{cor: jewelled repr}, $(2\pi - \theta_2,\varepsilon_2 \pi) \in A(Y_2)$. Therefore, up to switching $\theta_2$ with $2\pi - \theta_2$, we can assume that $\theta_1+\theta_2 \in \pi \mathbb{Z}$.
        In particular, by Theorem \ref{teo: alg for C3}, $\eta \in A(Y_1) + A(Y_2) =A(Y)$.

        Let us assume that $\varepsilon_2=0$ and therefore $\eta \in \overline{H_{3,0}}\setminus H_{3,0}$. This implies that $\eta(h)=1$.
        The \eqref{eq: per dim equazioni di A(Yi)} implies that for $i \in \{1,2\}$,
        \[
            \theta_i = \frac{2\pi k_i}{p_i} \quad \text{with $k_i \in \{1,\cdots, p_i\}$}.
        \]
        As we stated before, $\theta_1+\theta_2 \in \pi \mathbb{Z}$.
        Without loss of generality, we can say that either $\theta_1+\theta_2 = 2 \pi$ or $\theta_1 + \theta_2 = \pi$.
        We study these two cases separately.
        If $\theta_1+ \theta_2 = 2 \pi$, then
        \[
            \theta_2= 2\pi - \theta_1 = 2\pi \frac{p_1 - k_1}{p_1} = \frac{2\pi k_2}{p_2}.
        \]
        Therefore, for $i \in \{1,2\}$, we obtain that $\theta_i = \frac{2\pi j_i}{g}$ where $g \coloneq\gcd(p_1,p_2)$ and $j_i \in \{1,\cdots, g \}$.
        Since $\theta_i \notin \pi \mathbb{Z}$, then $g \ge 3$.

        If $\theta_1+ \theta_2 = \pi$, then
        \[
            \theta_2= \pi - \theta_1 = 2\pi \frac{p_1 - k_1}{2p_1} = \frac{2\pi k_2}{p_2}.
        \]
        Therefore, $\theta_i = \frac{2\pi j_i}{g'}$ where $g' = \gcd(2p_1,p_2)$ and $j_i \in \{1,\cdots, g' \}$.
        As before, since $\theta_i \notin \pi \mathbb{Z}$, we get that $g' \ge 3$.
        If $g' \equiv_2 1$, then $g=g'$ and we get the conclusion. If $g' \equiv_2 0$, then $g' \ge 4$ and $g \ge 2$.

        Let us assume that $\varepsilon_2=1$ and therefore $\eta \in \overline{H_{3,\pi}}\setminus H_{3,\pi}$. This implies that $\eta(h)=-1$. Again, the \eqref{eq: per dim equazioni di A(Yi)}
        implies that for $i \in \{1,2\}$
        \[
            \theta_i = \frac{ \pi}{p_i} + \frac{2\pi k_i}{p_i} = 2\pi \frac{2k_i +1}{2p_i}  \quad \text{with $k_i \in \{1,\cdots, p_i\}$}.
        \]
        If $\theta_1+ \theta_2 \in \delta \pi$, with $\delta \in \{1,2\}$ then
        \[
            \theta_2 = \delta \pi - \theta_1 = 2\pi \frac{p_1 \delta - 2k_i - 1}{2p_1} = 2\pi \frac{2k_2 +1}{2p_2}.
        \]
        As before, this implies that $\theta_i = \frac{2\pi j_i}{g'}$ where $g' = \gcd(2p_1,2p_2)=2g$ and $j_i \in \{1,\cdots, g'\}$.  
        Since $\theta_i \notin \pi \mathbb{Z}$, then $2g \ge 3$. Hence $g \ge 2$.
    \end{proof}
\end{lemma}

\begin{cor}\label{cor: if eta in H(Y)-H(Y) then P(Y)}
    Let $Y$ be as in Lemma \ref{lemma: eta in H(Y)-H(Y) in SFS with 2 sf}.
    If $\eta \in \{0,\pi\}^2 \subset \R{\partial Y}$ and $\eta \in \overline{H(Y)} \setminus H(Y)$, then $\eta \in P(Y)$.
    \begin{proof}
        We are going to use the same notation of Lemma \ref{lemma: eta in H(Y)-H(Y) in SFS with 2 sf}.
        As we proved in Lemma \ref{lemma: eta in H(Y)-H(Y) in SFS with 2 sf}, if $\eta \in \overline{H(Y)} \setminus H(Y)$, then $\eta \in A(Y)$. In particular, we can assume $g \ge 2$. If $g \ge 3$, then the conclusion is a consequence of \cite[Lemma 6.4]{Mino}.

        Let us assume that $g =2$. We write $p_1=2n_1$ and $p_2=2n_2$ with $\gcd(n_1,n_2)=1$.
        As shown in \cite[Lemma 7.13]{Mino}, if $(\theta_3,\psi_3) \in \overline{H_{3,0}}\setminus H_{3,0}$, then $2 \cos (\theta_3) \neq \pm 2$. Since by hypothesis $\eta \in \{0,\pi\}^2$, this implies that
        $\eta$ cannot be in $\overline{H_{3,0}}\setminus H_{3,0}$. 
        Therefore we can assume without loss of generality that
        $\eta(h)=-1$ and $\eta \in \overline{H_{3,\pi}}\setminus H_{3,\pi}$.
        According to Theorem \ref{teo: alg for C3},
        there exist $\theta_1,\theta_2 \in [0,2\pi]$ such that $\theta_i \notin \pi\mathbb{Z}$,
        \begin{align*}
            \pm 2=\Tr \eta(x_1x_2) \in \partial I(2\cos \theta_1,2\cos \theta_2) = \left\{2\cos (\theta_1+\theta_2),2\cos (\theta_1-\theta_2)\right\},
        \end{align*}
        and $(\theta_i,\varepsilon_2 \pi) \in A(Y_i) \cap \{\psi_i = \varepsilon_2\pi\}$.
        As we said in the proof of Lemma \ref{lemma: eta in H(Y)-H(Y) in SFS with 2 sf}, we can assume that
        \[
            \theta_1+\theta_2 \in \pi \mathbb{Z}.
        \]
        Since $\eta(h)=-1$, then for $i \in \{1,2\}$, the \eqref{eq: per dim equazioni di A(Yi)} implies that
        \[
        \theta_i = \frac{\pi}{p_i} + \frac{2\pi k_i}{p_i}.
        \]
        Let $o_Y$ be the order of the rational longitude of $Y$.
        We recall, as a consequence of Lemma \cite[Lemma 5.1]{Mino} that $o_Y \equiv_2 0$ if and only if $p_1/2+p_2/2 \equiv_2 0$.
        Since we assumed that $\eta \in \{0,\pi\}^2$ and $\eta \in \overline{H_{3,\pi}}\setminus H_{3,\pi}$, then, without loss of generality, we can assume that $\theta_1+\theta_2= \delta \pi$, with $\delta \in \{1,2\}$.
        Thus,
        \[
            \theta_1+\theta_2 = \pi\frac{2k_1p_2+p_2+2k_2p_1+p_1}{p_1p_2}= \pi\frac{2k_1n_2+n_2+2k_2n_1+n_1}{2n_1n_2}= \delta \pi.
        \]
        If $\theta_1+\theta_2=\delta \pi$, then $2k_1n_2+n_2+2k_2n_1+n_1$ is a multiple of $2n_1n_2$ and in particular it is even. Hence $n_1+n_2=p_1/2+p_2/2 \equiv_20$ and then $o_Y=2$. We have just proven that if $\eta \in \overline{H_{3,\pi}}\setminus H_{3,\pi}$ and $g=2$, then $o_Y=2$ and $\eta(h)=-1$. The conclusion holds by \cite[Lemma 6.4]{Mino}.
    \end{proof}
\end{cor}

For the rest of the section we prove that if $Y$ is an $SU(2)$-abelian graph manifold rational homology $3$-sphere with first homology of small order, then $C_2$ does not appear in the decomposition of Corollary \ref{cor: the useful decomposition}.

\begin{lemma}\label{lemma: gamma non intersecta le linee orizzontali}
Let $Y=Y_1\cup_{\Sigma_1'}C_2$ be a graph manifold and $\alpha \subset \Sigma_2$ a slope such that $Y(\alpha)$ is a rational homology $3$-sphere. If there exists
irreducible representation $\rho \colon \fund{Y_1} \to SU(2)$ such that $\rho(h_1)=\pm 1$, where $h_1 \subset \Sigma_1 = \partial Y_1$ is the regular fiber of $C_2$, then $Y(\alpha)$ is not $SU(2)$-abelian.
.
\begin{proof}
    As usual, we parameterize $\R{\Sigma_1'}$ and $\R{\Sigma_2'}$ according to the the bases $\{x_1,h_1\}$ and $\{x_2,h_2\}$ as in Section \ref{sec: two important pieces}.
    Therefore, the restriction $\restr{\rho}{\fund{\Sigma_1'}}$ has coordinates $(\theta_1,\varepsilon) \in \R{\Sigma_1'}$, with  $\varepsilon \in \{0,\pi\}$.
    According to Theorem \ref{teo: alg for C2}, since $(\theta_1,\varepsilon) \in H(Y_1)$, the set $\{\psi_2 = \varepsilon \}$ is contained in $H(Y) \subset \R{\Sigma_2'}$.

    Let $n,m \in \mathbb{Z}$ be integers such that
    \[
        \alpha = nx_2+mh_2 \in \fund{\Sigma_2'}.
    \]
    By Lemma \ref{lemma: rational longitude of C2}, the regular fiber $h_2 \subset \Sigma_2'$ is the rational longitude of $Y$.
    Since $Y(\alpha)$ is a rational homology $3$-sphere, $n=\geo{\alpha}{h_2}$ is nonzero. Hence the line
    \[
    L_\alpha = \left\{ n \theta_2+ m \psi_2 =0 \right\} \subset \R{\Sigma_2'}
    \]
    and $\{\psi_2=\varepsilon \}$ have an intersection. This implies that $H(Y) \cap L_\alpha \neq \emptyset$.
    This implies that there exists an irreducible representation $\rho \colon \fund{Y} \to SU(2)$ that kills the slope $\alpha \in \fund{\Sigma_2}'$. 
    Therefore, since
    \[
        \fund{Y(\alpha)} = \frac{\fund{Y}}{\normalsubgroup{\alpha}},
    \]
    the representation $\rho$ gives an irreducible $SU(2)$-representation of $\fund{Y(\alpha)}$.
    \end{proof}
    \end{lemma}

\begin{figure}[t]
    \centering
    \begin{subfigure}[b]{0.4\textwidth}
    \centering
        \begin{tikzpicture}[>=latex,scale=1.45]
            \draw[->] (-.5,0) -- (2.5,0) node[right] {$\theta_1$};
            \foreach \x /\n in {1/$\pi$,2/$2\pi$} \draw[shift={(\x,0)}] (0pt,2pt) -- (0pt,-2pt) node[below] {\tiny \n};
            \draw[->] (0,-.5) -- (0,2.5) node[below right] {$\psi_1$};
            \foreach \y /\n in {1/$\pi$,2/$2\pi$}
            \draw[shift={(0,\y)}] (2pt,0pt) -- (-2pt,0pt) node[left] {\tiny \n};
            \node[below left] at (0,0) {\tiny $0$};
             \draw (0,0) rectangle (2,2);
             \draw (1,0)--(1,2);
             \draw (0,1)--(2,1);

             \draw [blue, very thick] (0,0)--(0,2);
             \draw [blue, very thick] (2,0)--(2,2);

            \draw[red, thick] (0,0) --(2,1);
            \draw[red, thick] (0,1) --(2,2);
        \end{tikzpicture}
        \caption{$T_1 \subset \R{\Sigma_1}$.}
        \label{Figure: teo per C2 sigma 1}
        \end{subfigure}
    \quad
    \begin{subfigure}[b]{0.4\textwidth}
    \centering
        \begin{tikzpicture}[>=latex,scale=1.45]
            \draw[->] (-.5,0) -- (2.5,0) node[right] {$\theta_2$};
            \foreach \x /\n in {1/$\pi$,2/$2\pi$} \draw[shift={(\x,0)}] (0pt,2pt) -- (0pt,-2pt) node[below] {\tiny \n};
            \draw[->] (0,-.5) -- (0,2.5) node[below right] {$\psi_2$};
            \foreach \y /\n in {1/$\pi$,2/$2\pi$}
            \draw[shift={(0,\y)}] (2pt,0pt) -- (-2pt,0pt) node[left] {\tiny \n};
            \node[below left] at (0,0) {\tiny $0$};
             \draw (0,0) rectangle (2,2);
             \draw (1,0)--(1,2);
             \draw (0,1)--(2,1);
             \draw [red, thick] (1,0)--(1,2);
             
             \draw[blue, very thick] (0,0)--(2,0);
             \draw[blue, very thick] (0,1)--(2,1);
             \draw[blue,  very thick] (0,2)--(2,2);
             
            \draw [red] (1,0) circle (1.5pt);
            \draw [red] (1,1) circle (1.5pt);
            \draw [red] (1,2) circle (1.5pt);

            \draw[red, thick] (1,1) --(2,0.5);
            \draw[red, thick] (1,0) --(0,0.5);

            \draw[red, thick] (1,2) --(2,1.5);
            \draw[red, thick] (1,1) --(0,1.5);

            \draw[cyan, dotted] (0,2)--(2,0); 
        \end{tikzpicture}
        \caption{$T_2 \subset \R{\Sigma_1}$.}
         \label{Figure: teo per C2 sigma 2}
    \end{subfigure}
    \caption{Proof of Proposition \ref{prop: C2 and o1=o2=1 then not SUa}. The reducible and irreducible points are in blue and red respectively. The circle indicates that the corresponding point is not in $H(Y)$.}
    \label{Figure: theorem for C2}
\end{figure}
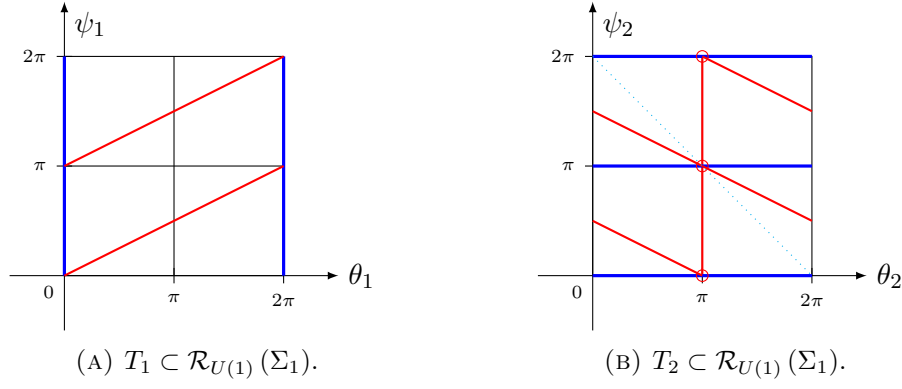

\begin{rmk}
    Let $\Sigma$ be a closed surface and we denote by $\Sigma_0$ the surface obtained by removing a small open disk from it. Let $Y\neq S^1 \times \mathbb{D}^2$ be a Seifert fibered manifold with torus boundary, that admits a Seifert fibration with $n$ singular fibers and whose base space is $\Sigma_0$. Let $\gamma$ be a slope in $\partial Y$. If $\gamma$ is a regular fiber for a Seifert fibration of $Y$, then $Y(\gamma)$ is a reducible manifold. If $\gamma$ is not a regular fiber, then $Y(\gamma)$ is a closed Seifert fibered manifold with base space $\Sigma$ and with $n+1$ singular fibers. 
    More precisely, $Y(\gamma)$ has $n$ singular fibers with the same orders as $Y$ and one additional singular fiber of order $\Delta(h,\gamma)$, where $h\subset \partial Y$ is a regular fiber.
    \label{rmk: order of the additional fiber of a dehn filling}
\end{rmk}

\begin{prop}\label{prop: C2 and o1=o2=1 then not SUa}
    Let $Y=Y_1\cup_{\Sigma_1'} C_2 \cup_{\Sigma_2'} Y_2$ be a rational homology $3$-sphere such that neither $Y_1$ not $Y_2$ is a solid torus. If $Y$ is a graph manifold such that $o_1=o_2=1$, then $Y$ is not $SU(2)$-abelian.
    \begin{proof}
        Since neither $Y_1$ nor $Y_2$ is a solid torus, they are both LIPA according to Theorem \ref{teo: admissible pieces have the property P}.
    
        We recall that, as a consequence of Lemma \ref{lemma: rational longitude of C2}, the rational longitude of $Y_1 \cup_{\Sigma_1'} C_2$ is $h_2 \subset \Sigma_2'$. Similarly, the rational longitude of $C_2 \cup_{\Sigma_2'} Y_2$ is $h_1 \subset \Sigma_1'$.
        Since $Y$ is a rational homology $3$-sphere, the manifold
        \[
            N \coloneqq C_2(\Sigma_1',\Sigma_2'; \lambda_1,\lambda_2)
        \]
        is a rational homology $3$-sphere. In particular, $N$ is a Seifert fibered space with base space $\mathbb{RP}^2$.
        If $N$ is not $SU(2)$-abelian, then $Y$ is not $SU(2)$-abelian by Corollary \ref{cor: Y SU(2) allora i due pezzi sono SU(2)-abelian}. Without loss of generality, we can assume that $N$ is $SU(2)$-abelian.
        According to \cite[Proposition 3.5]{Menangerie} and Remark \ref{rmk: order of the additional fiber of a dehn filling}, $N$ is $SU(2)$-abelian if and only if it is either a lens space or $\mathbb{RP}^3 \# \mathbb{RP}^3$. Therefore, $N$ is $SU(2)$-abelian if and only if it has at most one singular fiber. As an application of Remark \ref{rmk: order of the additional fiber of a dehn filling}, either $\geo{\lambda_1}{h_1}=1$ or $\geo{\lambda_2}{h_2}=1$.
        Up to swapping $Y_1$ with $Y_2$, we assume that $\geo{\lambda_1}{h_1}=1$. We consider $\{\lambda_1,h_1\}$ as a basis for $\fund{\Sigma_1'}$ and we give $(\theta_1,\psi_1)$ coordinates $\R{\Sigma_1'}$.
        We give $\R{\Sigma_2'}$ coordinates $(\theta_2,\psi_2)$ as explained in Remark \ref{rmk: bases for C2}.
        
        Let $\gamma \colon (-1,1) \to H(Y_1)$ be a LIPA-path of $Y_1$. Let $\Gamma=\ima \gamma \subset H(Y_1)$ be its support. Let $\Gamma^t$ be the set $\Gamma$ translated by the vector $(0,\pi)$. Explicitly,
        \[
            \Gamma^t \coloneq \left\{(\theta_1,\psi_1 + \pi) \in \R{\Sigma_1} \middle| (\theta_1,\psi_1) \in \Gamma \right\}.
        \]
        Lemma \ref{lemma: translation} implies that $\Gamma^t \subset H(Y_1)$. We define $T_1 \subset T(Y_1,\partial Y_1)$ as
        \[
            T_1=A(Y_1)\cup\Gamma\cup\Gamma^t \subset\R{\Sigma_1}.
        \]
        Since $o_1=1$, the space $A(Y_1)$ equals $\{\theta_1=0\} \subseteq\R{\Sigma_1}$.
        We recall that, since $Y=Y_1\cup_{\Sigma_1'}C_2 \cup_{\Sigma_2'}Y_2$ is a rational homology $3$-sphere, then $\geo{\lambda_2}{h_2}\neq0$ by Lemma \ref{lemma: rational longitude of C2}.
        If either $\Gamma$ or $\Gamma'$ intersects $\{\psi_1=0\}\cup \{\psi_1=\pi\}$, then the manifold $Y_1 \cup_{\Sigma_1}'C_2(\Sigma_2';\lambda_2)$ is not $SU(2)$-abelian as a consequence of Lemma \ref{lemma: gamma non intersecta le linee orizzontali} and this implies that $Y=Y_1\cup_{\Sigma_1'} C_2 \cup_{\Sigma_2'} Y_2$ is not $SU(2)$-abelian. Therefore we assume that neither $\Gamma$ nor $\Gamma'$ intersects $\{\psi_1=0\}\cup \{\psi_1=\pi\}$. Figure \ref{Figure: teo per C2 sigma 1} shows $T_1 \subset T(Y_1,\partial Y_1)\subset \R{\Sigma_1}$.
        Theorem \ref{teo: alg for C2} defines a subset $T_2$ of $H(Y_1\cup_{\Sigma_1'}C_2)\subset \R{\Sigma_2}$ corresponding to $T_1$. Figure \ref{Figure: teo per C2 sigma 2} shows the set $T_2$ and $A(Y_2)$.

        We define $\partial \Gamma\subset \R{\Sigma_1'}$ as the two points $\overline{\Gamma} \setminus \Gamma \subset \R{\Sigma_1'}$. We define $\partial \Gamma^t\subset \R{\Sigma_1'}$ similarly.
        We recall that $\{0,\pi\}^2 \subset \R{\Sigma_1'}$ is the set of the four ``central" points of $\R{\Sigma_1'}$.
        By definition of $\Gamma^t$ we obtain that
        \[
            \left| \partial \Gamma \cap \{0,\pi\}^2\right| = \left| \partial \Gamma^t \cap \{0,\pi\}^2\right|.
        \]
        In other words, $\Gamma$ has either two, one, or none endpoints in $\{0,\pi\}^2 \subset \R{\Sigma_1'}$ if and only if $\Gamma^t$ has the same number of endpoints in $\{0,\pi\}^2 \subset \R{\Sigma_1'}$.
        Figure \ref{Figure: teo per C2 sigma 2} shows that
        if the two points $\partial \Gamma$
        are not both in $\{0,\pi\}^2$, then the only straight line of $\R{\Sigma_2'}$ containing the origin that does not intersect $T_2 \subseteq H(Y_1\cup_{\Sigma_1'}C_2) \subset \R{\Sigma_2'}$ is $\{\psi_2=0\} \subset \R{\Sigma_2'}$. This implies that if if the two points $\partial \Gamma$
        are not both in $\{0,\pi\}^2$ and $\alpha \subset\Sigma_2'$ is a slope different from the regular fiber $h_2 \subset \Sigma_2'$, then $T_2$ has a nonempty intersection with the line $L_\alpha$,
        where
        \[
            L_\alpha \coloneq \left\{\eta \in \R{\Sigma_2'} \middle| \eta(\alpha)=1\right\}.
        \]
        Therefore, if $\alpha \subset\Sigma_2'$ is a slope as above, Theorem \ref{teo: M SU(2)-abeliano se e solo se i pezzi sono empty} implies that
        the manifold
        \[
            Y_1 \cup_{\Sigma_1'}C_2(\Sigma_2';\alpha),
        \]
        is not $SU(2)$-abelian. Therefore, if the two points $\partial \Gamma$
        are not both in $\{0,\pi\}^2$, then $Y=Y_1\cup_{\Sigma_1'}C_2 \cup_{\Sigma_2'}Y_2$ is not $SU(2)$-abelian by Lemma \ref{lemma: rational longitude of C2} and Corollary \ref{cor: Y SU(2) allora i due pezzi sono SU(2)-abelian}.
        
        Without loss of generality, we assume that both points in $\partial \Gamma$ are in $\{0,\pi\}^2 \subset \R{\Sigma_1'}$.
        Figure \ref{Figure: teo per C2 sigma 2} also represents in light blue and dotted the line
        \[
           L_{\nicefrac{1}{-1}}= \left\{(\theta_2,\psi_2) \in \R{\Sigma_2'} \middle| \theta_2-\psi_2 =0\right\}.
        \]
        Since $o_2=1$ and as a consequence of Corollary \ref{cor: A(Y) with rational longitude}, the set $A(Y_2)$ is equal to the line $L_{\nicefrac{1}{-1}}$ above if and only if $\lambda_2= x_2 -h_2$ as a slope of $\Sigma_2'$.
        
        As shown in Figure \ref{Figure: teo per C2 sigma 2}, if $\lambda_2$ is not equal to $x_2 -h_2$ as a slope of $\Sigma_2'$, then $T_2 \subseteq H(Y_1\cup_{\Sigma_1'}C_2)$ has an intersection with $A(Y_2)$.
        Moreover, Figure \ref{Figure: teo per C2 sigma 2} shows that every line of $\R{\Sigma_2'}$ connecting two points of the line $L_{\nicefrac{1}{-1}}$ intersects $T_2 \cup A(Y_1\cup_{\Sigma_1'}C_2)$.
        This implies that either $H(Y_1\cup_{\Sigma_1'}C_2)\cap A(Y_2) \neq \emptyset$ or a LIPA-path of $Y_2$ intersect $H(Y_1\cup_{\Sigma_1'}C_2) \cup A(Y_1\cup_{\Sigma_1'}C_2)$ and  therefore 
        \[
        \left(A(Y_1\cup_{\Sigma_1'}C_2)\cup H(Y_1\cup_{\Sigma_1'}C_2)\right) \cap H(Y_2) \neq \emptyset.
        \]
        In both cases, the conclusion is given by Theorem \ref{teo: M SU(2)-abeliano se e solo se i pezzi sono empty}.
    \end{proof}
\end{prop}

\begin{lemma}\label{lemma: C2 implies sigma 16}
    Let $Y=Y_1\cup_{\Sigma_1'} C_2 \cup_{\Sigma_2'} Y_2$ be a rational homology $3$-sphere such that neither $Y_1$ not $Y_2$ is a solid torus. If $Y$ is $SU(2)$-abelian graph manifold, then $|H_1(Y;\mathbb{Z})| \ge 16$.
    \begin{proof}
    As a consequence of Proposition \ref{prop: C2 and o1=o2=1 then not SUa}, either $o_1 \ge 2$ or $o_2 \ge2$. We assume that $o_1 \ge 2$. Remark \ref{rmk: order of toroidal} implies that the order of the torsion subgroup of $H_1(Y_1;\mathbb{Z})$ is divided by $o_1$.

    Let us focus on $C_2 \cup_{\Sigma_2'}Y_2$. Lemma \ref{lemma: rational longitude of C2} and Corollary \ref{cor: C2 obstruction longitudes} imply that $h_1 \in \fund{\Sigma_1'}$ is the rational longitude of $C_2 \cup_{\Sigma_2'}Y_2$ and it has order two. 
    Let us denote by $t_2'$ the order of the torsion subgroup of $H_1(C_2 \cup_{\Sigma_2'}Y_2;\mathbb{Z})$. Since the order of the rational longitude of $C_2 \cup_{\Sigma_2'}Y_2$ is two, $t_2'$ is divided by two and, in particular, $t_2' \ge 2$ and $o_2t_2' \ge 4$.
    Therefore, Remark \ref{rmk: order of toroidal} applied to the torus $\Sigma_1' \subset Y$ implies that
    \[
        \left| H_1\left(Y_1 \cup_{\Sigma_1'}C_2 \cup_{\Sigma_2'}Y_2;\mathbb{Z}\right) \right| = \geo{\lambda_1}{h_1}o_2t_2't_1o_1 \ge 4o_1^2 \ge 16.
    \]
    \end{proof}
\end{lemma}

\begin{lemma}\label{lemma: the last dec}
    Let $Y$ be an $SU(2)$-abelian graph manifold rational homology $3$-sphere with nontrivial JSJ decomposition. If $|H_1(Y;\mathbb{Z})| \le 7$, then there exists a system of disjoint tori $\tau \subset Y$ such that
    \[
        \overline{Y \setminus \tau} = \left( \bigcup C_3\right) \cup \left(\bigcup \mathbb{D}^2(\nicefrac{p_\ast}{q_\ast})\right),
    \]
    where $\mathbb{D}^2(\nicefrac{p_\ast}{q_\ast})$ is a Seifert fibered manifold with torus boundary and disk base space.
    \begin{proof}
    According to Corollary \ref{cor: the useful decomposition} and Lemma \ref{lemma: C2 implies sigma 16}, there exists a system of disjoint tori $\tau \subset Y$ such that
    \begin{align}
    \label{eq: una volta dec senza C2}
        \overline{Y \setminus \tau} = \left( \bigcup C_3\right) \cup \left(\bigcup M_i\right),
    \end{align}
    where $M_i$ is a Seifert fibered manifold with torus boundary.
    Note that in \eqref{eq: una volta dec senza C2} we did not assume that the manifold $M_i$ has a disk base space: a priori, according to Corollary \ref{cor: the useful decomposition}, $M_i$ could be of the form $C_2 \cup S^1 \times \mathbb{D}^2$, which is a Seifert fibered space with a punctured $\mathbb{RP}^2$ base space. In fact, we shall prove that this cannot occur by showing that $M_i$ cannot have such a base space.

    Let $\Sigma \in \tau$ be a torus that splits 
    $Y$ into a graph manifold $Y_1$ and a Seifert fibered manifold $Y_2$. Thus,
    $Y=Y_1 \cup_\Sigma Y_2$. We prove the conclusion by showing that $Y_2$ admits a fibration over a disk.
    Let $\sigma \in\mathbb{N}$ be the order of $H_1(Y;\mathbb{Z})$. Since $Y$ is a rational homology $3$-sphere, then $\sigma \neq 0$.
        According to Remark \ref{rmk: order of toroidal},
        \begin{align*}
            \sigma \coloneq \left|H_1(Y;\mathbb{Z})\right|=o_1o_2t_1t_2 \geo{\lambda_1}{\lambda_2},
        \end{align*}
    where $o_i$ is the order of $\lambda_i$ in $H_1(Y_i,\mathbb{Z})$ and $t_i$ is the order of the torsion subgroup of $H_1(Y_i,\mathbb{Z})$. Since $Y$ is a rational homology $3$-sphere, $\geo{\lambda_1}{\lambda_2} \ge 1$.

    As explained in Section \ref{section: the topology set up}, $Y_2$ has base space that is either a disk or a punctured $\mathbb{RP}^2$.
    We give the conclusion by showing that if $Y_2$ does not admit a Seifert fibration over a disk, then $|H(Y;\mathbb{Z})| \ge 8$.
    
    Let us assume that $Y_2$ does not admit a Seifert fibration with disk base space, and therefore that admits a Seifert fibration whose base space is a punctured $\mathbb{RP}^2$.
    
    If $Y_2$ has no singular fibers, then $Y_2$ admits a Seifert fibration over a disk as shown in \cite[Proposition 10.4.16]{martelli}.
    Therefore, we can assume that $Y_2$ has one or more fiberrs.
    This implies that
    \[
        \fund{Y_2} = \langle x_1,\cdots,x_n,y,z,h \,|\, [x_i,h], x_1\cdots x_nyz^2, zhz^{-1}h, x_i^{p_i}h^{q_i}  \rangle.
    \]
    where $n\ge 1$, $p_i \ge 2$, and $\gcd(p_i,q_i)=1$. Details can be found in \cite[Section 2.2]{boyer2017foliations}.
    It is straightforward to see that $h$ coincide with the rational longitude of $Y_2$ and it has order $2$.
    Therefore,
    \[
        H_1(Y_2;\mathbb{Z}) = \coker \begin{bmatrix}
            p_1 & 0 & \cdots & 0 & 0 & 1\\
            0 & p_2 & \cdots & \cdots & \cdots & \cdots \\
             \cdots & \cdots & \cdots& \cdots& \cdots & \cdots \\
            0 & 0 & \cdots& p_n & 0 & 1 \\
            q_1 & q_2 & \cdots & q_n& 2& 0 \\
            0 & \cdots & \cdots & \cdots & 0 &1 \\
            0 & 0 & \cdots& 0& 0& 2 
        \end{bmatrix}.
    \]
    Thus, the Smith algorithm applied to the matrix above shows that the torsion subgroup of $H_1(Y_2;\mathbb{Z})$ has order $t_2=2 p_1 \cdots p_n \ge 2^{n+1}$.
    
    Therefore $o_2=2$ and $t_2 \ge 4$. Remark \ref{rmk: order of toroidal} applied to $Y=Y_1\cup_{\Sigma} Y_2$ implies that $\sigma = |H_1(Y;\mathbb{Z})| \ge 8$. This gives the conclusion. 
    \end{proof}
\end{lemma}

\section{The main result}
In this section we prove Theorem \ref{teo: sigma small then not SU(2)-abelian} that positively answers \cite[Conjecture 1.5]{InstantoLSpaceAndSplicing} in the case of graph manifolds.

In the following $Y_2$ is a Seifert fibered space with disk base space and two singular fibers. Therefore $Y_2$ is of the form $\mathbb{D}^2(\nicefrac{p_1}{q_1},\nicefrac{p_2}{q_2})$ with $2 \le p_1 \le p_2$. Following the notation of \cite{Mino}, we denote by $g_2$ the greatest common division of the orders of these fibers, namely $g_2=\gcd(p_1,p_2)$. We recall that $g_2$ equals the cardinality of the torsion subgroup of $H_1(Y_2;\mathbb{Z})$.

\begin{cor}\label{cor: sigma=2 o1=delta=1 and o2=2}
    Let $Y_1$ be a graph manifold rational homology solid torus and $Y_2$ a Seifert fibered space with two singular fibers and disk base space. If $Y=Y_1 \cup_{\Sigma} Y_2$ is a rational homology $3$-sphere such that $\Delta(\lambda_1,\lambda_2) = 1$, $g_2=o_2=1$ and $o_1=2$, then $Y$ is not $SU(2)$-abelian.
    \begin{proof}
        Since $\geo{\lambda_1}{\lambda_2}=1$ the set $\{\lambda_1,\lambda_2\}$ is a basis for $\fund{\Sigma}$.
        We give to $\R{\Sigma}$ coordinates $(\theta,\psi)$ with respect to this ordered basis as in \eqref{eq: general coordinates}. Since $o_1=2$ and $o_2=1$, $A(Y_1)= \{\theta=0\} \cup \{\theta=\pi\}$ and $A(Y_2)=\{\psi=0\}$ by Corollary \ref{cor: A(Y) with rational longitude}.

        Let $\gamma_2$ be a LIPA-path and $\gamma_2^t$ its translation of $(\pi,0)$.
        Lemma \ref{lemma: translation} implies that $\gamma_2^t$ is a LIPA-path for $Y_2$.
        With an abuse of notation, we identify $\ima{\gamma_2}$ and $\ima{\gamma_2^t}=\ima{\gamma_2}+(\pi,0)$ with $\gamma_2$ and $\gamma_2^t$ respectively.
        The path $\gamma_2$ connects two points of $A(Y_2)$, namely
        \[
            (\theta_1,0) \in A(Y_2) \quad \text{and} \quad (\theta_2,0) \in A(Y_2).
        \]
        Therefore, $\gamma_2^t$ connects the points
        \[
            (\theta_1 +\pi,0) \in A(Y_2) \quad \text{and} \quad (\theta_2+\pi,0) \in A(Y_2).
        \]
        If either $\gamma_2$ or $\gamma_2^t$ intersects $A(Y_1)$, then $A(Y_1) \cap H(Y_2) \neq \emptyset$. As an application of Theorem \ref{teo: M SU(2)-abeliano se e solo se i pezzi sono empty}, the manifold $Y=Y_1 \cup_{\Sigma}Y_2$ is not $SU(2)$-abelian.
        
        Without loss of generality, we can assume that neither $\gamma_2$ nor $\gamma_2^t$ intersects $A(Y_1)$. As we saw in the proof of Proposition \ref{prop: sigma =1 ish}, this implies that either $\theta_1,\theta_2 \in [0,\pi]$ or $\theta_1,\theta_2 \in [\pi,2\pi]$.
        As an application of \cite[Lemma 7.1]{Mino}, we have that $P(Y_2)= \emptyset$. By Corollary \ref{cor: if eta in H(Y)-H(Y) then P(Y)}, we have that the accumulation points of $\gamma_2$ (resp. $\gamma_2^t$) are disjoint from $\{0,\pi\}^2 \subset \R{\Sigma}$. Therefore, either $\theta_1,\theta_2 \in (0,\pi)$ or $\theta_1,\theta_2 \in (\pi,2\pi)$.
        Figure \ref{Figure o1=1 e o2=2} shows $\gamma_2\cup\gamma_2^t$, $A(Y_1)$ and $A(Y_2)$ in the chosen coordinates.
        In particular, Figure \ref{Figure o1=1 e o2=2} illustrates that a LIPA-path of $Y_1$ intersects either $A(Y_2)$, $\gamma_2$, or $\gamma_2^t$. Therefore, 
        \[
        H(Y_1) \cap (H(Y_2) \cup A(Y_2)) \neq \emptyset.
        \]
        The conclusion is given by Theorem \ref{teo: M SU(2)-abeliano se e solo se i pezzi sono empty}.
    \end{proof}
\end{cor}

\begin{rmk}\label{rmk: change of basis}
Let $\Sigma$ be a torus and $\{x_0,y\}$ a basis for $\fund{\Sigma}$. Let $m,n \in \mathbb{Z}$ be two coprime nonzero integers. We define $\alpha \in \fund{\Sigma}$ as
\[
    \alpha = m x_0 + n y.
\]
For every $k \in \mathbb{Z}$, the set $\{x_k,y\}$ is a basis for $\fund{\Sigma}$, where $x_k = x_0 + k y$. In this basis,
\[
    \alpha = m x_k + (n-mk) y.
\]
In particular, there exists a $k \in \mathbb{Z}$ such that $-|m| \le -|n-mk| \le 0$
\end{rmk}

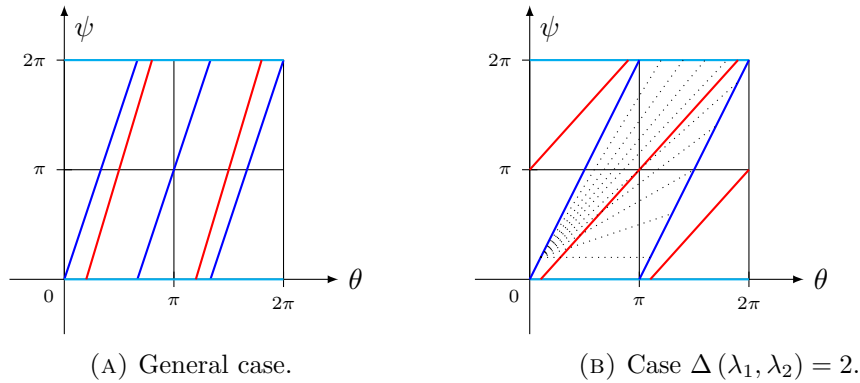
\begin{figure}[b]
    \centering
    \begin{subfigure}[b]{0.4\textwidth}
    \centering
        \begin{tikzpicture}[>=latex,scale=1.45]
            \draw[->] (-.5,0) -- (2.5,0) node[right] {$\theta$};
            \foreach \x /\n in {1/$\pi$,2/$2\pi$} \draw[shift={(\x,0)}] (0pt,2pt) -- (0pt,-2pt) node[below] {\tiny \n};
            \draw[->] (0,-.5) -- (0,2.5) node[below right] {$\psi$};
            \foreach \y /\n in {1/$\pi$,2/$2\pi$}
            \draw[shift={(0,\y)}] (2pt,0pt) -- (-2pt,0pt) node[left] {\tiny \n};
            \node[below left] at (0,0) {\tiny $0$};
             \draw (0,0) rectangle (2,2);
             \draw (1,0)--(1,2);
             \draw (0,1)--(2,1);

             \draw [blue, thick] (0,0)--(2/3,2);
             \draw [blue, thick] (2/3,0)--(4/3,2);
             \draw [blue, thick] (4/3,0)--(2,2);
             \draw [cyan, thick] (0,0)--(2,0);
              \draw [cyan, thick] (0,2)--(2,2);
             \draw [red, thick] (0.2,0)--(0.8,2);
             \draw [red, thick] (1.2,0)--(1.8,2);
        \end{tikzpicture}
        \caption{General case.}
        \label{Figure: sigma small at most 11}
    \end{subfigure}
    \quad
    \begin{subfigure}[b]{0.4\textwidth}
        \begin{tikzpicture}[>=latex,scale=1.45]
            \draw[->] (-.5,0) -- (2.5,0) node[right] {$\theta$};
            \foreach \x /\n in {1/$\pi$,2/$2\pi$} \draw[shift={(\x,0)}] (0pt,2pt) -- (0pt,-2pt) node[below] {\tiny \n};
            \draw[->] (0,-.5) -- (0,2.5) node[below right] {$\psi$};
            \foreach \y /\n in {1/$\pi$,2/$2\pi$}
            \draw[shift={(0,\y)}] (2pt,0pt) -- (-2pt,0pt) node[left] {\tiny \n};
            \node[below left] at (0,0) {\tiny $0$};
             \draw (0,0) rectangle (2,2);
             \draw (1,0)--(1,2);
             \draw (0,1)--(2,1);

             \draw [cyan, thick] (0,0)--(2,0);
             \draw [cyan, thick] (0,2)--(2,2);

             \draw [blue, thick] (0,0)--(1,2);
             \draw [blue, thick] (1,0)--(2,2);


            \draw [red, thick] (0.1,0)--(1.9,2);
            \draw [red, thick] (1.1,0)--(2,1);
            \draw [red, thick] (0,1)--(0.9,2);

            \draw[dotted] (0.1,0.2)--(1.9,1.8);
            \draw[dotted] (0.1,0.2)--(1.7,1.4);
            \draw[dotted] (0.1,0.2)--(1.5,1);
            \draw[dotted] (0.1,0.2)--(1.3,0.6);
            \draw[dotted] (0.1,0.2)--(1.1,0.2);

            \draw[dotted] (0.1,0.2)--(1.8,2);
            \draw[dotted] (0.1,0.2)--(1.6,2);
            \draw[dotted] (0.1,0.2)--(1.4,2);
            \draw[dotted] (0.1,0.2)--(1.2,2);
             
        \end{tikzpicture}
        \caption{Case $\geo{\lambda_1}{\lambda_2}=2$.}
        \label{Figure. Case d=2, c=pi}
    \end{subfigure}
    \caption{In Cyan $A(Y_1)$, in red the LIPA path for $Y_1$. The set in blue is $A(Y_2)$.}
    \end{figure}

\begin{lemma}\label{lemma: o2g2o1=1 and sigma <= 5}
    Let $Y_1$ be a graph manifold rational homology solid torus and $Y_2$ a Seifert fibered space with disk base space and two singular fibers.
    Let us further assume that $Y_1 \neq S^1 \times \mathbb{D}^2$, and let $Y=Y_1 \cup_{\Sigma} Y_2$.
    If $3 \le  \Delta(\lambda_1,\lambda_2) \le 5$ and $o_1=o_2=g_2=1$, then $Y$ is not $SU(2)$-abelian.
     \begin{proof}
     The conclusion requires three claims that are proven at the end of this proof.
        Let $\gamma \colon (-1,1)\to \R{\Sigma}$ be a path.
        With an abuse of notation, we identify $\gamma$ with its image in $\R{\Sigma}$.
            
        Let $\xi \subset$ be a slope such that $\{\xi,\lambda_1\}$ is an ordered basis for $\fund{\Sigma}$.
        We give the space $\R{\Sigma}$ coordinates $(\theta,\psi)$ as in \eqref{eq: general coordinates} with respect to this ordered basis.
        Let $n \in \mathbb{Z}$ be an integer such that
        \[
            \lambda_2= \geo{\lambda_1}{\lambda_2} \xi+ n \lambda_1 \subset \Sigma.
        \]
        According to Remark \ref{rmk: change of basis}, we can suppose without loss of generality that $ -\geo{\lambda_1}{\lambda_2}<n<0$.
        In this case, we write $\lambda_2=\nicefrac{\geo{\lambda_1}{\lambda_2}}{n}$.
        Corollary \ref{cor: A(Y) with rational longitude} implies that
        \[
            A(Y_2) = \left\{ \geo{\lambda_1}{\lambda_2} \theta+ n \psi =0\right\} \subset \R{\Sigma} \quad \text{and} \quad A(Y_1) = \left\{ \psi=0 \right\} \subset \R{\Sigma}.
        \]
        
        Since $Y_1$ is not a solid torus, then it is LIPA by Theorem \ref{teo: admissible pieces have the property P}.
        Let $\gamma_1 \colon (-1,1) \to H(Y_1)$ be a LIPA path for $Y_1$, we denote by $\gamma_1^t$ the path $\gamma_1$ translated by $(\pi,0)$. The path $\gamma_1^t$ is a LIPA-path for $Y_1$ by Lemma \ref{lemma: translation}. Figure \ref{Figure: sigma small at most 11} shows $\R{\Sigma}$, $A(Y_1)$, $A(Y_2)$ in the case $\lambda_2=-\nicefrac{3}{1}$, $\gamma_1$, and $\gamma_1^t$ in the chosen coordinates.

        If $A(Y_2)$ intersect $\gamma_1 \cup \gamma_1^t$, then $H(Y_1) \cap A(Y_2) \neq \emptyset$ and $Y$ is not $SU(2)$-abelian by Theorem \ref{teo: M SU(2)-abeliano se e solo se i pezzi sono empty}. Therefore, we assume that $A(Y_2)$ does not intersect the paths $\gamma_1$ and $\gamma_1^t$.

       Let $\mu$ be the weight as in Definition \ref{degn weight mu}.
       By Lemma \ref{lemma: weight of D2(p1,p2)}, $\mu(Y_2) \ge \nicefrac{2}{3}$. Let $L \subset \R{\Sigma}$ be the line containing a LIPA-path of $Y_2$. Let $\nicefrac{p}{q} \in \mathbb{Q} \cup \{\nicefrac{1}{0}\}$ be the slope of $L$. We recall that this means that there exists a $c \in [0,2\pi]$ such that
       \begin{equation}
            \label{eq: una volta Linea del LIPA Y2}
            L=\left\{(\theta,\psi) \in \R{\Sigma} \middle| p \theta+q \psi = c \right\}.
        \end{equation}
        Without loss of generality, we assume that $L$ is the line that realizes the maximal weight $\mu(Y_2)$, namely
        \[
            \mu(L)= \frac{h(L \cap H(Y_2))}{h(L)}= \mu(Y_2) \ge \frac{2}{3}.
        \]
        Here $h$ is the Lebesgue measure on $L$ induced by the coordinates $(\theta,\psi)$ as in Section \ref{sec. utensilis}.
        
        As a consequence of Fact \ref{Fact: centralizers}, if $\eta \colon \fund{\Sigma} \to \Ui$ is a representation in
        $H(Y_2)$, then $\eta(h_2) \in \{\pm 1\}$ where $h_2 \subset \partial Y_2$ is a regular fiber of $Y_2$. This implies that if $\eta' \colon \fund{\Sigma} \to \Ui$ is a representation in $L$, then $\eta'(h_2) \in \{\pm 1\}$. Thus, $L$ contains two points of $\{0,\pi\}^2 \subset \R{\Sigma}$. In other words, either $c =0$ or $c=\pi$.
        
        We define 
        \begin{equation}\label{eq: defn S}
            S \coloneqq \{\nicefrac{1}{0},\nicefrac{0}{1}, \pm \nicefrac{1}{1}, \pm \nicefrac{2}{1}, \pm \nicefrac{1}{2}\} \subseteq \mathbb{Q} \cup \{\nicefrac{1}{0}\},
        \end{equation}
        \begin{claim}\label{claim: se non in S allora non SU(2)a}
            Let $\nicefrac{p}{q}$ be the slope of $L$.
            If $\nicefrac{p}{q} \notin S$, then $Y=Y_1\cup_{\Sigma} Y_2$ is not $SU(2)$-abelian.
        \end{claim}

        We prove the conclusion by showing that for every  $\nicefrac{p}{q}\in S$, either the slope of the line $L$ is not $\nicefrac{p}{q}$ or, if it does, $L \cap H(Y_2)$ intersects $A(Y_1) \cup H(Y_1)$. In the latter case, $Y=Y_1\cup_{\Sigma}Y_2$ is not $SU(2)$-abelian by Theorem \ref{teo: M SU(2)-abeliano se e solo se i pezzi sono empty}.
        
       According to Corollary \ref{cor: A(Y) with rational longitude}, the space $A(Y_2) \subset\R{\Sigma}$ is a line of slope $\geo{\lambda_1}{\lambda_2}/n$.
       We recall that a LIPA-path of $Y_2$ starts and ends at $A(Y_2)$ and it is not contained in $A(Y_2)$. Therefore, since $L$ contains a LIPA-path of $Y_2$, the line $L$ contains at least one point of $A(Y_2)$. We define $m \in \mathbb{Z}$ as
       \[
       m \coloneq |A(Y_2) \cap L| = |pn-\geo{\lambda_1}{\lambda_2}q|= m \ge 1.
       \]
       We remark that $L$ is not contained in $A(Y_2)$ by definition of LIPA-path. Therefore $|A(Y_2) \cap L|< \infty$ and the integer $m$ is well-defined.
       In particular, the points of the intersection $A(Y_2) \cap L$ are equally spaced in $L$.
        Therefore $A(Y_2)$ divides the line $L$ in $m$ segments $I_1, \cdots, I_m$ of equal length:
        \[
            \frac{h(I_j)}{h(L)}= \frac{1}{m} \quad \text{and} \quad L = \bigcup_{i=1}^m I_i.
        \]
        By Lemma \ref{lemma: eta in H(Y)-H(Y) in SFS with 2 sf}, $\overline{H(Y_2)} \setminus H(Y_2) \subset A(Y_2)$.
        Thus, the intersection $L\cap H(Y_2)$ is such that
        \[
             \overline{L\cap H(Y_2)} \setminus L\cap H(Y_2) \subset A(Y_2).
        \]
        Therefore $\overline{L\cap H(Y_2)}$ is the union of segments $I_i$.
        
        Let $l \in \mathbb{N}$ be defined as
        \[
            l \coloneqq \left| L \cap A(Y_2) \cap \{0,\pi\}^2\right| \in \mathbb{N}.
        \]
        In other words, the lines $L$ and $A(Y_2)$ share $l$ ``central" points $\{0,\pi\}^2\subset \R{\Sigma}$.
        \begin{claim}\label{claim mu(m-2-l/m)}
            $\mu(L) \le \max\left\{0, \frac{m-2-l}{m}\right\}$.
        \end{claim}
        We recall that the line $L$ admits the description as in \eqref{eq: una volta Linea del LIPA Y2} with $c \in \{0,\pi\}$.
       The quantity $\frac{m-2-l}{m}$ depends of the choice of $c \in \{0,\pi\}$ as the integer $l$ depends on $c$.

       \begin{claim}\label{claim: togliere segments}
            If a segment $I_j$ contains a point in $\{0,\pi\}^2$, then it is not contained in $\overline{L \cap H(Y_2)}$.
       \end{claim}
        
       We now give an idea of the resolution procedure. Given a $n \in \{-\geo{\lambda_1}{\lambda_2}-1,\cdots,-1\}$ with $\gcd(n,\geo{\lambda_1}{\lambda_2})=1$ we consider 
       \[
       \lambda_2=\frac{\geo{\lambda_1}{\lambda_2}}{n} \coloneqq \geo{\lambda_1}{\lambda_2} \xi + n \lambda_1 \subset \Sigma.
       \]
       For each slope in $S$, for each $n$ as above and $c \in \{0,\pi\}$ we compute the quantity $\frac{m-2-l}{m}$, if this is smaller than $\nicefrac{2}{3}$ then, since 
       \[
       \mu(L)=\mu(Y_2) \ge \frac{2}{3}, \quad \text{and} \quad \mu(Y_2) \le \frac{m-2-l}{m}
       \]
       by Lemma \ref{lemma: weight of D2(p1,p2)} and Claim \ref{claim mu(m-2-l/m)},
       the line $L$ cannot have this slope. If $\frac{m-2-l}{m}$ is greater than or equal to $\frac{2}{3}$, then we study the case in Figure \ref{Figure finales p1} and Figure \ref{Figure finales p2}. From these it will result that the manifold $Y_1 \cup_{\Sigma}Y_2$ is not $SU(2)$-abelian.
        Studying all possible combinations we conclude that either $Y_1 \cup_{\Sigma} Y_2$ is not $SU(2)$-abelian, or the slope of the line $L$ is not in $S$, and Claim \ref{claim: se non in S allora non SU(2)a} gives the conclusion.
        We separate the proof in two cases: case $\geo{\lambda_1}{\lambda_2}=3$ and
        case $\geo{\lambda_1}{\lambda_2} \in \{4,5\}$.  
        
        \underline{Case $\geo{\lambda_1}{\lambda_2}=3$}.
        In this case $\lambda_2= \nicefrac{3}{n}$, where $n \in \{-1,-2\}$.
        Table \ref{table: intersections pt 1} and Table \ref{table: intersections pt 2} imply that $\frac{m-2-l}{m} < \frac{2}{3}$ in all cases but the one for which $\lambda_2=\nicefrac{3}{2}$, $\nicefrac{p}{q}=-\nicefrac{1}{2}$, and $c = \pi$.

        Let us study this case which is represented by Figure \ref{Figure d=3}.
        Since we assumed that $H(Y_1)\cap A(Y_2) = \emptyset$, in Figure \ref{Figure d=3} the LIPA-path $\gamma_1$ is placed in a way that $\gamma_1 \cup \gamma_1^t$ is disjoint from $A(Y_2)$ and in a minimal position in the following sense.
        By Claim \ref{claim: togliere segments}, if a segment $I_j$ contains one of the points in $\{0,\pi\}^2$, then it is not contained in $\overline{H(Y_2) \cap L}$.
        Similarly, if the LIPA-path $\gamma_1$ passes though a segment $I_{j'}$, then either it is not contained in $H(Y_2) \cap L$ or $Y=Y_1 \cup_{\Sigma} Y_2$ is not $SU(2)$-abelian according to Theorem \ref{teo: M SU(2)-abeliano se e solo se i pezzi sono empty}.
        Therefore, we draw the LIPA-path $\gamma_1$ in a way that it passes through the subsegment $I_j$ so we do not have to exclude an additional subsegment, and therefore to make $H(Y_2) \cap L$ the largest possible.
        
        A direct computation in the considered case shows that $m=|A(Y_2)\cap L|=8$.
        Figure \ref{Figure d=3} shows that there are at most two segments $I_1,I_2$ of $L$ each connecting two points of $A(Y_2)$, that do not intersect $\gamma_2 \cup \gamma_2^t$ and whose accumulation points are not in $\{0,\pi\}^2$.
        Furthermore,
        \[
        \frac{h(I_1)}{h(L)}=
        \frac{h(I_2)}{h(L)}=\frac{1}{m}=\frac{1}{8}.
        \]
        This implies that
        \[
            \frac{h\left(I_1 \cup I_2\right)}{h(L)} = \frac{h(I_1)}{h(L)} + \frac{h(I_2)}{h(L)}=\frac{1}{4}< \frac{2}{3}.
        \]
        Since $\mu(Y_2)=\mu(L) \ge \nicefrac{2}{3}$ by Lemma \ref{lemma: weight of D2(p1,p2)}, if $\lambda_2=\nicefrac{3}{2}$, $\nicefrac{p}{q}=-\nicefrac{1}{2}$, and $c = \pi$,
        then $\overline{L \cap H(Y_2)}$ properly contains $I_1 \cup I_2$ and at least one segment $I_j$ in $\overline{L \setminus (I_1 \cup I_2})$. Since the accumulation points of $L \cap H(Y_2)$ are disjoint from $\{0,\pi\}^2$, the segment $I_j$ is disjoint from $\{0,\pi\}^2$ as well. Figure \ref{Figure d=3} shows that $I_j \cap (\gamma_1 \cup \gamma_1^t) \neq \emptyset$. Thus, $H(Y_1) \cap H(Y_2) \neq \emptyset$. Therefore $Y=Y_1 \cup_{\Sigma} Y_2$ is not $SU(2)$-abelian by Theorem \ref{teo: M SU(2)-abeliano se e solo se i pezzi sono empty}.
        
        Thus, either $Y=Y_1 \cup_{\Sigma} Y_2$ is not $SU(2)$-abelian or  the slope of the line $L$ is not in $S$. Claim \ref{claim: se non in S allora non SU(2)a} gives the conclusion.

        \underline{Case $\geo{\lambda_1}{\lambda_2} \in \{4,5\}$}. These cases follow from the technique showed in the previous cases. Table \ref{table: intersections pt 1} and Table \ref{table: intersections pt 2} show the cases for which $\frac{m-2-l}{m} \ge \frac{2}{3}$. These cases are studied in Figure \ref{Figure finales p1} and Figure \ref{Figure finales p2}.
        All the figures are drawn such that $\gamma_1\cup\gamma_1^t$ is in a minimal position as explained above.
        By studying all the cases reported in Figure \ref{Figure finales p1} and Figure \ref{Figure finales p2} we conclude that either the slope of $L$ is not in $S$, and therefore $Y=Y_1 \cup_{\Sigma} Y_2$ is not $SU(2)$-abelian by Claim \ref{claim: se non in S allora non SU(2)a}, or $L\cap H(Y_2)$ intersects $A(Y_1)\cup H(Y_2)$ which implies that $Y=Y_1\cup_{\Sigma} Y_2$ is not $SU(2)$-abelian by Theorem \ref{teo: M SU(2)-abeliano se e solo se i pezzi sono empty}.
    \end{proof}
\end{lemma}
 
\begin{proof}[Proof of Claim \ref{claim: se non in S allora non SU(2)a}]
    If $\nicefrac{p}{q} \notin S$, then either $|p| \ge 3$ or $|p| \in \{1,2\}$ and $|q| \ge 3$. We prove first that if $|p| \ge 3$, then we get the conclusion. Therefore, we prove that if $|q| \ge 3$ and $|p| \in \{1,2\}$, then we get the conclusion.

                Let $L'$ be a rational line of $\R{\Sigma}$ of slope $\nicefrac{p'}{q'}$. We start by showing that if $|pq'-p'q| \ge 3$, then $L'\cap H(Y_2) \neq \emptyset$.
            
                Let $\mu \subset \Sigma$ be the meridian of the chosen fibration of $Y_2=\mathbb{D}^2(\nicefrac{p_1}{q_1},\nicefrac{p_2}{q_2})$ as in \cite[Definition 2.1.1]{Mino}. Therefore, $\{\mu,h_2\}$ is a basis for $\fund{\Sigma}=\fund{\partial Y_2}$ and if $\fund{Y_2}$ is presented as
                \[
                    \fund{Y_2}=\fund{\mathbb{D}^2(\nicefrac{p_1}{q_1},\nicefrac{p_2}{q_2})}= \left\langle x_1,x_2,h_2\middle| x_1^{p_1}h_2^{q_1}, x_2^{p_2}h_2^{q_2},[x_1,h_2],[x_2,h_2] \right\rangle,
                \]
                then $\mu=x_1x_2 \in \fund{Y_2}$.
                Let us assume that $\R{\Sigma}$ has coordinates $(\theta_2,\psi_2)$ as in \eqref{eq: general coordinates} with respect to the basis $\{\mu,h_2\}$.
                As an application of Fact \ref{Fact: centralizers}, if $\eta \in H(Y_2)$, then $\eta(h_2)= \pm 1$. Since $ \mu(L) \neq 0$, there exists $\varepsilon \in \{0,\pi\}$ such that
                \[
                    L= \left\{\psi_2= \varepsilon \right\} \subset \R{\Sigma}.
                \]
                Therefore, there exist $a_1,a_2 \in [0,2\pi]$, with $a_1 <a_2$ such that
                \begin{equation}
                \left\{ (x,\varepsilon) \in \R{\Sigma} \middle| x \in (a_1,a_2)\right\}\subset L \cap H(Y_2).
                 \end{equation}
                 According to Corollary \ref{cor: jewelled repr},
                 \begin{equation}
                \left\{ (x,\varepsilon) \in \R{\Sigma} \middle| x \in (a_1,a_2) \cup (2\pi - a_2,2\pi a_1)\right\}\subset L \cap H(Y_2).
                \label{eq: una volta Claim Lcap s ge 3}
                 \end{equation}
                 Since $\mu(L)=\mu(Y_2) \ge \nicefrac{2}{3}$ by Lemma \ref{lemma: weight of D2(p1,p2)}, we can assume that $|a_2-a_1| \ge \nicefrac{2\pi}{3}$.
                Since $g_2=1$, \cite[Lemma 7.1]{Mino} implies that $P(Y_2)=\emptyset$. By Lemma
                \ref{lemma: eta in H(Y)-H(Y) in SFS with 2 sf},
                $\overline{H(Y_2)}\setminus H(Y_2)$ cannot contain any point of $\{0,\pi\}^2$. Therefore, we can assume that $a_1,a_2 \in (0,\pi)$.
                
                Let $n = |pq'-p'q|$. Therefore, the line $L$ intersects the line $L'$ in $n$ equidistant points. If $n \ge 3$, the \eqref{eq: una volta Claim Lcap s ge 3} implies that at least one of this intersection point lies in $L \cap H(Y_2)$.

                We recall that $A(Y_1) \subset \R{\Sigma}$ is the line $\{\psi=0\}\subset \R{\Sigma}$ that has slope $\nicefrac{0}{1}$.
                If $|p| \ge 3$, then $A(Y_1) \cap H(Y_2) \neq \emptyset$. The conclusion is given by Theorem \ref{teo: M SU(2)-abeliano se e solo se i pezzi sono empty}.

            We can assume that $|p| \in \{1,2\}$.
            Let us assume that $|q| \ge 3$, we shall prove that $Y=Y_1\cup_{\Sigma}Y_2$ is not $SU(2)$-abelian.
            Without loss of generality, we assume that $p \in \{1,2\}$, therefore, either
            \[
                \frac{p}{q} = \frac{1}{q} \quad \text{or} \quad \frac{p}{q}=\frac{2}{q},
            \]
            with $|q| \ge 3$.

            We recall that $\gamma_1 \colon (-1,1)\to \R{\Sigma}$ is a LIPA-path for $Y_1$.
            The closure $\overline{\gamma_1}$ is homotopic in $\R{\Sigma}$ to the circle $\{\theta=0\}$, see Figure \ref{Figure: sigma small at most 11}.
            Since $\nicefrac{p}{q}$ is either $\nicefrac{1}{q}$ or $\nicefrac{2}{q}$ by assumption,
            we obtain that
                \[
                    \left|L \cap \overline{\gamma_1} \right| = |q|.
                \]
                In particular, the points in $L \cap \overline{\gamma}$ are equally spaced in $L$.
                As we saw before, this implies the one is in $H(Y_2)$ and therefore in $\gamma \cap H(Y_2)$. Hence, if $|q|\ge 3$ then $\gamma \cap H(Y_2) \neq \emptyset$. Therefore $H(Y_1) \cap H(Y_2) \neq \emptyset$. The conclusion holds by Theorem \ref{teo: M SU(2)-abeliano se e solo se i pezzi sono empty}.
            \end{proof}
    \begin{proof}[Proof of Claim \ref{claim mu(m-2-l/m)}]
            By Theorem \ref{teo: admissible pieces have the property P}, $\mu(Y_2)=\mu(L) > 0$.
            As we mentioned in the proof of Lemma \ref{lemma: o2g2o1=1 and sigma <= 5}, the line $L$ contains two points $P_1,P_2 \in \{0,\pi\}^2 \subset\R{\Sigma}$.
            As we said before, the points $\overline{L \cap H(Y_2)} \setminus L \cap H(Y_2)$ are in $A(Y_2)$.
            We recall that $A(Y_2)$ divides $L$ in $m$ subintervals $I_1, \cdots, I_m$ each of the length.
            As an application of Lemma \ref{lemma: eta in H(Y)-H(Y) in SFS with 2 sf}, the accumulation points of $H(Y_2)$ are disjoint from $\{P_1,P_2\}$. Therefore, $L \cap H(Y_2)$ is contained in $\widetilde{L} \subset L$, where $\widetilde{L}$ is the union of the subintervals $I_j$ whose accumulation points are disjoint from $\{P_1,P_2\}$.
            There are at most $m-2$ such subintervals. 
            Notice that, if one (resp. two) point $\{P_1,P_2\}$ is in the intersection
            $L \cap A(Y_2)$, then there are exactly two (resp. four) subinterval $I_j$ with one endpoint in $\{P_1,P_2\}$.
            Thus, when it is not empty, $\widetilde{L}$ is made of at most $m-2-l>0$ subinterval, each of length $h(I_j)=\nicefrac{1}{m}$.
            This implies that
            \[
            \mu(L)= \mu\left( L \cap H(Y_2)\right) \le \frac{h\left( \widetilde{L}\right)}{h(L)} = (m-2-l) \frac{h(I_j)}{h(L)} = \frac{m-2-l}{m}.
            \]
    \end{proof}

     \begin{proof}[Proof of Claim \ref{claim: togliere segments}]
     According to \cite[TBA]{Mino}, since $g_2 =1$, then $P(Y_2) = \emptyset$.
     By Corollary \ref{cor: if eta in H(Y)-H(Y) then P(Y)}, the accumulation points of $H(Y_2)$ are disjoint from $\{0,\pi\}^2$. Therefore, if the segment $I_j$ contains one of the point in $\{0,\pi\}^2$, then it cannot be contained in $\overline{L \cap H(Y_2)}$.
     \end{proof}
\begin{lemma}\label{lemma: o2g2=1 o2=2 and sigma = 2}
    Let $Y_1$ be a graph manifold rational homology solid torus and $Y_2$ a Seifert fibered space with disk base space and two singular fibers.
    Let us further assume that $Y_1 \neq S^1 \times \mathbb{D}^2$, and let $Y=Y_1 \cup_{\Sigma} Y_2$.
    If $\Delta(\lambda_1,\lambda_2) =2$, $o_2=o_1=1$, and either $g_2=1$ or $g_2 =2$, then $Y$ is not $SU(2)$-abelian.
    \begin{proof}
        The proof uses a similar notation to Lemma \ref{lemma: o2g2o1=1 and sigma <= 5}.
        Let $\gamma \colon (-1,1)\to \R{\Sigma}$ be a path.
        With an abuse of notation, we identify $\gamma$ with its image in $\R{\Sigma}$.
            
        Let $\xi \subset \Sigma$ be a slope such that $\{\xi,\lambda_1\}$ is an ordered basis for $\fund{\Sigma}$.
        We give the space $\R{\Sigma}$ coordinates $(\theta,\psi)$ as in \eqref{eq: general coordinates} with respect to this ordered basis.
        According to Remark \ref{rmk: change of basis}, we can assume that the slope $\xi$ is such that
        \[
            \lambda_2=2 \xi -\lambda_1 \subset \Sigma.
        \]
        Therefore by Corollary \ref{cor: A(Y) with rational longitude},
        \[
            A(Y_2) = \left\{ 2 \theta- \psi =0\right\} \subset \R{\Sigma} \quad \text{and} \quad A(Y_1) = \left\{ \psi=0 \right\} \subset \R{\Sigma}.
        \]
        
        Since neither $Y_1$ nor $Y_2$ is a solid torus, then they are both LIPA by Theorem \ref{teo: admissible pieces have the property P}.
        For $i \in \{1,2\}$, we denote by $\gamma_i \colon (-1,1) \to\R{\Sigma}$ a LIPA-path for $Y_i$. We denote $\gamma_1^t$ (resp. $\gamma_2^t$) the path $\gamma_i$ translated by $(\pi,0)$ (resp. $(0,\pi)$). For $i \in \{1,2\}$, the path $\gamma_i^t$ is a LIPA-path of $Y_i$ by Lemma \ref{lemma: translation}.
        If either $\gamma_1$ or $\gamma_1^t$ intersects $A(Y_2)$, then $H(Y_1) \cap A(Y_2) \neq \emptyset$. As an application of Theorem \ref{teo: M SU(2)-abeliano se e solo se i pezzi sono empty}, the manifold $Y$ is not $SU(2)$-abelian.
        Therefore, we can assume that neither $\gamma_1$ not $\gamma_1^t$ intersects $A(Y_2)$.
        Figure \ref{Figure. Case d=2, c=pi} shows $\R{\Sigma}$, $A(Y_1)$, $A(Y_2)$ $\gamma_1$, and $\gamma_2^t$ in the chosen coordinates.

        Let $L \subset \R{\Sigma}$ be a line that contains $\gamma_2$, we recall that $L \not \subset A(Y_2)$. According to \cite[Lemma 7.1]{Mino}, $P(Y_2)=\emptyset$ and therefore the accumulation points of $H(Y_2)$ are disjoint from $\{0,\pi\}^2 \subset \R{\Sigma}$ by Lemma \ref{lemma: eta in H(Y)-H(Y) in SFS with 2 sf}. This implies that $H(Y_2) \cap L$ contains a open segment $I$ such that
        \[
            I \not \subset A(Y_2), \quad \partial I \in A(Y_2), \quad \text{and} \quad \overline{I}\cap \{0,\pi\}^2 = \emptyset.
        \]
        Figure \ref{Figure. Case d=2, c=pi} shows that every such a segment intersect either $A(Y_1)$ or $\gamma_1 \cup \gamma_1^t \subseteq H(Y_1)$. This implies that
        \[
            \left(H(Y_2) \cap L\right) \cap \left(A(Y_1) \cup H(Y_1)\right) \neq \emptyset \quad \text{and therefore} \quad H(Y_2) \cap \left(A(Y_1) \cup H(Y_1)\right) \neq \emptyset.
        \]
        The conclusion holds by Theorem \ref{teo: M SU(2)-abeliano se e solo se i pezzi sono empty}.
    \end{proof}
    \end{lemma}
\begin{table}
    \centering
    \begin{tabular}{c||c|c|c|c|c|c|c|c}
        $\lambda_2$& $-\nicefrac{2}{1}$& $\nicefrac{2}{1}$ & $-\nicefrac{1}{1}$ & $\nicefrac{1}{1}$ & $-\nicefrac{1}{2}$ & $\nicefrac{1}{2}$ & $\nicefrac{0}{1}$ & $\nicefrac{1}{0}$\\
        \hline  
        $-\nicefrac{3}{1}$ & $0$  & $\nicefrac{2}{5}$ & $0$ & $0$ & $\nicefrac{1}{2}$ & $\nicefrac{4}{7}$ & $0$ & $0$ \\
        $-\nicefrac{3}{2}$& $0$  & $\nicefrac{4}{7}$ & $0$ & $\nicefrac{2}{5}$ & $0$ & $\nicefrac{1}{2}$ & $0$ & $0$ \\
        $-\nicefrac{4}{1}$& $0$  & $\nicefrac{1}{3}$ & $0$ & $\nicefrac{2}{5}$ & $\nicefrac{4}{7}$ & \fcolorbox{red}{white}{$\nicefrac{2}{3}$} & $0$ & $0$ \\
        $-\nicefrac{4}{3}$& $0$  & $\nicefrac{3}{5}$ & $0$ & $\nicefrac{4}{7}$ & $\nicefrac{2}{5}$ & \fcolorbox{red}{white}{$\nicefrac{8}{11}$} & $0$ & $0$ \\
        $-\nicefrac{5}{1}$& $0$  & $\nicefrac{4}{7}$ & $0$ & $\nicefrac{1}{3}$ & \fcolorbox{red}{white}{$\nicefrac{2}{3}$} & \fcolorbox{red}{white}{$\nicefrac{8}{11}$} & $\nicefrac{2}{5}$ & $0$ \\
        $-\nicefrac{5}{2}$& $0$  & \fcolorbox{red}{white}{$\nicefrac{2}{3}$} & $0$ & $\nicefrac{4}{7}$ & $\nicefrac{1}{2}$ & \fcolorbox{red}{white}{$\nicefrac{2}{3}$} & $\nicefrac{2}{5}$ & $0$ \\
        $-\nicefrac{5}{3}$& $0$  & \fcolorbox{red}{white}{$\nicefrac{8}{11}$} & $0$ & $\nicefrac{1}{2}$ & $\nicefrac{4}{7}$ & \fcolorbox{red}{white}{$\nicefrac{10}{13}$} & $\nicefrac{2}{5}$ & $0$ \\
        $-\nicefrac{5}{4}$& $0$  & \fcolorbox{red}{white}{$\nicefrac{10}{13}$} & $0$ & \fcolorbox{red}{white}{$\nicefrac{2}{3}$} & $\nicefrac{1}{3}$ & \fcolorbox{red}{white}{$\nicefrac{5}{7}$} & $\nicefrac{2}{5}$ & $0$
    \end{tabular}
    \caption{We compute $\nicefrac{m-2-l}{m}$ in the different cases of $\nicefrac{p}{q}\in S$ and $c = 0$. In red are the cases for which $\nicefrac{m-2-l}{m} \ge \nicefrac{2}{3}$.}
    \label{table: intersections pt 1}
\end{table}
\begin{table}
    \centering
    \begin{tabular}{c||c|c|c|c|c|c|c|c}
        $\lambda_2$& $-\nicefrac{2}{1}$& $\nicefrac{2}{1}$ & $-\nicefrac{1}{1}$ & $\nicefrac{1}{1}$ & $-\nicefrac{1}{2}$ & $\nicefrac{1}{2}$ & $\nicefrac{0}{1}$ & $\nicefrac{1}{0}$\\
        \hline  
        $-\nicefrac{3}{1}$& $0$  & $\nicefrac{2}{5}$ & $0$ & $\nicefrac{1}{2}$ & $\nicefrac{2}{5}$ & $\nicefrac{4}{7}$ & $0$ & $0$ \\
        $-\nicefrac{3}{2}$& $0$  & $\nicefrac{4}{7}$ & $0$ & $\nicefrac{2}{5}$ & $\nicefrac{1}{2}$ & \fcolorbox{red}{white}{$\nicefrac{3}{4}$} & $0$ & $0$ \\
        $-\nicefrac{4}{1}$& $0$  & \fcolorbox{red}{white}{$\nicefrac{2}{3}$} & $0$ & $\nicefrac{2}{5}$ & $\nicefrac{4}{7}$ & \fcolorbox{red}{white}{$\nicefrac{2}{3}$} & $\nicefrac{1}{2}$ & $0$ \\
        $-\nicefrac{4}{3}$& $0$  & \fcolorbox{red}{white}{$\nicefrac{4}{5}$} & $0$ & $\nicefrac{4}{7}$ & $\nicefrac{2}{5}$ & \fcolorbox{red}{white}{$\nicefrac{8}{11}$} & $\nicefrac{1}{2}$ & $0$ \\
        $-\nicefrac{5}{1}$& $0$  & $\nicefrac{4}{7}$ & $\nicefrac{1}{2}$ & \fcolorbox{red}{white}{$\nicefrac{2}{3}$} & \fcolorbox{red}{white}{$\nicefrac{2}{3}$} & \fcolorbox{red}{white}{$\nicefrac{8}{11}$} & $\nicefrac{2}{5}$ & $0$ \\
        $-\nicefrac{5}{2}$& $0$  & \fcolorbox{red}{white}{$\nicefrac{2}{3}$} & $0$ & $\nicefrac{4}{7}$ & \fcolorbox{red}{white}{$\nicefrac{3}{4}$} & \fcolorbox{red}{white}{$\nicefrac{5}{6}$} & $\nicefrac{2}{5}$ & $0$ \\
        $-\nicefrac{5}{3}$& $0$  & \fcolorbox{red}{white}{$\nicefrac{8}{11}$} & $0$ & \fcolorbox{red}{white}{$\nicefrac{3}{4}$} & $\nicefrac{4}{7}$ & \fcolorbox{red}{white}{$\nicefrac{10}{13}$} & $\nicefrac{2}{5}$ & $0$ \\
        $-\nicefrac{5}{4}$& $0$  & \fcolorbox{red}{white}{$\nicefrac{10}{13}$} & $0$ & \fcolorbox{red}{white}{$\nicefrac{2}{3}$} & \fcolorbox{red}{white}{$\nicefrac{2}{3}$} & \fcolorbox{red}{white}{$\nicefrac{6}{7}$} & $\nicefrac{2}{5}$ & $\nicefrac{1}{2}$
    \end{tabular}
    \caption{We compute $\nicefrac{m-2-l}{m}$ in the different cases of $\nicefrac{p}{q}\in S$ and $c = \pi$. In red are the cases for which $\nicefrac{m-2-l}{m} \ge \nicefrac{2}{3}$.}
    \label{table: intersections pt 2}
\end{table}

\begin{rmk}\label{rmk: the pieces}
Let $Y$ is an $SU(2)$-abelian graph manifold rational homology $3$-sphere with $|H_1(Y;\mathbb{Z})| \le 5$, let $\tau \subset Y$ be a system as in Lemma \ref{lemma: the last dec}. If $Y_i$ is a Seifert fibered manifold with torus boundary in $\overline{Y \setminus \tau}$, then it has disk base space by Lemma \ref{lemma: the last dec} and, as a consequence
of Corollary \ref{cor: Y SU(2) allora i due pezzi sono SU(2)-abelian}, it admits an $SU(2)$-abelian Dehn filling.
As a consequence of \cite[Theorem 1.2]{Menangerie} and Remark \ref{rmk: order of the additional fiber of a dehn filling}, either $Y_i=S^1\times \mathbb{D}^2$ or $Y_i$ is one of the following manifold
\begin{equation}\label{eq: le varie SFS su un disco}
    \mathbb{D}^2\left(\frac{p_1}{q_1},\frac{p_2}{q_2}\right), \quad \mathbb{D}^2\left(\frac{2}{1},\frac{4}{q_2},\frac{4}{q_3}\right), \quad \mathbb{D}^2\left(\frac{3}{q_1},\frac{3}{q_2},\frac{3}{q_3}\right), \quad \mathbb{D}^2\left(\frac{2}{1}, \cdots, \frac{2}{1},\frac{p_n}{q_n}\right),
\end{equation}
where $n \ge 3$ and $\gcd(p_i,q_i)=1$.
\end{rmk}

\begin{rmk}\label{rmk: orders and torsions}
    As we mentioned before, the torsion subgroup of $H_1(\mathbb{D}^2(\nicefrac{p_1}{q_1},\nicefrac{p_2}{q_2});\mathbb{Z})$ is equal to $\gcd(p_1,p_2)$. 
    Remark \ref{lemma orders and torsion of C3} implies that the three rightmost manifolds in \eqref{eq: le varie SFS su un disco} have first homology with torsion subgroup of order bigger then or equal to two.
\end{rmk}

\begin{repteo}{teo: sigma small then not SU(2)-abelian}
    Let $Y$ be a toroidal graph manifold rational homology $3$-sphere. If $|H_1(Y;\mathbb{Z})| \le 5$, then $Y$ is not $SU(2)$-abelian.
    \begin{proof}
        Since $Y$ is toroidal, it is neither a lens space nor $\mathbb{RP}^3 \# \mathbb{RP}^3$.
        Therefore, by \cite[Theorem 1.2]{Menangerie}, $Y$ admits a Seifert fibration over $S^2(p_1,p_2,p_3)$ with $(p_1,p_2,p_3)\in \{(2,4,4),(3,3,3)\}$.
        Thus,, there exists a triple $(q_1,q_2,q_3) \in \mathbb{Z}^3$ with $\gcd(p_i,q_i)=1$ such that
        \begin{align}
            \label{eq: una volta fund}
            \fund{Y} = \langle x_1,x_2,x_3,h \,|\, [x_i,h], x_i^{p_i}h^{q_i}, x_1x_2x_3\rangle.
         \end{align}
        A direct computation on \eqref{eq: una volta fund} shows that  
        $|H_1(Y;\mathbb{Z})| \ge 9$. Therefore, if $Y$ is a Seifert fibered manifold with $|H_1(Y;\mathbb{Z})| \le 5$, then $Y$ is not $SU(2)$-abelian.
        
        Let us assume that $Y$ is not a Seifert fibered space.
        By Lemma \ref{lemma: the last dec}, there exists a system of disjoint tori $\tau \subset Y$ such that
        \[
        \overline{Y \setminus \tau} = \left( \bigcup C_3\right) \cup \left(\bigcup M_i\right),
        \]
        where $M_i$ is a Seifert fibered manifold with torus boundary, disk base space and $n \ge 2$ singular fibers.
        Let $\Sigma \in \tau$ be a torus such that
        \[
            \overline{Y \setminus \Sigma} = Y_1 \cup Y_2,
        \]
        where $Y_2$ is a Seifert fibered manifold with disk base space and $Y_1$ is a graph manifold. As shown in Lemma \ref{lemma: system for S2}, up to changing $\Sigma$ with a vertical torus of $Y_2$, we can assume this latter to have exactly two singular fibers.
        
        Let $\sigma \in\mathbb{N}$ be the order of $H_1(Y;\mathbb{Z})$. Since $Y$ is a rational homology $3$-sphere, then $\sigma \neq 0$.
        According to Remark \ref{rmk: order of toroidal},
        \begin{align}
            \label{eq: la uso una volta sigma}
            \sigma \coloneq \left|H_1(Y;\mathbb{Z})\right|=o_1o_2t_1t_2 \geo{\lambda_1}{\lambda_2},
        \end{align}
        where $o_i$ is the order of $\lambda_i$ in $H_1(Y_i,\mathbb{Z})$ and $t_i$ is the order of the torsion subgroup of $H_1(Y_i,\mathbb{Z})$. Clearly, $o_i \ge 1$ and it divides $t_i$.
        We recall that if $Y_2=\mathbb{D}^2(\nicefrac{p_1}{q_1},\nicefrac{p_2}{q_2})$, then $t_2=\gcd(p_1,p_2)=g_2$.
        
        If $\sigma=1$, then the conclusion holds by Proposition \ref{prop: sigma =1 ish}.
        If $\sigma \in \{2,3,5\}$, then the identity in \eqref{eq: la uso una volta sigma} implies that one of the following holds:
        \begin{enumerate}[label=\roman*)]
            \item $o_1o_2t_1=1$, $t_2= \sigma$, and $\geo{\lambda_1}{\lambda_2}=1$;
            \item $o_1o_2t_2=1$, $t_1= \sigma$, and $\geo{\lambda_1}{\lambda_2}=1$;
            \item $o_1o_2t_1t_2=1$ and $\geo{\lambda_1}{\lambda_2}=\sigma$.
        \end{enumerate}
        If either case (\romannum{1}) or case (\romannum{2}) holds, then we get the conclusion by Proposition \ref{prop: sigma =1 ish}. If case (\romannum{3}) holds, then Lemma \ref{lemma: o2g2o1=1 and sigma <= 5} and Lemma \ref{lemma: o2g2=1 o2=2 and sigma = 2} imply the conclusion.

        Let us assume that $\sigma=4$.
        The relation \eqref{eq: la uso una volta sigma} implies that one of the following holds:
        \begin{enumerate}
            \item $o_1o_2t_1t_2=1$ and $\geo{\lambda_1}{\lambda_2}=4$;
            \item $o_1o_2=1$, $t_1=2$, $t_2=1$, and $\geo{\lambda_1}{\lambda_2}=2$;
            \item $o_1o_2=1$, $t_1=1$, $t_2=2$, and $\geo{\lambda_1}{\lambda_2}=2$;
            \item $o_1o_2=1$, $t_1t_2=4$, and $\geo{\lambda_1}{\lambda_2}=1$;
            \item $o_1=t_1=2$, $o_2t_2=1$, and $\geo{\lambda_1}{\lambda_2}=1$;
            \item $o_1t_1=1$, $o_2=t_2=2$, and $\geo{\lambda_1}{\lambda_2}=1$;
        \end{enumerate}
        Lemma \ref{lemma: o2g2o1=1 and sigma <= 5} yields the conclusion if case $(1)$ holds.
        Similarly, Lemma \ref{lemma: o2g2=1 o2=2 and sigma = 2} implies the conclusion if either case $(2)$ or case $(3)$ holds.
        If case $(4)$ holds, then the conclusion is implied by Proposition \ref{prop: sigma =1 ish}. Corollary \ref{cor: sigma=2 o1=delta=1 and o2=2} implies the conclusion under the condition that case $(5)$ holds.
        
        Let us assume that case $(6)$ holds. In particular $t_1=o_1=1$.
        If $Y_1$ is a Seifert fibered manifold, then it is one of the manifolds in \eqref{eq: le varie SFS su un disco}. According to Remark \ref{rmk: the pieces}, since $t_1=1$, the manifold $Y_1$ has two singular fibers.
        In this case the conclusion holds by \cite[Theorem 1.2]{Mino}: if $Y=Y_1 \cup_\Sigma Y_2$ is $SU(2)$-abelian, then it lies in one of the $7$ classes of \cite[Theorem 1.2]{Mino}, it can be proven that neither of them have first homology of order $5$ or less.

        Let us assume that $Y_1$ is not a Seifert fibered manifold.
        By Lemma \ref{lemma: the last dec} there are two tori $\Sigma_1,\Sigma_2 \subset Y_1$ such that
        \[
            Y_1= Y_l \cup_{\Sigma_1} C_3 \cup_{\Sigma_2} Y_r,
        \]
        where $Y_l$ and $Y_r$ are graph manifolds. Since we assumed that $Y_1$ is not a Seifert fibered manifold, then either $Y_l$ or $Y_r$ is not a solid torus. Without loss of generality, we assume that $Y_l$ is not a solid torus.
        Let $t_l$ and $t_r$ the orders of the the torsion subgroups of $H_1(Y_l;\mathbb{Z})$ and $H_1(Y_r;\mathbb{Z})$ respectively.
        As an application of Lemma \ref{lemma orders and torsion of C3}, since $t_1=1$, then $t_l=t_r=1$. Up to iterating this progress, we can assume that $Y_l$ is a Seifert fibered manifold with disk base space and that $t_l=1$.

        Therefore, $Y_l$ is one of the manifold in \eqref{eq: le varie SFS su un disco}. As $t_l=1$, Remark \ref{rmk: orders and torsions} implies that $Y_l$ has two singular fibers. Therefore, we call $Y_4 \coloneq Y_l$, $Y_3 \coloneqq\overline{Y \setminus Y_4}$, and we consider 
        \[
        Y=Y_3 \cup_{\Sigma'} Y_4.
        \]

        As we just shown, $Y_3\neq S^1 \times \mathbb{D}^2$ is a graph manifold with torus boundary and $Y_4$ is a Seifert fibered manifold with disk base space, two singular fibers and $t_l=1$.
        This configuration is not contained in case $(6)$. Therefore this fall into either case $(1)$, case $(2)$, case $(3)$, case $(4)$, or case $(5)$, which are all proven above and this concludes the proof.
    \end{proof}
\end{repteo}

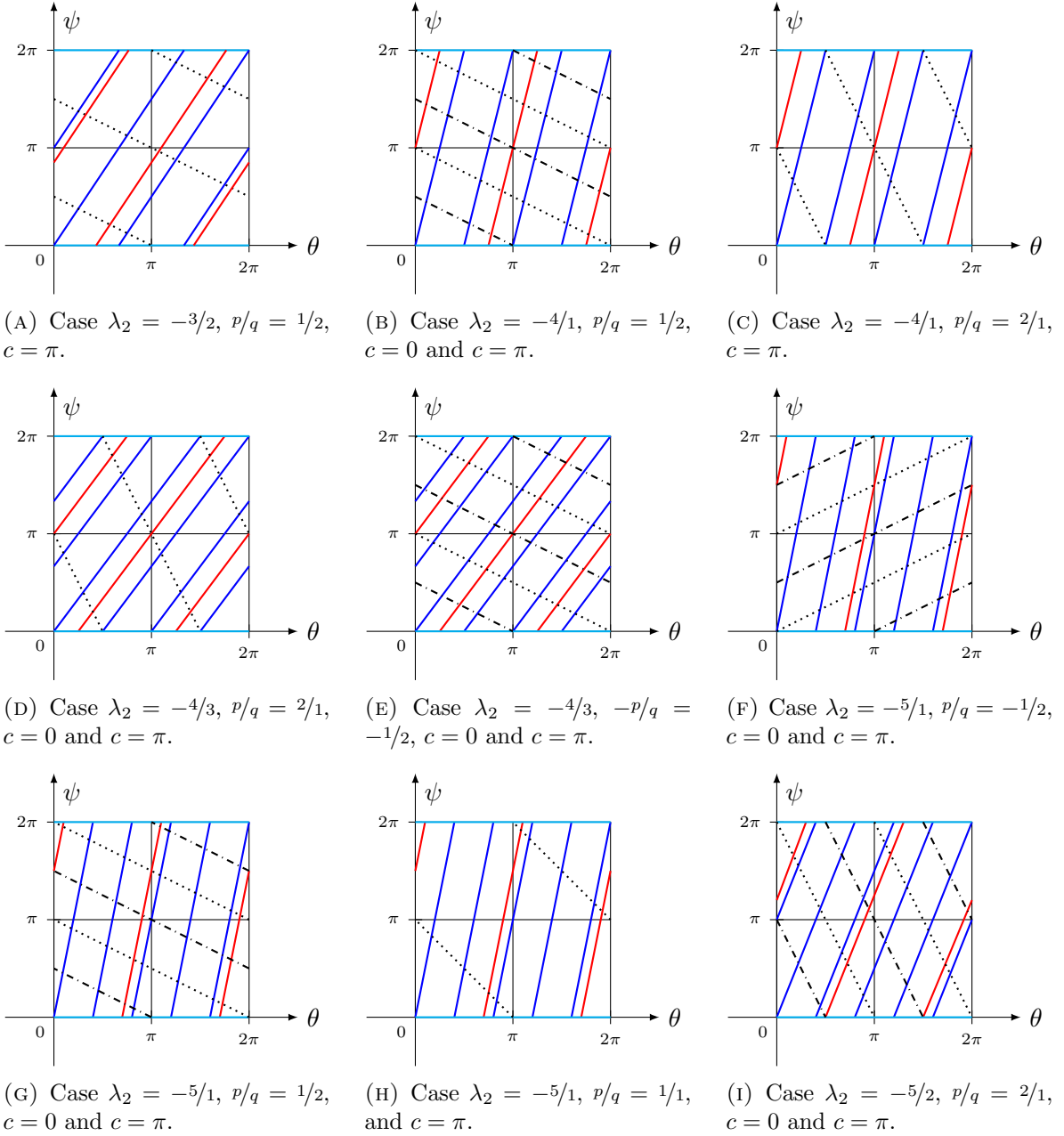
\begin{figure}[t]
    \centering
    \begin{subfigure}[b]{0.3\textwidth}
    \centering
        \begin{tikzpicture}[>=latex,scale=1.45]
            \draw[->] (-.5,0) -- (2.5,0) node[right] {$\theta$};
            \foreach \x /\n in {1/$\pi$,2/$2\pi$} \draw[shift={(\x,0)}] (0pt,2pt) -- (0pt,-2pt) node[below] {\tiny \n};
            \draw[->] (0,-.5) -- (0,2.5) node[below right] {$\psi$};
            \foreach \y /\n in {1/$\pi$,2/$2\pi$}
            \draw[shift={(0,\y)}] (2pt,0pt) -- (-2pt,0pt) node[left] {\tiny \n};
            \node[below left] at (0,0) {\tiny $0$};
             \draw (0,0) rectangle (2,2);
             \draw (1,0)--(1,2);
             \draw (0,1)--(2,1);

             \draw [cyan, thick] (0,0)--(2,0);
              \draw [cyan, thick] (0,2)--(2,2);

            \draw [blue, thick] (0,0)--(4/3,2);
             \draw [blue, thick] (4/3,0)--(2,1);
             \draw [blue, thick] (0,1)--(2/3,2);
             \draw [blue, thick] (2/3,0)--(2,2);
             
            \draw[dotted, thick] (1,2) --(2,1.5);
            \draw[dotted,thick] (0,1.5) --(2,0.5);
            \draw[dotted,thick] (0,0.5) --(1,0);

            \draw[red, thick] (4/3+0.1,0)--(2,1-0.15);
            \draw[red, thick] (0,1-0.15) -- (2/3 +0.1,2);
             \draw[red, thick] (4/3+0.1-1,0)--(2/3 +0.1+1,2);
        \end{tikzpicture}
        \caption{Case $\lambda_2= -\nicefrac{3}{2}$, $\nicefrac{p}{q}=\nicefrac{1}{2}$, $c=\pi$.}
        \label{Figure d=3}
    \end{subfigure}
    \quad
    \begin{subfigure}[b]{0.3\textwidth}
    \centering
        \begin{tikzpicture}[>=latex,scale=1.45]
            \draw[->] (-.5,0) -- (2.5,0) node[right] {$\theta$};
            \foreach \x /\n in {1/$\pi$,2/$2\pi$} \draw[shift={(\x,0)}] (0pt,2pt) -- (0pt,-2pt) node[below] {\tiny \n};
            \draw[->] (0,-.5) -- (0,2.5) node[below right] {$\psi$};
            \foreach \y /\n in {1/$\pi$,2/$2\pi$}
            \draw[shift={(0,\y)}] (2pt,0pt) -- (-2pt,0pt) node[left] {\tiny \n};
            \node[below left] at (0,0) {\tiny $0$};
             \draw (0,0) rectangle (2,2);
             \draw (1,0)--(1,2);
             \draw (0,1)--(2,1);

             \draw[blue, thick] (0,0)--(1/2,2);
             \draw[blue, thick] (1/2,0)--(1,2);
             \draw[blue, thick] (1.5,2)--(1,0);
             \draw[blue, thick] (1.5,0)--(2,2);

             \draw[red, thick] (0.75,0)--(1.25,2);
             \draw[red, thick] (1.75,0)--(2,1);
             \draw[red, thick] (0,1)--(0.25,2);

             \draw [cyan, thick] (0,0)--(2,0);
              \draw [cyan, thick] (0,2)--(2,2);

            \draw[dotted, thick] (0,2) --(2,1);
            \draw[dotted,thick] (0,1) --(2,0);

            \draw[ dashdotted, thick] (1,2) --(2,1.5);
            \draw[ dashdotted, thick] (0,1.5) --(2,0.5);
            \draw[ dashdotted, thick] (0,0.5) --(1,0);

        \end{tikzpicture}
        \caption{Case $\lambda_2= -\nicefrac{4}{1}$, $\nicefrac{p}{q}=\nicefrac{1}{2}$, $c=0$ and $c=\pi$.}
        \label{Figure d=4, -1/2, c=0}
    \end{subfigure}
    \quad
    \begin{subfigure}[b]{0.3\textwidth}
    \centering
         \begin{tikzpicture}[>=latex,scale=1.45]
            \draw[->] (-.5,0) -- (2.5,0) node[right] {$\theta$};
            \foreach \x /\n in {1/$\pi$,2/$2\pi$} \draw[shift={(\x,0)}] (0pt,2pt) -- (0pt,-2pt) node[below] {\tiny \n};
            \draw[->] (0,-.5) -- (0,2.5) node[below right] {$\psi$};
            \foreach \y /\n in {1/$\pi$,2/$2\pi$}
            \draw[shift={(0,\y)}] (2pt,0pt) -- (-2pt,0pt) node[left] {\tiny \n};
            \node[below left] at (0,0) {\tiny $0$};
             \draw (0,0) rectangle (2,2);
             \draw (1,0)--(1,2);
             \draw (0,1)--(2,1);

             \draw[blue, thick] (0,0)--(1/2,2);
             \draw[blue, thick] (1/2,0)--(1,2);
             \draw[blue, thick] (1.5,2)--(1,0);
             \draw[blue, thick] (1.5,0)--(2,2);

             \draw[red, thick] (0.75,0)--(1.25,2);
             \draw[red, thick] (1.75,0)--(2,1);
             \draw[red, thick] (0,1)--(0.25,2);

             \draw [cyan, thick] (0,0)--(2,0);
              \draw [cyan, thick] (0,2)--(2,2);

            \draw[dotted, thick] (0,1) --(0.5,0);
            \draw[dotted, thick] (0.5,2) --(1.5,0);
            \draw[dotted, thick] (1.5,2) --(2,1);

        \end{tikzpicture}
       \caption{Case $\lambda_2= -\nicefrac{4}{1}$, $\nicefrac{p}{q}=\nicefrac{2}{1}$, $c=\pi$.}
        \label{Figure: d=4, -2/1, c=pi}
    \end{subfigure} \\[8pt]

    \centering
    \begin{subfigure}[b]{0.3\textwidth}
    \centering
         \begin{tikzpicture}[>=latex,scale=1.45]
            \draw[->] (-.5,0) -- (2.5,0) node[right] {$\theta$};
            \foreach \x /\n in {1/$\pi$,2/$2\pi$} \draw[shift={(\x,0)}] (0pt,2pt) -- (0pt,-2pt) node[below] {\tiny \n};
            \draw[->] (0,-.5) -- (0,2.5) node[below right] {$\psi$};
            \foreach \y /\n in {1/$\pi$,2/$2\pi$}
            \draw[shift={(0,\y)}] (2pt,0pt) -- (-2pt,0pt) node[left] {\tiny \n};
            \node[below left] at (0,0) {\tiny $0$};
             \draw (0,0) rectangle (2,2);
             \draw (1,0)--(1,2);
             \draw (0,1)--(2,1);

             \draw[blue, thick] (0,0)--(1.5,2);
             \draw[blue, thick] (1.5,0)--(2,2/3);
             \draw[blue, thick] (0,2/3)--(1,2);
             \draw[blue, thick] (1,0)--(2,4/3);
             \draw[blue, thick] (0,4/3) --(0.5,2); 
             \draw[blue, thick] (0.5,0) --(2,2); 

             \draw[red, thick] (0.25,0)--(1.75,2);
             \draw[red, thick] (1.25,0)--(2,1);
             \draw[red, thick] (0,1)--(0.75,2);
             
             \draw [cyan, thick] (0,0)--(2,0);
              \draw [cyan, thick] (0,2)--(2,2);

            \draw[dotted, thick] (0,1) --(0.5,0);
            \draw[dotted, thick] (0.5,2) --(1.5,0);
            \draw[dotted, thick] (1.5,2) --(2,1);
              
        \end{tikzpicture}
        \caption{Case $\lambda_2= -\nicefrac{4}{3}$, $\nicefrac{p}{q}=\nicefrac{2}{1}$, $c=0$ and $c=\pi$.}
         
    \end{subfigure}
    \quad
    \begin{subfigure}[b]{0.3\textwidth}
    \centering
        \begin{tikzpicture}[>=latex,scale=1.45]
            \draw[->] (-.5,0) -- (2.5,0) node[right] {$\theta$};
            \foreach \x /\n in {1/$\pi$,2/$2\pi$} \draw[shift={(\x,0)}] (0pt,2pt) -- (0pt,-2pt) node[below] {\tiny \n};
            \draw[->] (0,-.5) -- (0,2.5) node[below right] {$\psi$};
            \foreach \y /\n in {1/$\pi$,2/$2\pi$}
            \draw[shift={(0,\y)}] (2pt,0pt) -- (-2pt,0pt) node[left] {\tiny \n};
            \node[below left] at (0,0) {\tiny $0$};
             \draw (0,0) rectangle (2,2);
             \draw (1,0)--(1,2);
             \draw (0,1)--(2,1);

             \draw[blue, thick] (0,0)--(1.5,2);
             \draw[blue, thick] (1.5,0)--(2,2/3);
             \draw[blue, thick] (0,2/3)--(1,2);
             \draw[blue, thick] (1,0)--(2,4/3);
             \draw[blue, thick] (0,4/3) --(0.5,2); 
             \draw[blue, thick] (0.5,0) --(2,2); 

             \draw[red, thick] (0.25,0)--(1.75,2);
             \draw[red, thick] (1.25,0)--(2,1);
             \draw[red, thick] (0,1)--(0.75,2);

             \draw [cyan, thick] (0,0)--(2,0);
              \draw [cyan, thick] (0,2)--(2,2);

            \draw[dotted, thick] (0,2) --(2,1);
            \draw[dotted,thick] (0,1) --(2,0);

            \draw[dashdotted, thick] (1,2) --(2,1.5);
            \draw[ dashdotted, thick] (0,1.5) --(2,0.5);
            \draw[ dashdotted, thick] (0,0.5) --(1,0);
              
        \end{tikzpicture}
        \caption{Case $\lambda_2= -\nicefrac{4}{3}$, $-\nicefrac{p}{q}=-\nicefrac{1}{2}$, $c=0$ and $c=\pi$.}
    \end{subfigure}
    \quad
    \begin{subfigure}[b]{0.3\textwidth}
    \centering
         \begin{tikzpicture}[>=latex,scale=1.45]
            \draw[->] (-.5,0) -- (2.5,0) node[right] {$\theta$};
            \foreach \x /\n in {1/$\pi$,2/$2\pi$} \draw[shift={(\x,0)}] (0pt,2pt) -- (0pt,-2pt) node[below] {\tiny \n};
            \draw[->] (0,-.5) -- (0,2.5) node[below right] {$\psi$};
            \foreach \y /\n in {1/$\pi$,2/$2\pi$}
            \draw[shift={(0,\y)}] (2pt,0pt) -- (-2pt,0pt) node[left] {\tiny \n};
            \node[below left] at (0,0) {\tiny $0$};
             \draw (0,0) rectangle (2,2);
             \draw (1,0)--(1,2);
             \draw (0,1)--(2,1);

             \draw[blue, thick] (0,0)--(2/5,2);
             \draw[blue, thick] (2/5,0)--(4/5,2);
             \draw[blue, thick] (4/5,0)--(6/5,2);
             \draw[blue, thick] (6/5,0)--(8/5,2);
             \draw[blue, thick] (8/5,0)--(2,2);

             \draw[red, thick] (8/5+0.1,0)--(2,1.5);
             \draw[red, thick] (0,1.5)--(0.1,2);
            \draw[red, thick] (4/5-0.1,0)--(1.1,2);

             \draw [cyan, thick] (0,0)--(2,0);
              \draw [cyan, thick] (0,2)--(2,2);

            \draw[dotted, thick] (0,0) --(2,1);
            \draw[dotted, thick] (0,1) --(2,2);

            \draw[dashdotted, thick] (1,0)--(2,0.5);
            \draw[dashdotted, thick] (0,0.5)--(2,1.5);
            \draw[dashdotted, thick] (0,1.5)--(1,2);
           
        \end{tikzpicture}
       \caption{Case $\lambda_2= -\nicefrac{5}{1}$, $\nicefrac{p}{q}=-\nicefrac{1}{2}$, $c=0$ and $c=\pi$.}
    \end{subfigure}
    \\[8pt]
    \centering
    \begin{subfigure}[b]{0.3\textwidth}
    \centering
         \begin{tikzpicture}[>=latex,scale=1.45]
            \draw[->] (-.5,0) -- (2.5,0) node[right] {$\theta$};
            \foreach \x /\n in {1/$\pi$,2/$2\pi$} \draw[shift={(\x,0)}] (0pt,2pt) -- (0pt,-2pt) node[below] {\tiny \n};
            \draw[->] (0,-.5) -- (0,2.5) node[below right] {$\psi$};
            \foreach \y /\n in {1/$\pi$,2/$2\pi$}
            \draw[shift={(0,\y)}] (2pt,0pt) -- (-2pt,0pt) node[left] {\tiny \n};
            \node[below left] at (0,0) {\tiny $0$};
             \draw (0,0) rectangle (2,2);
             \draw (1,0)--(1,2);
             \draw (0,1)--(2,1);

             \draw[blue, thick] (0,0)--(2/5,2);
             \draw[blue, thick] (2/5,0)--(4/5,2);
             \draw[blue, thick] (4/5,0)--(6/5,2);
             \draw[blue, thick] (6/5,0)--(8/5,2);
             \draw[blue, thick] (8/5,0)--(2,2);

             \draw[red, thick] (8/5+0.1,0)--(2,1.5);
             \draw[red, thick] (0,1.5)--(0.1,2);
            \draw[red, thick] (4/5-0.1,0)--(1.1,2);

             \draw [cyan, thick] (0,0)--(2,0);
              \draw [cyan, thick] (0,2)--(2,2);

            \draw[dotted, thick] (2,0) --(0,1);
            \draw[dotted, thick] (2,1) --(0,2);

            \draw[dashdotted, thick] (1,0)--(0,0.5);
            \draw[dashdotted, thick] (2,0.5)--(0,1.5);
            \draw[dashdotted, thick] (2,1.5)--(1,2);
               
        \end{tikzpicture}
       \caption{Case $\lambda_2= -\nicefrac{5}{1}$, $\nicefrac{p}{q}=\nicefrac{1}{2}$, $c=0$ and $c=\pi$.}         
    \end{subfigure}
    \quad
    \begin{subfigure}[b]{0.3\textwidth}
    \centering
         \begin{tikzpicture}[>=latex,scale=1.45]
            \draw[->] (-.5,0) -- (2.5,0) node[right] {$\theta$};
            \foreach \x /\n in {1/$\pi$,2/$2\pi$} \draw[shift={(\x,0)}] (0pt,2pt) -- (0pt,-2pt) node[below] {\tiny \n};
            \draw[->] (0,-.5) -- (0,2.5) node[below right] {$\psi$};
            \foreach \y /\n in {1/$\pi$,2/$2\pi$}
            \draw[shift={(0,\y)}] (2pt,0pt) -- (-2pt,0pt) node[left] {\tiny \n};
            \node[below left] at (0,0) {\tiny $0$};
             \draw (0,0) rectangle (2,2);
             \draw (1,0)--(1,2);
             \draw (0,1)--(2,1);

             \draw[blue, thick] (0,0)--(2/5,2);
             \draw[blue, thick] (2/5,0)--(4/5,2);
             \draw[blue, thick] (4/5,0)--(6/5,2);
             \draw[blue, thick] (6/5,0)--(8/5,2);
             \draw[blue, thick] (8/5,0)--(2,2);

             \draw[red, thick] (8/5+0.1,0)--(2,1.5);
             \draw[red, thick] (0,1.5)--(0.1,2);
            \draw[red, thick] (4/5-0.1,0)--(1.1,2);

             \draw [cyan, thick] (0,0)--(2,0);
              \draw [cyan, thick] (0,2)--(2,2);

            \draw[dotted, thick] (1,0) --(0,1);
            \draw[dotted, thick] (2,1) --(1,2);
               
        \end{tikzpicture}
       \caption{Case $\lambda_2= -\nicefrac{5}{1}$, $\nicefrac{p}{q}=\nicefrac{1}{1}$, and $c=\pi$.}
    \end{subfigure}
    \quad
    \begin{subfigure}[b]{0.3\textwidth}
    \centering
         \begin{tikzpicture}[>=latex,scale=1.45]
            \draw[->] (-.5,0) -- (2.5,0) node[right] {$\theta$};
            \foreach \x /\n in {1/$\pi$,2/$2\pi$} \draw[shift={(\x,0)}] (0pt,2pt) -- (0pt,-2pt) node[below] {\tiny \n};
            \draw[->] (0,-.5) -- (0,2.5) node[below right] {$\psi$};
            \foreach \y /\n in {1/$\pi$,2/$2\pi$}
            \draw[shift={(0,\y)}] (2pt,0pt) -- (-2pt,0pt) node[left] {\tiny \n};
            \node[below left] at (0,0) {\tiny $0$};
             \draw (0,0) rectangle (2,2);
             \draw (1,0)--(1,2);
             \draw (0,1)--(2,1);

             \draw[blue, thick] (0,0)--(4/5,2);
             \draw[blue, thick] (4/5,0)--(8/5,2);
             \draw[blue, thick] (8/5,0)--(2,1);
             \draw[blue, thick] (0,1)--(2/5,2);
             \draw[blue, thick] (2/5,0)--(6/5,2);
             \draw[blue, thick] (6/5,0)--(2,2);

             \draw[red, thick] (3/5-0.1,0)--(7/5-0.1,2);
             \draw[red, thick] (8/5-0.1,0)--(2,1.2);
            \draw[red, thick] (0,1.2)--(2/5-0.1,2);

             \draw [cyan, thick] (0,0)--(2,0);
              \draw [cyan, thick] (0,2)--(2,2);

            \draw[dotted, thick] (0,2) --(1,0);
            \draw[dotted, thick] (1,2) --(2,0);

            \draw[dashdotted, thick] (0,1) --(0.5,0);
            \draw[dashdotted, thick] (0.5,2) --(1.5,0);
            \draw[dashdotted, thick] (1.5,2) --(2,1);
        \end{tikzpicture}
       \caption{Case $\lambda_2= -\nicefrac{5}{2}$, $\nicefrac{p}{q}=\nicefrac{2}{1}$, $c=0$ and $c=\pi$.}
    \end{subfigure}
   \caption{$A(Y_1)$ and $A(Y_2)$ are in cyan and blue. The LIPA-path of $Y_1$ is in red and the dotted and dash-dotted lines represent $L$ in the cases $c=0$ and $c=\pi$.}
   \label{Figure finales p1}
    \end{figure}
    
    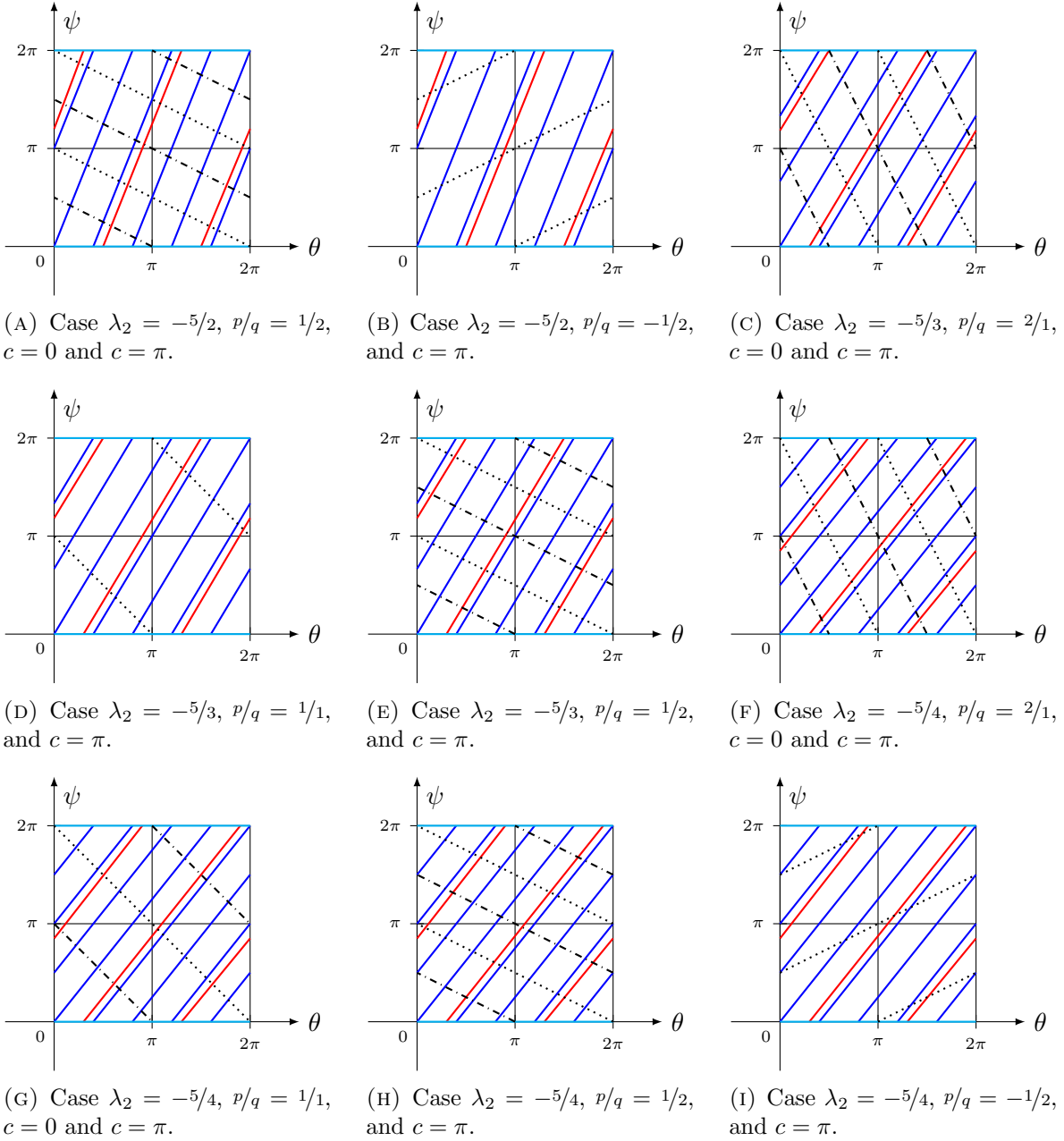
\begin{figure}
    \begin{subfigure}[b]{0.3\textwidth}
    \centering
         \begin{tikzpicture}[>=latex,scale=1.45]
            \draw[->] (-.5,0) -- (2.5,0) node[right] {$\theta$};
            \foreach \x /\n in {1/$\pi$,2/$2\pi$} \draw[shift={(\x,0)}] (0pt,2pt) -- (0pt,-2pt) node[below] {\tiny \n};
            \draw[->] (0,-.5) -- (0,2.5) node[below right] {$\psi$};
            \foreach \y /\n in {1/$\pi$,2/$2\pi$}
            \draw[shift={(0,\y)}] (2pt,0pt) -- (-2pt,0pt) node[left] {\tiny \n};
            \node[below left] at (0,0) {\tiny $0$};
             \draw (0,0) rectangle (2,2);
             \draw (1,0)--(1,2);
             \draw (0,1)--(2,1);

             \draw[blue, thick] (0,0)--(4/5,2);
             \draw[blue, thick] (4/5,0)--(8/5,2);
             \draw[blue, thick] (8/5,0)--(2,1);
             \draw[blue, thick] (0,1)--(2/5,2);
             \draw[blue, thick] (2/5,0)--(6/5,2);
             \draw[blue, thick] (6/5,0)--(2,2);

            \draw[red, thick] (3/5-0.1,0)--(7/5-0.1,2);
             \draw[red, thick] (8/5-0.1,0)--(2,1.2);
            \draw[red, thick] (0,1.2)--(2/5-0.1,2);
             
             \draw [cyan, thick] (0,0)--(2,0);
              \draw [cyan, thick] (0,2)--(2,2);

            \draw[dotted, thick] (2,0) --(0,1);
            \draw[dotted, thick] (2,1) --(0,2);

            \draw[dashdotted, thick] (1,0)--(0,0.5);
            \draw[dashdotted, thick] (2,0.5)--(0,1.5);
            \draw[dashdotted, thick] (2,1.5)--(1,2);
               
        \end{tikzpicture}
       \caption{Case $\lambda_2= -\nicefrac{5}{2}$, $\nicefrac{p}{q}=\nicefrac{1}{2}$, $c=0$ and $c=\pi$.}         
    \end{subfigure}
    \quad
    \begin{subfigure}[b]{0.3\textwidth}
    \centering
         \begin{tikzpicture}[>=latex,scale=1.45]
            \draw[->] (-.5,0) -- (2.5,0) node[right] {$\theta$};
            \foreach \x /\n in {1/$\pi$,2/$2\pi$} \draw[shift={(\x,0)}] (0pt,2pt) -- (0pt,-2pt) node[below] {\tiny \n};
            \draw[->] (0,-.5) -- (0,2.5) node[below right] {$\psi$};
            \foreach \y /\n in {1/$\pi$,2/$2\pi$}
            \draw[shift={(0,\y)}] (2pt,0pt) -- (-2pt,0pt) node[left] {\tiny \n};
            \node[below left] at (0,0) {\tiny $0$};
             \draw (0,0) rectangle (2,2);
             \draw (1,0)--(1,2);
             \draw (0,1)--(2,1);

             \draw[blue, thick] (0,0)--(4/5,2);
             \draw[blue, thick] (4/5,0)--(8/5,2);
             \draw[blue, thick] (8/5,0)--(2,1);
             \draw[blue, thick] (0,1)--(2/5,2);
             \draw[blue, thick] (2/5,0)--(6/5,2);
             \draw[blue, thick] (6/5,0)--(2,2);

            \draw[red, thick] (3/5-0.1,0)--(7/5-0.1,2);
             \draw[red, thick] (8/5-0.1,0)--(2,1.2);
            \draw[red, thick] (0,1.2)--(2/5-0.1,2);

             \draw [cyan, thick] (0,0)--(2,0);
              \draw [cyan, thick] (0,2)--(2,2);

            \draw[dotted, thick] (1,0)--(2,0.5);
            \draw[dotted, thick] (0,0.5)--(2,1.5);
            \draw[dotted, thick] (0,1.5)--(1,2);
               
        \end{tikzpicture}
       \caption{Case $\lambda_2=- \nicefrac{5}{2}$, $\nicefrac{p}{q}=-\nicefrac{1}{2}$, and $c=\pi$.}
    \end{subfigure}
    \quad
    \begin{subfigure}[b]{0.3\textwidth}
    \centering
         \begin{tikzpicture}[>=latex,scale=1.45]
            \draw[->] (-.5,0) -- (2.5,0) node[right] {$\theta$};
            \foreach \x /\n in {1/$\pi$,2/$2\pi$} \draw[shift={(\x,0)}] (0pt,2pt) -- (0pt,-2pt) node[below] {\tiny \n};
            \draw[->] (0,-.5) -- (0,2.5) node[below right] {$\psi$};
            \foreach \y /\n in {1/$\pi$,2/$2\pi$}
            \draw[shift={(0,\y)}] (2pt,0pt) -- (-2pt,0pt) node[left] {\tiny \n};
            \node[below left] at (0,0) {\tiny $0$};
             \draw (0,0) rectangle (2,2);
             \draw (1,0)--(1,2);
             \draw (0,1)--(2,1);

             \draw[blue, thick] (0,0)--(6/5,2);
             \draw[blue, thick] (6/5,0)--(2,4/3);
             \draw[blue, thick] (0,4/3)--(2/5,2);
             \draw[blue, thick] (2/5,0)--(8/5,2);
             \draw[blue, thick] (8/5,0)--(2,2/3);
             \draw[blue, thick] (0,2/3)--(4/5,2);
             \draw[blue, thick] (4/5,0)--(2,2);

             \draw[red, thick] (1/5+0.1,0)--(7/5+0.1,2);
             \draw[red, thick] (6/5+0.1,0)--(2,4/3-0.15);
             \draw[red, thick] (2/5+0.1,2)--(0,4/3-0.15);

             \draw [cyan, thick] (0,0)--(2,0);
              \draw [cyan, thick] (0,2)--(2,2);

            \draw[dotted, thick] (0,2) --(1,0);
            \draw[dotted, thick] (1,2) --(2,0);

            \draw[dashdotted, thick] (0,1) --(0.5,0);
            \draw[dashdotted, thick] (0.5,2) --(1.5,0);
            \draw[dashdotted, thick] (1.5,2) --(2,1);
               
        \end{tikzpicture}
       \caption{Case $\lambda_2= -\nicefrac{5}{3}$, $\nicefrac{p}{q}=\nicefrac{2}{1}$, $c=0$ and $c=\pi$.}
    \end{subfigure} \\[8pt]

    \centering
    \begin{subfigure}[b]{0.3\textwidth}
    \centering
         \begin{tikzpicture}[>=latex,scale=1.45]
            \draw[->] (-.5,0) -- (2.5,0) node[right] {$\theta$};
            \foreach \x /\n in {1/$\pi$,2/$2\pi$} \draw[shift={(\x,0)}] (0pt,2pt) -- (0pt,-2pt) node[below] {\tiny \n};
            \draw[->] (0,-.5) -- (0,2.5) node[below right] {$\psi$};
            \foreach \y /\n in {1/$\pi$,2/$2\pi$}
            \draw[shift={(0,\y)}] (2pt,0pt) -- (-2pt,0pt) node[left] {\tiny \n};
            \node[below left] at (0,0) {\tiny $0$};
             \draw (0,0) rectangle (2,2);
             \draw (1,0)--(1,2);
             \draw (0,1)--(2,1);

             \draw[blue, thick] (0,0)--(6/5,2);
             \draw[blue, thick] (6/5,0)--(2,4/3);
             \draw[blue, thick] (0,4/3)--(2/5,2);
             \draw[blue, thick] (2/5,0)--(8/5,2);
             \draw[blue, thick] (8/5,0)--(2,2/3);
             \draw[blue, thick] (0,2/3)--(4/5,2);
             \draw[blue, thick] (4/5,0)--(2,2);

             \draw[red, thick] (1/5+0.1,0)--(7/5+0.1,2);
             \draw[red, thick] (6/5+0.1,0)--(2,4/3-0.15);
             \draw[red, thick] (2/5+0.1,2)--(0,4/3-0.15);

             \draw [cyan, thick] (0,0)--(2,0);
              \draw [cyan, thick] (0,2)--(2,2);

            \draw[dotted, thick] (1,2) --(2,1);
            \draw[dotted, thick] (0,1) --(1,0);

        \end{tikzpicture}
       \caption{Case $\lambda_2= -\nicefrac{5}{3}$, $\nicefrac{p}{q}=\nicefrac{1}{1}$, and $c=\pi$.}         
    \end{subfigure}
    \quad
    \begin{subfigure}[b]{0.3\textwidth}
    \centering
         \begin{tikzpicture}[>=latex,scale=1.45]
            \draw[->] (-.5,0) -- (2.5,0) node[right] {$\theta$};
            \foreach \x /\n in {1/$\pi$,2/$2\pi$} \draw[shift={(\x,0)}] (0pt,2pt) -- (0pt,-2pt) node[below] {\tiny \n};
            \draw[->] (0,-.5) -- (0,2.5) node[below right] {$\psi$};
            \foreach \y /\n in {1/$\pi$,2/$2\pi$}
            \draw[shift={(0,\y)}] (2pt,0pt) -- (-2pt,0pt) node[left] {\tiny \n};
            \node[below left] at (0,0) {\tiny $0$};
             \draw (0,0) rectangle (2,2);
             \draw (1,0)--(1,2);
             \draw (0,1)--(2,1);

             \draw[blue, thick] (0,0)--(6/5,2);
             \draw[blue, thick] (6/5,0)--(2,4/3);
             \draw[blue, thick] (0,4/3)--(2/5,2);
             \draw[blue, thick] (2/5,0)--(8/5,2);
             \draw[blue, thick] (8/5,0)--(2,2/3);
             \draw[blue, thick] (0,2/3)--(4/5,2);
             \draw[blue, thick] (4/5,0)--(2,2);

             \draw[red, thick] (1/5+0.1,0)--(7/5+0.1,2);
             \draw[red, thick] (6/5+0.1,0)--(2,4/3-0.15);
             \draw[red, thick] (2/5+0.1,2)--(0,4/3-0.15);

             \draw [cyan, thick] (0,0)--(2,0);
              \draw [cyan, thick] (0,2)--(2,2);

            \draw[dotted, thick] (2,0) --(0,1);
            \draw[dotted, thick] (2,1) --(0,2);

            \draw[dashdotted, thick] (1,0)--(0,0.5);
            \draw[dashdotted, thick] (2,0.5)--(0,1.5);
            \draw[dashdotted, thick] (2,1.5)--(1,2);
               
        \end{tikzpicture}
       \caption{Case $\lambda_2= -\nicefrac{5}{3}$, $\nicefrac{p}{q}=\nicefrac{1}{2}$, and $c=\pi$.}
    \end{subfigure}
    \quad
    \begin{subfigure}[b]{0.3\textwidth}
    \centering
         \begin{tikzpicture}[>=latex,scale=1.45]
            \draw[->] (-.5,0) -- (2.5,0) node[right] {$\theta$};
            \foreach \x /\n in {1/$\pi$,2/$2\pi$} \draw[shift={(\x,0)}] (0pt,2pt) -- (0pt,-2pt) node[below] {\tiny \n};
            \draw[->] (0,-.5) -- (0,2.5) node[below right] {$\psi$};
            \foreach \y /\n in {1/$\pi$,2/$2\pi$}
            \draw[shift={(0,\y)}] (2pt,0pt) -- (-2pt,0pt) node[left] {\tiny \n};
            \node[below left] at (0,0) {\tiny $0$};
             \draw (0,0) rectangle (2,2);
             \draw (1,0)--(1,2);
             \draw (0,1)--(2,1);

             \draw[blue, thick] (0,0)--(8/5,2);
             \draw[blue, thick] (8/5,0)--(2,0.5);
             \draw[blue, thick] (0,0.5)--(6/5,2);
             \draw[blue, thick] (6/5,0)--(2,1);
             \draw[blue, thick] (0,1)--(4/5,2);
             \draw[blue, thick] (4/5,0)--(2,1.5);
             \draw[blue, thick] (0,1.5)--(2/5,2);
             \draw[blue, thick] (2/5,0)--(2,2);

             \draw[red, thick] (6/5+0.1,0)--(2,1-0.15);
             \draw[red, thick] (1/5+0.1,0)--(9/5+0.1,2);
            \draw[red, thick] (0,1-0.15)--(4/5+0.1,2);

             \draw [cyan, thick] (0,0)--(2,0);
              \draw [cyan, thick] (0,2)--(2,2);

            \draw[dotted, thick] (0,2) --(1,0);
            \draw[dotted, thick] (1,2) --(2,0);

            \draw[dashdotted, thick] (0,1) --(0.5,0);
            \draw[dashdotted, thick] (0.5,2) --(1.5,0);
            \draw[dashdotted, thick] (1.5,2) --(2,1);
        \end{tikzpicture}
       \caption{Case $\lambda_2= -\nicefrac{5}{4}$, $\nicefrac{p}{q}=\nicefrac{2}{1}$, $c=0$ and $c=\pi$.}
    \end{subfigure}
    \\[8pt]
    \centering
    \begin{subfigure}[b]{0.3\textwidth}
    \centering
         \begin{tikzpicture}[>=latex,scale=1.45]
            \draw[->] (-.5,0) -- (2.5,0) node[right] {$\theta$};
            \foreach \x /\n in {1/$\pi$,2/$2\pi$} \draw[shift={(\x,0)}] (0pt,2pt) -- (0pt,-2pt) node[below] {\tiny \n};
            \draw[->] (0,-.5) -- (0,2.5) node[below right] {$\psi$};
            \foreach \y /\n in {1/$\pi$,2/$2\pi$}
            \draw[shift={(0,\y)}] (2pt,0pt) -- (-2pt,0pt) node[left] {\tiny \n};
            \node[below left] at (0,0) {\tiny $0$};
             \draw (0,0) rectangle (2,2);
             \draw (1,0)--(1,2);
             \draw (0,1)--(2,1);

             \draw[blue, thick] (0,0)--(8/5,2);
             \draw[blue, thick] (8/5,0)--(2,0.5);
             \draw[blue, thick] (0,0.5)--(6/5,2);
             \draw[blue, thick] (6/5,0)--(2,1);
             \draw[blue, thick] (0,1)--(4/5,2);
             \draw[blue, thick] (4/5,0)--(2,1.5);
             \draw[blue, thick] (0,1.5)--(2/5,2);
             \draw[blue, thick] (2/5,0)--(2,2);

             \draw[red, thick] (6/5+0.1,0)--(2,1-0.15);
             \draw[red, thick] (1/5+0.1,0)--(9/5+0.1,2);
            \draw[red, thick] (0,1-0.15)--(4/5+0.1,2);

             \draw [cyan, thick] (0,0)--(2,0);
              \draw [cyan, thick] (0,2)--(2,2);

              \draw[dotted, thick] (0,2) --(2,0);

            \draw[dashdotted, thick] (1,2) --(2,1);
            \draw[dashdotted, thick] (0,1) --(1,0);

        \end{tikzpicture}
       \caption{Case $\lambda_2= -\nicefrac{5}{4}$, $\nicefrac{p}{q}=\nicefrac{1}{1}$, $c=0$ and $c=\pi$.}         
    \end{subfigure}
    \quad
    \begin{subfigure}[b]{0.3\textwidth}
    \centering
         \begin{tikzpicture}[>=latex,scale=1.45]
            \draw[->] (-.5,0) -- (2.5,0) node[right] {$\theta$};
            \foreach \x /\n in {1/$\pi$,2/$2\pi$} \draw[shift={(\x,0)}] (0pt,2pt) -- (0pt,-2pt) node[below] {\tiny \n};
            \draw[->] (0,-.5) -- (0,2.5) node[below right] {$\psi$};
            \foreach \y /\n in {1/$\pi$,2/$2\pi$}
            \draw[shift={(0,\y)}] (2pt,0pt) -- (-2pt,0pt) node[left] {\tiny \n};
            \node[below left] at (0,0) {\tiny $0$};
             \draw (0,0) rectangle (2,2);
             \draw (1,0)--(1,2);
             \draw (0,1)--(2,1);

              \draw[blue, thick] (0,0)--(8/5,2);
             \draw[blue, thick] (8/5,0)--(2,0.5);
             \draw[blue, thick] (0,0.5)--(6/5,2);
             \draw[blue, thick] (6/5,0)--(2,1);
             \draw[blue, thick] (0,1)--(4/5,2);
             \draw[blue, thick] (4/5,0)--(2,1.5);
             \draw[blue, thick] (0,1.5)--(2/5,2);
             \draw[blue, thick] (2/5,0)--(2,2);

             \draw[red, thick] (6/5+0.1,0)--(2,1-0.15);
             \draw[red, thick] (1/5+0.1,0)--(9/5+0.1,2);
            \draw[red, thick] (0,1-0.15)--(4/5+0.1,2);

             \draw [cyan, thick] (0,0)--(2,0);
              \draw [cyan, thick] (0,2)--(2,2);

            \draw[dotted, thick] (2,0) --(0,1);
            \draw[dotted, thick] (2,1) --(0,2);

            \draw[dashdotted, thick] (1,0)--(0,0.5);
            \draw[dashdotted, thick] (2,0.5)--(0,1.5);
            \draw[dashdotted, thick] (2,1.5)--(1,2);
               
        \end{tikzpicture}
       \caption{Case $\lambda_2= -\nicefrac{5}{4}$, $\nicefrac{p}{q}=\nicefrac{1}{2}$, and $c=\pi$.}
    \end{subfigure}
    \quad
    \begin{subfigure}[b]{0.3\textwidth}
    \centering
         \begin{tikzpicture}[>=latex,scale=1.45]
            \draw[->] (-.5,0) -- (2.5,0) node[right] {$\theta$};
            \foreach \x /\n in {1/$\pi$,2/$2\pi$} \draw[shift={(\x,0)}] (0pt,2pt) -- (0pt,-2pt) node[below] {\tiny \n};
            \draw[->] (0,-.5) -- (0,2.5) node[below right] {$\psi$};
            \foreach \y /\n in {1/$\pi$,2/$2\pi$}
            \draw[shift={(0,\y)}] (2pt,0pt) -- (-2pt,0pt) node[left] {\tiny \n};
            \node[below left] at (0,0) {\tiny $0$};
             \draw (0,0) rectangle (2,2);
             \draw (1,0)--(1,2);
             \draw (0,1)--(2,1);

             \draw[blue, thick] (0,0)--(8/5,2);
             \draw[blue, thick] (8/5,0)--(2,0.5);
             \draw[blue, thick] (0,0.5)--(6/5,2);
             \draw[blue, thick] (6/5,0)--(2,1);
             \draw[blue, thick] (0,1)--(4/5,2);
             \draw[blue, thick] (4/5,0)--(2,1.5);
             \draw[blue, thick] (0,1.5)--(2/5,2);
             \draw[blue, thick] (2/5,0)--(2,2);

             \draw[red, thick] (6/5+0.1,0)--(2,1-0.15);
             \draw[red, thick] (1/5+0.1,0)--(9/5+0.1,2);
            \draw[red, thick] (0,1-0.15)--(4/5+0.1,2);

             \draw [cyan, thick] (0,0)--(2,0);
              \draw [cyan, thick] (0,2)--(2,2);

            \draw[dotted, thick] (1,0) --(2,0.5);
            \draw[dotted, thick] (0,0.5) --(2,1.5);
            \draw[dotted, thick] (0,1.5) --(1,2);

        \end{tikzpicture}
       \caption{Case $\lambda_2= -\nicefrac{5}{4}$, $\nicefrac{p}{q}=-\nicefrac{1}{2}$, and $c=\pi$.}
    \end{subfigure}
    \caption{$A(Y_1)$ and $A(Y_2)$ are in cyan and blue. The LIPA-path of $Y_1$ is in red and the dotted and dash-dotted lines represent $L$ in the cases $c=0$ and $c=\pi$.}
    \label{Figure finales p2}
    \end{figure}
\printbibliography
\end{document}